%% file: BEPV_Wondertopes_revision.tex
\documentclass{amsart} 

\usepackage[paperheight=11in, 
    paperwidth=8.5in,
    outer=1.25in,
    inner=1.25in,
    bottom=1.4in,
    top=1.4in]{geometry} 

\usepackage{amsmath, amsxtra, amsthm, amsfonts, amssymb}
\usepackage{mathrsfs}
\usepackage{graphicx}
\usepackage{url}
\usepackage{color}
\usepackage{bbm}
\usepackage{tikz-cd}

\usepackage[T1]{fontenc}
\usepackage{enumitem}
\usepackage{standalone}

 \usepackage[hidelinks,backref=page]{hyperref} 
 \hypersetup{
    colorlinks,
    citecolor=magenta,
    filecolor=magenta,
    linkcolor=blue,
    urlcolor=black
}

\usepackage{cleveref}

\theoremstyle{definition}
\newtheorem{thm}{Theorem}[section]
\newtheorem{cor}[thm]{Corollary}
\newtheorem{lem}[thm]{Lemma}
\newtheorem{prop}[thm]{Proposition}
\newtheorem{defn}[thm]{Definition}

\newtheorem{eg}[thm]{Example}
\newtheorem{rem}[thm]{Remark}

\newtheorem{ques}[thm]{Question}

\newcommand{\RR}{\mathbb R}
\newcommand{\CC}{\mathbb C}
\newcommand{\PP}{\mathbb P}

\newcommand{\LL}{\mathscr L}
\renewcommand{\AA}{\mathbb A}

\newcommand{\TT}{\mathcal{T}}
\newcommand{\FF}{\mathcal{B}}
\newcommand{\NN}{\mathcal{N}}

\newcommand{\Res}{\operatorname{Res}}
\newcommand{\Bl}{\operatorname{Bl}}
\newcommand{\intr}{\operatorname{int}}

\newcommand{\build}{\mathcal B}
\newcommand{\poly}{\mathcal{P}}

\newcommand{\arrangement}{\mathcal{A}}
\newcommand{\BB}{\mathcal{B}}

\keywords{Positive geometries, De Concini--Procesi wonderful compactifications, moduli space of pointed stable rational curves, hyperplane arrangements, polytopes}

\title{Wondertopes}
\author{Sarah Brauner}
\email{sarah\_brauner@brown.edu}
\address{Division of Applied Mathematics, Brown University, Providence, RI, 02912}
\author{Christopher Eur}
\email{ceur@cmu.edu}
\address{Department of Mathematical Sciences, Carnegie Mellon University, Pittsburgh, PA, 15217}
\author{Elizabeth Pratt}
\email{epratt@berkeley.edu}
\address{Department of Mathematics, University of California, Berkeley, Berkeley, CA, 94720}
\author{Raluca Vlad}
\email{raluca\_vlad@brown.edu}
\address{Department of Mathematics, Brown University, Providence, RI, 02912}

\begin{document}

\maketitle

\begin{abstract}
Positive geometries were introduced by Arkani-Hamed--Bai--Lam in their study of scattering amplitudes in theoretical physics.  We show that a positive geometry from a polytope admits a log resolution of singularities to another positive geometry.  Our result states that the regions in a wonderful compactification of a hyperplane arrangement complement, which we call \emph{wondertopes}, are positive geometries.  A familiar wondertope is the curvy associahedron which tiles the moduli space of $n$-pointed stable rational curves $\overline{M}_{0,n}$.  Thus our work generalizes the known positive geometry $\overline{M}_{0,n}$.
\end{abstract}

\section{Introduction}

Positive geometries were introduced by Arkani-Hamed--Bai--Lam \cite{NimaLam17} as a mathematical framework that encodes integrands arising in the computation of scattering amplitudes. Various quantum field theories give rise to various positive geometries. A positive geometry is a pair consisting of a complex projective variety $X$ and a semialgebraic subset $X_{\geq 0}$ of its real points, together with a uniquely determined meromorphic form $\Omega(X,X_{\geq 0})$ satisfying the following technical definition.

\begin{defn}\label{defn:positive-geometry}
A \emph{positive geometry} is a pair $(X, X_{\geq 0})$ consisting of an $n$-dimensional irreducible complex normal projective variety $X$ defined over $\RR$, and a semialgebraic subset $X_{\geq 0} \subset X(\RR)$, along with a unique nonzero meromorphic $n$-form $\Omega(X,X_\geq 0)$, satisfying the following properties.
\begin{itemize}
\item The interior $X_{>0} := \intr(X_{\geq 0})$ of $X_{\geq 0}$ is a smooth oriented real $n$-manifold, and the Euclidean closure of $X_{>0}$ is $X_{\geq 0}$.
\item Let $\partial X$ be the Zariski closure of $X_{\geq 0} \setminus X_{> 0}$, whose codimension 1 irreducible components are $C_1, \dotsc, C_k$.
Then, with $C_{i, \geq 0}$ as the Euclidean closure of $\intr(C_i \cap (X_{\geq 0}\setminus X_{>0}))$ in $C_i(\RR)$, the $n$-form $\Omega(X,X_{\geq 0})$ satisfies the following recursive property:
\begin{enumerate}[label=(\alph*)]
    \item if $n=0$, then $X = X_{\geq 0}$ is a point and $\Omega(X, X_{\geq 0}) = \pm 1$, depending on orientation; and 
    \item if $n > 0$, the form $\Omega(X, X_{\geq 0})$ has simple poles along each $C_i$ and no other poles, and every $(C_i, C_{i, \geq 0})$ is a positive geometry, with the orientation on $C_{i, >0}$ inherited from that of $X_{>0}$, such that the residue of $\Omega(X, X_{\geq 0})$ along $C_i$ satisfies
\[
\Res_{C_i} \Omega(X, X_{\geq 0}) = \Omega(C_i, C_{i, \geq 0}).
\]
\end{enumerate}
\end{itemize}
The $n$-form $\Omega(X,X_{\geq 0})$ is called the \emph{canonical form} of $(X,X_{\geq 0})$, and $X_{\geq 0}$ the \emph{nonnegative part}.  We call $(C_i, C_{i, \geq 0})$ a (codimension 1) \emph{boundary component} of $(X, X_{\geq 0})$.
\end{defn}

Many well-studied spaces are positive geometries with naturally arising nonnegative regions.
Examples include flag varieties and toric varieties \cite{NimaLam17}, the moduli space $\overline{M}_{0,n}$ of $n$-pointed stable rational curves \cite{AHLcluster}, and del Pezzo surfaces and their moduli spaces \cite{delpezzo}.
Conjecturally, Grassmann polytopes and amplituhedra are nonnegative regions of positive geometries \cite[Conjecture 6.3]{NimaLam17}; see \cite{RST} for a recent progress on the amplituhedron conjecture.
 For a friendly introduction to positive geometries, see \cite{Lam22}.

\medskip
By definition, a positive geometry $(X, X_{\geq 0})$  defines a pair $(X,\partial X^1)$ of a variety and a divisor where $\partial X^1$ denotes the union of components of $\partial X$ of codimension 1 in $X$.
Thus, one may naturally ask:  (when) does a positive geometry $(X, X_{\geq 0})$ admit a version of log resolution of singularities?  More precisely, we ask:

\begin{ques}\label{ques:logres}
For a positive geometry $(X, X_{\geq 0})$, (when) is
there a positive geometry $(\widetilde X, \widetilde X_{\geq 0})$, along with a birational map $\pi: \widetilde X \to X$ restricting to a diffeomorphism $\widetilde X_{>0} \cong X_{>0}$, such that $\widetilde X$ is smooth and the boundary $\partial \widetilde X^1$ is a simple normal crossing divisor?
\end{ques}

Such a pair $(\widetilde X, \widetilde X_{\geq 0})$ may not always exist; see Example~\ref{eg:complex}.
However, we affirmatively answer the question for a prototypical example of a positive geometry, namely, a polytope in a projective space.
We achieve this by introducing a new family of positive geometries, whose nonnegative regions we call \emph{wondertopes}.  Let us now describe the construction.

\medskip
Let $V \cong \RR^{n+1}$ be a real vector space.
Let $\poly$ be a full-dimensional (convex) polytope in $\RR^n$, considered as a subset of the projective space $\PP V$ by the identification $\PP V \setminus \PP H \cong \RR^n$ for a hyperplane $H$ in $V$. 
We will abuse notation to write $\PP V$ also for the complex projective space $\PP (V \otimes_\RR \CC)$, as we trust that the field will be clear from context; as per the types in Definition~\ref{defn:positive-geometry}, the varieties considered throughout the paper are complex, while the semialgebraic sets live inside the real points of these varieties. The pair $(\PP V, \poly)$ is a positive geometry \cite{NimaLam17}, whose canonical form is well-studied for its connection to adjoint hypersurfaces \cite{KohnRanestad19} and the volume of dual polytopes~\cite[~\S~7.4]{NimaLam17}.

\smallskip
A wondertope is constructed from a polytope $\poly \subset \PP V$ and a set $\build$ of proper linear subvarieties of $\PP V$ satisfying the following properties:
\begin{enumerate}[label = (\arabic*)]
\item for every $F \in \build$, the intersection $F\cap \poly$ is a (possibly empty) face of $\poly$,
\item\label{item:hyperplane-condition} every hyperplane whose intersection with $\poly$ is a facet of $\poly$ is in $\build$, and
\item $\build$ is a \emph{building set} in the sense of \cite[Theorem 2.3.(2)]{de1995wonderful}; that is, for every intersection $L = \bigcap_{F\in \mathcal F} F$ of a subset $\mathcal F$ of $\build$, the subset
$$
\mathcal F_{L} = \{F\in \build \mid \text{$F$ is minimal (by inclusion) among elements of $\build$ that contain $L$}\}
$$
of $\build$ satisfies 
$
\displaystyle \sum_{F\in \mathcal F_{L}} \operatorname{codim}_{\PP V} F = \operatorname{codim}_{\PP V} L
$.
\end{enumerate}
See Definition \ref{def:buildingset_matroid} for an equivalent combinatorial definition of a building set in the context of hyperplane arrangements, and \cite[Definition 2.2]{pagaria2023hodge} in the context of subspace arrangements.

\begin{defn}\label{defn:wondertope}
With notation as above, choose any ordering of $\build = \{F_1, \dotsc, F_k\}$ such that $i\leq j$ if $F_i \subseteq F_j$, and define $X^\build$, also denoted $\PP V^\build$, to be the sequential blow-up
\[
X^{\FF} := \Bl_{F_k}(\dotsb (\Bl_{F_2}(\Bl_{F_1} \PP V))\dotsb )
\]
with the blow-down map $\pi_{\FF}: X^{\FF} \longrightarrow \PP V$.
Here, we repurpose notation to write $F_2$ for the strict transform of $F_2$ in $\operatorname{Bl}_{F_1}\PP V$, and similarly for $F_3, \dotsc, F_k$.
The \emph{wondertope} $\widetilde{\poly}^{\build}$ is defined as
\[
\widetilde{\poly}^{\build} = \text{the Euclidean closure of  } \pi_{\build}^{-1}(\text{the interior of }\poly) \text{ in $X^\build(\RR)$}.
\]
\end{defn}

That a wondertope is semialgebraic is verified in Corollary~\ref{cor:wondertopesemialg}.  Our main result is the following.

\begin{thm}\label{thm:intro_mainthm}
The pair $(X^\build, \widetilde \poly^\build)$ is a positive geometry whose canonical form is the pullback $\pi_\build^*\Omega(\PP V, \poly)$ of the canonical form of $(\PP V, \poly)$.
The boundary components of $(X^\build, \widetilde\poly^\build)$ are
(the strict transforms of)
the exceptional divisors $E_F$ for $F\in \build$ such that $\dim(\poly \cap F) = \dim F$.
\end{thm}

The variety $X^\build$ is known as the \emph{wonderful compactification} of the arrangement complement $C(\build) = \PP V \setminus (\bigcup \build)$, introduced by De Concini and Procesi in \cite{de1995wonderful}.
They showed that the boundary $X^\build \setminus C(\build)$, which consists of (the strict transforms of) the exceptional divisors of our sequential blow-up, is a simple normal crossing divisor.
Thus, Theorem~\ref{thm:intro_mainthm} implies that the pair $(X^\build, \widetilde \poly^\build)$ is a log resolution of singularities of the positive geometry $(\PP V, \poly)$ in the sense of Question~\ref{ques:logres}.

\smallskip
For $(X^\build, \widetilde\poly^\build)$ to be a positive geometry, but not necessarily with simple normal crossing boundary components, the condition on $\build$ can be relaxed slightly; see Theorem~\ref{thm:maingeneral}.
We caution however that removing any one of the three conditions on $\build$ above results in a failure of the conclusion of Theorem~\ref{thm:intro_mainthm}; see Examples~\ref{eg:triangle}, \ref{eg:pyramid2}, and \ref{eg:pyramid3}.

\medskip
One notable special case of Theorem \ref{thm:intro_mainthm} arises from the \emph{rank $(n-1)$ braid arrangement} $A_{n-1}$ comprised of hyperplanes $\{ x_i - x_j = 0 \}$ in $V = \RR^n / \RR\cdot (1, 1, \dotsc, 1)$ for $1 \leq i < j \leq n$.
The complement $C(A_{n-1}) = \PP V \setminus A_{n-1}$ is isomorphic to $M_{0,n+1}$, the moduli space of $n+1$ distinct points on $\PP^1$. For an appropriate choice of building set, the corresponding wonderful compactification is isomorphic to the Deligne--Knudsen--Mumford compactification $\overline{M}_{0,n+1}$ \cite{deligne1969irreducibility}. The polytope in this context is a simplex, and the corresponding wondertope $(\overline{M}_{0,n+1})_{\geq 0}$ is combinatorially isomorphic to the $n$-associahedron, which parameterizes triangulations of an $(n+1)$-gon---hence, it is sometimes called the ``curvy associahedron,'' and was first studied by Devadoss \cite{devadoss239tessellations}. Note that this wondertope $(\overline{M}_{0,n+1})_{\geq 0}$ is the \emph{tropical compactification} of $M_{0, n+1}$ in the sense of \cite{arkani2021stringy}. Theorem \ref{thm:intro_mainthm} then recovers the following result of Arkani-Hamed--He--Lam \cite[Proposition 8.2]{AHLcluster}, which uses work of Brown--Carr--Schneps \cite[Proposition 2.7]{brown2010algebra}.
See Section \ref{section:braid} for details.

\begin{cor}\label{cor:m0n-positivegeometry}
The pair $(\overline{M}_{0,n+1}, (\overline{M}_{0,n+1})_{\geq 0})$ is a positive geometry with canonical form given by the \emph{Parke-Taylor} form.
\end{cor}

\medskip
To any hyperplane arrangement, one can associate a matroid.
As a real hyperplane arrangement defines an oriented matroid, and a region in the arrangement can be understood as a maximal covector of that matroid, one may ask the following.

\begin{ques}\label{question:matroidalwondertope}
Is there a purely matroidal notion of wondertopes and their canonical forms?
\end{ques}

In Appendix~\ref{section:combinatorialproduct}, with a view towards future studies, we summarize the combinatorial counterparts in matroid theory of the geometry of wonderful compactifications.
For instance, our proof of Theorem \ref{thm:intro_mainthm} uses a key property that exceptional divisors in a wonderful compactification are naturally a product of two smaller wonderful compactifications \cite{de1995wonderful} (see Proposition~\ref{prop:wonderful}).
The following theorem is the combinatorial analogue for matroids; see Appendix~\ref{section:combinatorialproduct} for relevant terminologies and definitions.

%To any hyperplane arrangement, one can associate a matroid $M$ whose lattice of flats $\LL(M)$ is isomorphic to the intersection lattice of the arrangement. Then, for every building set $\build \subset \LL(M)$, the intersection lattice of the boundary divisors in the wonderful compactification $X^\build$ is encoded by the \emph{nested set complex} $\NN(\BB)$, a simplicial complex introduced by DeConcini--Procesi \cite{de1995wonderful} and generalized to the context of matroids by Feichtner--Kozlov \cite{feichtner2004incidence}. For every $F\in \build$, the \emph{link} $\mathcal{N}(\build)_F$ of the nested set complex $\mathcal{N}(\build)$ encodes the intersection of boundary divisors restricted to the exceptional divisor corresponding to $F$.

\begin{thm}\label{thm:nestedproduct}
Let $M$ be a matroid and $\build \subset \LL(M)$ a building set.
For any $F \in \FF$, 
we have an isomorphism
\[ \mathcal{N}(\build)_F \cong \mathcal{N}(\FF^{F}) \times \mathcal{N}(\FF_F),\]
where $\FF^{F}$ and $\FF_F$ are certain building sets defined on the restriction and contraction matroids $M^{F}$ and $M_F$, respectively. 
\end{thm}

Theorem \ref{thm:nestedproduct} generalizes a result of Zelevinsky for Boolean matroids \cite{zelevinsky2006nested}, but is not new:
a more general version of Theorem~\ref{thm:nestedproduct} can be found in \cite[Propositions 2.8.6-7]{bibby2022leray}. However, it takes some work to specialize their results to our context, so we provide a proof for completeness.

\subsection*{Related works}
In the study of bordifications of links of tropical moduli spaces \cite{Brown23}, wondertopes $\widetilde\poly^\build$ of a polytope $\poly$ with respect to an intersection-closed set $\build$ of linear subspaces appeared as ``blow-ups of polyhedral linear configurations.''
In a slightly different setting of stratified manifolds and real blow-ups\footnote{Loosely put, our notion of a \emph{blow-up} in the sense used in algebraic geometry \cite[Chapter II.7]{hartshorne1977algebraic} replaces the blown-up locus with the projectivization of its normal bundle, whereas the ``real blow-up'' known as the ``balls, beams, and plates construction'' in \cite{KT} replaces the blown-up locus with the sphere-bundle of its normal bundle.} \cite{Gai03}, wondertopes essentially appeared as compactifications of interiors of polytopes.
In a related direction, ``poset associahedra'' were introduced in \cite{Gal} as certain compactifications of the interior of Stanley's order polytope of a poset.  Lastly, while this paper was in preparation, we were informed of the forthcoming work \cite{acyclohedra} on acyclonestohedra.  It introduces the notion of oriented building sets on oriented matroids, which partially answers Question~\ref{question:matroidalwondertope}.
Since this paper was posted, a notion of canonical forms of oriented matroids was given in \cite{EL25}, which is another step towards answering Question~\ref{question:matroidalwondertope}. 
Brown-Dupont \cite{Brown-Dupont} propose an interpretation of positive geometries using Deligne's mixed Hodge theory for complex varieties; in their setting, invariance under blow-ups is immediate, whereas our paper shows that the situation is more subtle in the recursive setting (see non-examples in Sections \ref{section:nonexamples} and \ref{section:detailed-example}).

\subsection*{Organization}
The remainder of the paper proceeds as follows. Section \ref{section:semialgebraicsets} discusses properties of semialgebraic subsets, and proves that wondertopes are indeed semialgebraic. In Section \ref{section:blowingupsingleface} we show that blowing-up a single face in a polytope results in a positive geometry. We then pair this result with induction in Section \ref{section:blow-up-sequence} to prove Theorem \ref{thm:intro_mainthm}. Section \ref{section:examples} focuses on examples; we discuss $\overline{M}_{0,n}$ and Corollary \ref{cor:m0n-positivegeometry} in Section \ref{section:braid}, and feature several pathologies in Sections \ref{section:nonexamples} and \ref{section:detailed-example} that illustrate the necessity of the conditions imposed on $\build$.  Appendix \ref{section:combinatorialproduct} rephrases the structure of wonderful compactifications in the language of matroids, and proves Theorem \ref{thm:nestedproduct}.

\subsection*{Acknowledgements}
The authors are grateful to Bernd Sturmfels for suggesting this problem, and to Spencer Backman, Melody Chan, Nick Early, Thomas Lam, Chiara Meroni, Arnau Padrol, Roberto Pagaria, Daniel Plaumann, Frank Sottile, Simon Telen, Josephine Yu for useful discussions.
The first author is supported by US NSF DMS-2303060. The second author is supported by US NSF DMS-2246518. The third author is supported by an NSF Graduate Research Fellowship. The authors are also grateful to the two anonymous referees for their valuable suggestions.

\section{Semialgebraic subsets}\label{section:semialgebraicsets}

We define and record some properties of semialgebraic subsets, and show that wondertopes are semialgebraic.
Much of this section can be found in a standard reference in the general setting of \emph{real spectra}, which is beyond our needs.
For the convenience of the reader, we include here a specialized and self-contained account of semialgebraicity in a scheme over $\RR$, and verify its compatibility with the notion of semialgebraicity in a projective space given in \cite{NimaLam17}.
Readers willing to believe the semialgebraicity of a wondertope (Corollary~\ref{cor:wondertopesemialg}) may skip this section.

\begin{defn}
A polynomial $f\in \RR[x_1, \dotsc, x_n]$ defines a subset $\{x\in \RR^n: f(x) > 0\}$ of $\RR^n$.
The collection of \emph{finite boolean combinations} of such subsets, i.e.\ the smallest collection that contains all such subsets and is closed under finite unions, finite intersections, and complements, is called the set of \emph{semialgebraic subsets} of $\RR^n$.
\end{defn}

For example, a polytope and the interior of a polytope in $\RR^n$ are semialgebraic.
Semialgebraic subsets are known to satisfy the following properties, both of which follow from the Tarski--Seidenberg theorem.

\begin{prop}\label{prop:basicsemialg}
\cite[Chapter 2]{BochnakCosteRoy}
\begin{itemize}
\item Images and preimages of semialgebraic subsets under an algebraic map are semialgebraic.
\item The closure and the interior (with respect to the Euclidean topology on $\RR^n$) of a semialgebraic subset are semialgebraic.
\end{itemize}
\end{prop}

We now define semialgebraic subsets of schemes.
Let $A$ be a finitely generated $\RR$-algebra, and $Y = \operatorname{Spec} A$ its affine scheme.
Because the only automorphism of the field $\RR$ is the identity, an $\RR$-valued point $y \in Y(\RR)$, i.e.\ a ring map $y: A \to \RR$, is uniquely determined by the maximal ideal $\mathfrak m_y = y^{-1}(0)$ of $A$.  Thus, for the point $y$ in $\operatorname{Spec} A$ and an element $f\in A$, the value $f(y)$ is well-defined as the image of $f$ under $A/\mathfrak m_y \cong \RR$.
We may hence define the following.

\begin{defn}
For the affine scheme $Y = \operatorname{Spec} A$ of a finitely generated $\RR$-algebra $A$, a subset $S\subseteq Y(\RR)$ of its $\RR$-valued points is \emph{semialgebraic} if it is a finite boolean combination of subsets of the form $\{y \in Y(\RR) \mid f(y) >0\}$ for some $f\in A$.

For a scheme $X$ of finite type over $\RR$, a subset $S\subseteq X(\RR)$ is \emph{semialgebraic} if $S\cap U$ is semialgebraic for every affine open subset $U\subset X$.  By convention, we assume a scheme to be separated.
\end{defn}

We observe two features of this definition:
\begin{itemize}
\item The definition of semialgebraic subsets in an affine scheme is compatible with the one for subsets of $\RR^n$.  More precisely, for the affine scheme $Y = \operatorname{Spec} A$, choosing a generating set $x_1, \dotsc, x_n$ for the $\RR$-algebra $A$ gives a closed embedding $\iota: Y \hookrightarrow \AA^n(\RR)$.  Then, a subset $S\subseteq Y(\RR)$ is semialgebraic if and only if $\iota(S) \subseteq \RR^n$ is semialgebraic, since every $f\in A$ has a lift $\widetilde f \in \RR[x_1, \dotsc, x_n]$ such that any two such lifts restrict to the same function on $\iota(Y)(\RR)$.
\item For a scheme $X$ of finite type over $\RR$, it suffices to check the semialgebraicity condition for an open affine cover of $X$, since open subsets and finite unions of semialgebraic subsets are semialgebraic, and since finite type schemes and intersections of affine open subsets of (separated) schemes are quasi-compact.
\end{itemize}
Together, the two observations imply the following.

\begin{lem}
Let $\{(U_\alpha, \iota_\alpha)\}$ be an affine open cover of a scheme $X$ of finite type over $\RR$, with closed embeddings $\iota_\alpha: U_\alpha \hookrightarrow \RR^{n_\alpha}$.  Then, a subset $S\subseteq X(\RR)$ is semialgebraic if $\iota_\alpha(U_\alpha \cap S)$ is semialgebraic for all $\alpha$.
\end{lem}

Proposition~\ref{prop:basicsemialg} generalizes straightforwardly to the setting of semialgebraic subsets of a scheme, as follows; we omit the proof.

\begin{prop}\label{prop:basicsemialg2}
Let $X$ and $Y$ be schemes of finite type over $\RR$, and $\varphi: X\to Y$ a morphism (which is necessarily quasi-compact).
\begin{itemize}
\item Images and preimages of semialgebraic subsets are semialgebraic under $\varphi$.
\item The closure and the interior (with respect to the Euclidean topology) of a semialgebraic subset are semialgebraic.
\end{itemize}
\end{prop}

\begin{cor}\label{cor:wondertopesemialg}
For any sequential blow-up $\pi: X^\build \to \PP V$ of linear subspaces that do not intersect the interior of a polytope $\poly \subset \PP V$, the strict transform $\widetilde\poly^\build$, i.e. the Euclidean closure of $\pi^{-1}(\intr(\poly))$ in $X^\build(\RR)$, is semialgebraic.
\end{cor}

\begin{rem}
The authors of \cite{NimaLam17} define a semialgebraic subset of a projective space as follows:
A subset $S\subset \PP^n(\RR)$ is said to be semialgebraic if it is the image under $\RR^{n+1}\setminus \{0\} \to \PP^n(\RR)$ of a semialgebraic subset $S\subseteq \RR^{n+1}$, all of whose defining equalities and inequalities are given by homogeneous polynomials.

Let us verify that their definition agrees with ours.
Proposition~\ref{prop:basicsemialg2} implies that if $S\subset \PP^n(\RR)$ is semialgebraic in their sense, then it is also in our sense.
For the converse, since $\PP^n$ is covered by the (finitely many) coordinate affine charts, it suffices to show the converse for a subset $S$ of $\PP^n$ that avoids a coordinate hyperplane, say $\{x_0 = 0\}$.  In that case, the converse holds because under $\AA^n \cong U_0 = \PP^n\setminus \{x_0 = 0\}$, we have 
\[
\{ x \in \AA^n(\RR) : f(x) > 0 \} = \text{the image under $\RR^{n+1}\setminus \{x_0 = 0\} \to U_0$ of } \{(x_0, x) \in \RR^{n+1} : F(x) > 0 \}
\]
for any polynomial $f$ and any even degree homogenization $F$ of $f$ by the variable $x_0$.
\end{rem}

\begin{rem}
    More concretely, for a building set $\build$ as in the Introduction, consider the rational map
    $$\varphi : \PP V \dashrightarrow \PP V \times \prod_{\PP W \in \build} \PP(V/W),$$
    coming from the identity $\PP V \to \PP V$ on the first coordinate and the quotient $V\to V/W$ on the other coordinates. According to \cite{de1995wonderful}, the sequential blow-up $X^\build$ is isomorphic to the Zariski closure of the image of $\varphi$ inside $\PP V \times \prod_{\PP W \in \build} \PP(V/W)$, with the blow-down map $\pi$ being projection onto the first factor $\PP V$. Under this embedding, the interior of the wondertope $\widetilde\poly^\build$ is cut out by the multi-homogeneous equations cutting out $X^\build$, together with the strict linear inequalities in the coordinates of the first factor $\PP V$ corresponding to the facet hyperplanes of $\poly$. 
    Finding the multi-homogeneous equations and inequalitites that cut out the entire wondertope $\widetilde\poly^\build$ inside this product $\PP V \times \prod_{\PP W \in \build} \PP(V/W)$
    would be an interesting future direction of study.
\end{rem}

\section{Blowing up a single face}\label{section:blowingupsingleface}

As in the introduction, let $V$ be an $(n+1)$-dimensional real vector space, and let $\poly$ be a full-dimensional polytope in $\PP V$, which defines a positive geometry $(\PP V, \poly)$.  The goal of this section is to prove the following ``fundamental computation,'' from which the proof of Theorem~\ref{thm:intro_mainthm} will follow via induction in Section~\ref{section:blow-up-sequence}.

\begin{thm}\label{thm:fundcomp}
Let $W\subset V$ be a proper subspace such that $\PP W \cap \poly$ is a (possibly empty) face of $\poly$.
Let $\pi: X = \operatorname{Bl}_{\PP W} \PP V \to \PP V$ be the blow-up, and let the \emph{strict transform} $\widetilde\poly$ of $\poly$ be
\[
\widetilde \poly = \text{the Euclidean closure of }\pi^{-1}(\text{the interior of } \poly)
\text{ in } X(\RR).\]
Then, the pair $(X, \widetilde \poly)$ is a positive geometry whose canonical form is $\pi^*\Omega(\PP V, \poly)$.
\end{thm}

In Section~\ref{section:polytopes-triangulations}, we prepare by collecting known results about polytopes, normal cones, and triangulations.  Then, in Section~\ref{subsection:fundamental-computation}, we use these to prove Theorem~\ref{thm:fundcomp}.

\subsection{Polytopes, normal cones, and triangulations} \label{section:polytopes-triangulations}
A \emph{polytope} is the convex hull of a finite set of points in $\RR^n$.
Throughout, let $\poly$ be a \emph{(projective) polytope} in $\PP V$, by which we mean a polytope in $\PP V \setminus \PP H \cong \RR^n$ for some hyperplane $H$.
%By the definition of a polytope in $\PP V$, there is a hyperplane $H \subset V$ such that $\PP H \cap \poly$ is empty.
Choosing $H^+$ to be one of the two open half-spaces in $V$ defined by $H$, we obtain the following polyhedral cone in $V$:
\begin{align*}
\widehat\poly &= \RR_{\geq 0}\{ v\in H^+ \mid \text{image of $v$ in $\PP V$ lies in $\poly$}\}\\
&= \RR_{\geq 0}\{ v\in H^+ \mid \text{image of $v$ in $\PP V$ is a vertex of $\poly \subset (\PP V \setminus \PP H) \cong \AA^n(\RR)$}\},
\end{align*}
which is \emph{pointed} (or also called \emph{strongly convex}), i.e.\ $\widehat\poly \cap (-\widehat\poly) = \{0\}$.
Its image in $\PP V$ recovers the original polytope $\poly$.
If $\poly = \emptyset$, we set $\widehat P = \{0\}$.
We will freely switch between considering a projective polytope as a polytope in $\PP V \setminus \PP H \cong \RR^n$ for some hyperplane $H$ and as the image in $\PP V$ of a pointed polyhedral cone in $V$. In other words, we will freely switch between a projective polytope $\poly \subset \PP V$ and its corresponding pointed polyhedral cone $\widehat{\poly} \subset V$.

\begin{defn}
For a linear subspace $W\subset V$ such that $\widehat\poly \cap W$ is a face of $\widehat\poly$, we define two pointed polyhedral cones $\widehat{\poly_W} \subset W$ and $\widehat{\poly^W} \subset V/W$ by
\[
\widehat{\poly_W} = \widehat\poly \cap W \qquad\text{and}\qquad \widehat{\poly^W} = 
\begin{cases}
\{0\} & \text{if $\widehat\poly \cap W = \{0\}$}\\
\text{the image of $\widehat\poly$ under $V\to V/W$} & \text{otherwise}.
\end{cases}
\] 
We call the polytope $\poly_W \subset \PP W$ the \emph{face of $\poly$ relative to $W$}, and call the polytope $\poly^W \subset \PP(V/W)$ the \emph{normal polytope of $\poly$ relative to $W$}.
\end{defn}

The cone $\widehat{\poly^W}$ is sometimes known as the \emph{normal cone} of $\widehat\poly$ relative to $W$.
It is a pointed polyhedral cone 
(see for instance \cite[Lemma~4.7]{Brown23}), so that the normal polytope $\poly^W$ is a well-defined projective polytope. We caution that our definition of normal polytope is not the same as the polar dual polytope constructed from the normal fan.

We observe that these induced polytopes $\poly_W$ and $\poly^W$ arise in our context via the following Proposition; for a proof, see \cite[Proposition 5.10]{Brown23}.\footnote{\cite[Proposition 5.10]{Brown23} is stated for the more general case of a sequential blow-up of an intersection-closed collection of linear subspaces, which specializes easily to the case stated here.}
For a description of blow-ups and exceptional divisors appearing below, see \cite[II.7]{hartshorne1977algebraic}.

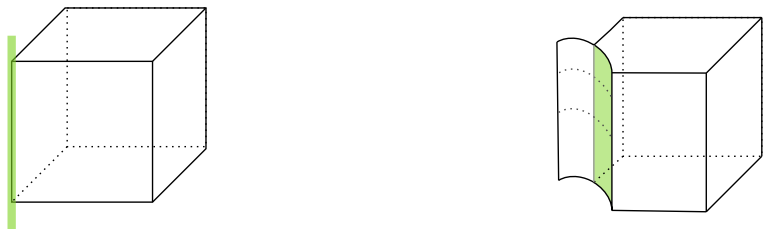
\begin{figure}[!h]
	\begin{minipage}{0.45\textwidth}
 \centering
    \scalebox{.7}{\input{./BEPV_Figures_arXiv/prodbefore.tex}}
    \end{minipage}
	\begin{minipage}{0.45\textwidth}
	\centering
 \scalebox{.7}
    {\input{./BEPV_Figures_arXiv/prodafter.tex}}
    \end{minipage}
    \caption{Left: a cube in $\PP^3$ and a $\PP^1$ (green) intersecting the cube along an edge $F$. Right: the cube after blowing-up $\PP^3$ along the $\PP^1$, with $E \cong \PP^1 \times \PP^1$ and $E_{\geq 0} \cong F \times \Delta^1$ (the latter in green). Here, $\Delta^1$ denotes the standard simplex of the projective line (see Proposition~\ref{prop:canonical-form-polytope}~\ref{prop:canonical-form-general-simplex} for a definition).}
    \label{fig:exceptional-product}
\end{figure}

\begin{prop}\label{prop:exceptionalproduct}
Let $\pi: \operatorname{Bl}_{\PP W} \PP V \to \PP V$ be the blow-up, and $E$ the exceptional divisor.
Suppose $\poly$ is full-dimensional.
Then, under the canonical isomorphism
$E \cong \PP W \times \PP(V/W)$,
the strict transform $\widetilde\poly$ satisfies $E\cap \widetilde\poly \cong \poly_W \times \poly^W$. See Figure~\ref{fig:exceptional-product}.
\end{prop}

We now describe how faces and normal polytopes interact with  triangulations of polytopes.

\begin{defn}\label{def:triangulation}
A \emph{triangulation} of $\poly$ is a finite collection $\TT$ of simplices in $\PP V$ of dimensions all equal to $\dim \poly$ such that
\begin{enumerate}[label = (\roman*)]
        \item[(i)] the simplices in $\TT$ cover $\poly$, i.e. $\bigcup_{T \in \TT} T = \poly$, and 
        \item[(ii)] any two simplices $T, T' \in \TT$ intersect in a (possibly empty) common face.
\end{enumerate}
\end{defn}

Let us fix a triangulation $\TT$ of $\poly$, and fix a linear subspace $W\subset V$ such that $\PP W$ intersects $\poly$ along a (possibly empty) face.  Note that $\PP W$ then intersects each simplex $T \in \TT$ along a (possibly empty) face, so we may consider $T_W$ and $T^W$.

\begin{lem}\label{lem:induced-triangulations}
One has triangulations of $\poly_W$ and $\poly^W$ obtained from $\TT$ as follows.
\begin{enumerate}[label = (\alph*)]
\item \label{lem:triangulation-of-face} The collection $\TT_W = \{T_W \mid T\in \TT \text{ satisfying } \dim T_W = \dim \poly_W\}$ is a triangulation of $\poly_W$.
\item \label{lem:triangulation-of-normal} Suppose further that $\dim \poly_W = \dim \PP W$.  Then, for any simplex $F$ in the triangulation $\TT_W$ of $\poly_W$, the collection $\TT^{W,F} = \{T^W \mid T\in \TT \text{ satisfying } T_W = F\}$ is a triangulation of $\poly^W$.
\end{enumerate}
\end{lem}

\begin{proof}
We omit the proof of the first statement \ref{lem:triangulation-of-face}, as it is straightforward.
We prove the second statement \ref{lem:triangulation-of-normal}, which is false without the assumption $\dim \poly_W = \dim \PP W$.

Every $T\in \TT^{W,F}$ is a full dimensional simplex in $\PP V$ which intersects $\PP W$ along $F$, a full-dimensional simplex in $\PP W$. Hence the normal face $T^W$ of $T$ relative to $\PP W$ is a full-dimensional simplex in $\PP(V/W)$.
With the cone $\widehat\poly$ obtained by a choice of an open half-space $H^+ \subset V$, we let $\widehat T$ be the cone of a simplex $T\in \mathcal T$ for the same half-space $H^+$.
We need show that $\{\widehat{T^W} \mid T^W\in \mathcal T^{W,F}\}$ is a triangulation of $\widehat{\poly^W}$.

Pick any non-zero $v \in \widehat{\poly} \subset V$ with image $\overline{v} \in \widehat{\poly^W} \subset V/W$. We will show that $\overline{v}$  belongs to $\widehat{T^W}$ for some $T^W \in \TT^{W,F}$. Consider a ray $w \in \widehat{\poly}$ whose projectivization lies in the interior of $F$. Then for any positive integer $c$, both $\frac{1}{c} \cdot v$ and $w$ lie in $ \widehat{\poly}$, so the sum $\frac{1}{c} \cdot v + w$ does as well. Since our original triangulation $\TT$ covers $\poly$, there must exist a $T \in \TT$ whose cone  $\widehat{T}$ contains $\frac{1}{c} \cdot v + w$ for infinitely many positive integers $c$.  This cone $\widehat{T} \subset V$ is closed in the real topology induced by any linear identification $V \cong \RR^{n+1}$, so $\widehat{T}$ must also contain the limit point $w$. By assumption, $w$ lies in the interior of $F$, so $T^W \in \TT^{W,F}$ by part \ref{lem:triangulation-of-face}. Its normal cone $\widehat{T^W}$ thus contains $\overline{v + c \cdot w} = \overline{v}$.
    
Finally, consider any simplices $T_1,T_2 \in \TT^{W,F}$. Since $\TT$ is a triangulation, $T_1$ intersects $T_2$ in a common face $S = T_1 \cap T_2$. We claim that $T_1^W \cap T_2^W = S^W$, which is a common face. By construction, both $T_1$ and $T_2$ intersect $\PP W$ along $F,$ which is a full-dimensional simplex in $\PP W$. Thus the same must be true for $S = T_1 \cap T_2$. It follows that $\widehat{S^W} \subseteq \widehat{T_1^W}\cap \widehat{T_2^W}$ is a common face. 
In fact, this containment is an equality since, by assumption, the cones $\widehat{T_1}$ and $\widehat{T_2}$ intersect only along $\widehat{S}$.
\end{proof}

In our context, the utility of considering triangulations arises from the following known results in positive geometry.

\begin{prop}\label{prop:canonical-form-polytope}
For a full-dimensional polytope $\poly \subset \PP V$, choose an orientation of its interior by an orientation of an affine chart $\AA^n(\RR) \cong \RR^n$ containing it.
\begin{enumerate}[label = (\alph*)]
\item \cite[Section 6.1]{NimaLam17} The pair $(\PP V, 
\poly)$ is a positive geometry.   
\item \cite[Section 3.1]{NimaLam17} \label{prop:canonical-form-triangulation} If $\TT$ is a triangulation of $\poly$, all of whose simplicies are oriented the same way as $\poly$, then the canonical forms of $(\PP V,\poly)$ and $(\PP V,T)$ for $T \in \TT$ satisfy
    $$\Omega(\PP V, \poly) = \sum_{T \in \TT} \Omega(\PP V, T).$$

\item \cite[Section 5.1]{NimaLam17}  \label{prop:canonical-form-general-simplex} Suppose $\poly$ is the standard simplex
$
\Delta^n = \{[X_0: \dotsc : X_n] \mid X_i \geq 0 \text{ for all } i\} \subset \PP^n(\RR).
$
Let $U$ be the affine chart $\AA^n \cong U = \{X_0 \neq 0\} \subset \PP^n$ with coordinates $x_1 = \frac{X_1}{X_0}, \dots, x_n = \frac{X_n}{X_0}$, in which the interior of $\Delta^n$ is identified with the positive orthant $\{x \in \AA^n(\RR) \mid x_i > 0 \text{ for all $i$}\}$.
If the standard orientation of $\AA^n(\RR) \cong U(\RR)$ given by the ordered coordinates $(x_1, \dotsc, x_n)$ agrees with the orientation of $\intr(\Delta^n)$, then the canonical form satisfies
$$\Omega(\PP^n, \Delta^n)|_{U} = \frac{dx_1}{x_1} \wedge \dotsm \wedge \frac{dx_n}{x_n}.$$

\end{enumerate}
\end{prop}

Lastly, we record the orientation and residue conventions from \cite[Appendix A]{NimaLam17} in the following remark, with a focus on  the case of the standard simplex $\Delta^n$ from Proposition~\ref{prop:canonical-form-polytope}\ref{prop:canonical-form-general-simplex}.

\begin{rem}\label{rem:orientation}
Given a boundary component $(C,C_{\geq0})$ of $(X,X_{\geq0})$, 
pick a (boundary) chart $U\subset X_{\geq0}$ diffeomorphic to $\mathbb{R}^{n-1}\times\mathbb{R}_{\geq0}$ such that the restriction $U\cap C_{>0}$ is identified with $\mathbb{R}^{n-1}\times\{0\}$. 
For instance, one can pick $U$ to be an analytic neighborhood around a point of $C_{>0}$.
We may choose the identification $U \cong \RR^{n-1} \times \RR_{\geq 0}$ such that the restricted identification between $U\cap X_{>0}$ and $\mathbb{R}^{n-1}\times\mathbb{R}_{>0}$ is orientation preserving, where $\mathbb{R}^{n-1}\times\mathbb{R}_{>0} \subset \mathbb{R}^n$ has the standard orientation. 
We orient $C_{>0}$ such that the identification between $U\cap C_{>0}$ and $\mathbb{R}^{n-1}\times\{0\}$ is orientation preserving, where $\mathbb{R}^{n-1}\times\{0\}$ has the standard orientation. We say that this is the orientation that $X_{>0}$ \emph{induces~on}~$C_{>0}$.

For the residue, if a 
divisor is locally given by an equation $f = 0$, then the residue of a form with a simple pole at $\{f = 0\}$, when locally written as $\alpha \wedge \frac{df}{f}$, is $\alpha|_{f = 0}$.
Let us explicitly describe these conventions in the case of the standard simplex $\Delta^n \subset \PP^n(\RR)$.
    
    Consider the boundary component $(\{X_n = 0\}, \Delta^n \cap \{X_n = 0\}) \cong ( \PP^{n-1}, \Delta^{n-1})$. It intersects our chart $\AA^n$ along $\AA^{n-1} \cong \{x_n = 0\} \subset \AA^n$. Since the standard Euclidean orientation on $\AA^n(\RR)$ induces the standard Euclidean orientation on the hyperplane $\AA^{n-1}(\RR) \cong \{x_n = 0\}$, then we consider $\Delta^{n-1}$ oriented by the standard orientation on the real points of our chart $\AA^{n-1}\cong \{X_0 \neq 0\} \subset \PP^{n-1}$. And indeed, the recursive formula in Definition~\ref{defn:positive-geometry} holds:
    $$\Res_{\{x_n = 0\}} \; \frac{dx_1}{x_1} \wedge \dotsm \wedge \frac{dx_n}{x_n} 
        = \frac{dx_1}{x_1} \wedge \dotsm \wedge \frac{dx_{n-1}}{x_{n-1}}.$$

    However, consider instead $(\{X_{n-1} = 0\}, \Delta^n \cap \{X_{n-1} = 0\}) \cong ( \PP^{n-1}, \Delta^{n-1})$. 
    Because the standard Euclidean orientation on $\AA^n(\RR)$ induces the opposite orientation on $\AA^{n-1}(\RR) \cong \{x_{n-1} =~0\}$, 
    we consider $\Delta^{n-1}$ oriented by the opposite orientation on our chart $\AA^{n-1}\cong \{X_0 \neq 0\} \subset \PP^{n-1}$. The recursive formula in Definition~\ref{defn:positive-geometry} catches this opposite orientation:
    \begin{align*}
        \Res_{\{x_{n-1} = 0\}} \frac{dx_1}{x_1} \wedge \dotsm \wedge \frac{dx_n}{x_n} & = \Res_{\{x_{n-1} = 0\}}  \left(- \frac{dx_1}{x_1} \wedge \dotsm \wedge \frac{dx_{n-2}}{x_{n-2}} \wedge \frac{dx_n}{x_n}\right) \wedge \frac{dx_{n-1}}{x_{n-1}}\\
        &= - \frac{dx_1}{x_1} \wedge \dotsm \wedge \frac{dx_{n-2}}{x_{n-2}} \wedge \frac{dx_n}{x_n}.
    \end{align*}
\end{rem}

\subsection{The fundamental computation}\label{subsection:fundamental-computation}

Let us recall the setting of Theorem~\ref{thm:fundcomp}. That is, let
\begin{itemize}
\item $\poly$ be a full-dimensional (oriented) polytope in $\PP V$,
\item $W\subset V$ be a proper linear subspace such that $\PP W \cap \poly$ is a (possibly empty) face of $\poly$,
\item $\pi: X = \operatorname{Bl}_{\PP W} \PP V \to \PP V$ be the blow-up, with exceptional divisor $E$, and
\item $\widetilde\poly \subset X(\RR)$ be the strict transform of $\poly$, i.e.\ the Euclidean closure of $\pi^{-1}(\intr(\poly))$.
\end{itemize}
Since $\pi$ restricts to an isomorphism $X\setminus E \overset\sim\to \PP V \setminus \PP W$, in particular a diffeomorphism $\pi^{-1}(\intr(\poly)) \to \intr(\poly)$, the interior $\intr(\widetilde\poly)$ of $\widetilde\poly$ is a real $n$-manifold, which we orient such that $\pi$ is orientation preserving.
Let $\partial \widetilde\poly$ be the Zariski closure of $\widetilde\poly \setminus \intr(\widetilde\poly)$ inside $X$.
The following two propositions constitute the key computations for the proof of Theorem~\ref{thm:fundcomp}.

\begin{prop}\label{prop:Elessdim}
The exceptional divisor $E$ is an irreducible component of $\partial \widetilde\poly$ if and only if $\dim \poly_W = \dim \PP W$.  When $\dim \poly_W < \dim \PP W$, the form $\pi^*\Omega(\PP V, \poly)$ has no pole along $E$.
\end{prop}

Let $E_{\geq 0} = E \cap \widetilde\poly$.  Under $E \cong \PP W \times \PP(V/W)$, recall that $E_{\geq 0} \cong \poly_W \times \poly^W$ (Proposition~\ref{prop:exceptionalproduct}).

\begin{prop}\label{prop:Eequaldim}
When $\dim \poly_W = \dim \PP W$, the form $\pi^*\Omega(\PP V, \poly)$ has a simple pole along $E$, and under the isomorphism $E \cong \PP V \times \PP(V/W)$, its residue along $E$ satisfies
\[
\operatorname{Res}_E(\pi^*\Omega(\PP V, \poly)) = \Omega(\PP W, \poly_W) \wedge \Omega(\PP(V/W), \poly^W).
\]
Moreover, the pair $(E, E_{\geq 0})$ is a positive geometry with the canonical form $\operatorname{Res}_E(\pi^*\Omega(\PP V, \poly))$.
\end{prop}

We record here the following Lemma, which, though standard, will be of use in proving Theorem~\ref{thm:fundcomp}, Lemma~\ref{lem:blow-up-face-of-simplex} and Theorem~\ref{thm:maingeneral}.

\begin{lem}\label{lem:uniqueextend}
For $X$ a smooth irreducible complex projective variety of dimension $n$, if two meromorphic $n$-forms on $X$ agree on a Zariski dense open subset, then they agree on all of $X$.
\end{lem}

\begin{proof}
As $X$ is projective, a meromorphic $n$-form on $X$ is a rational section of the canonical bundle $\omega_X$ \cite[pg.\ 170]{Griffiths--Harris}.
The lemma thus follows from the general fact that if two rational sections of a torsion-free sheaf on an integral scheme agree on a Zariski open subset, then they are equal.\end{proof}

We now prove Theorem~\ref{thm:fundcomp}, reproduced below, postponing the proofs of Propositions \ref{prop:Elessdim} and~\ref{prop:Eequaldim}.

\newtheorem*{thm:fundcomp}{Theorem~\ref{thm:fundcomp}}
\begin{thm:fundcomp}
Let $W\subset V$ be a proper subspace such that $\PP W \cap \poly$ is a (possibly empty) face of $\poly$.
Let $\pi: X = \operatorname{Bl}_{\PP W} \PP V \to \PP V$ be the blow-up, and let the \emph{strict transform} $\widetilde\poly$ of $\poly$ be
\[
\widetilde \poly = \text{the Euclidean closure of }\pi^{-1}(\text{the interior of } \poly)
\text{ in } X(\RR).\]
Then, the pair $(X, \widetilde \poly)$ is a positive geometry whose canonical form is $\pi^*\Omega(\PP V, \poly)$.
\end{thm:fundcomp}

\begin{proof}%[Proof of Theorem~\ref{thm:fundcomp}]
We induct on $\dim \PP V$.
The base case $\dim \PP V = 1$ holds trivially since $\PP W \subset \PP V$ has codimension at most 1 so that the blow-up is an isomorphism on all of $\PP V$.

\smallskip
For the general case, without loss of generality, suppose the codimension of $\PP W \subset \PP V$ is at least 2.
We must verify that the residues of $\pi^*\Omega(\PP V, \poly)$ are the canonical forms of the boundary components of $\partial \widetilde\poly$.
The boundary components are of two possible cases:  
(i) the strict transform $\widetilde {\PP H}$ of a hyperplane $\PP H\subset \PP V$ such that $\poly_H= \PP H \cap \poly$ is a facet of $\poly$, and
(ii) the exceptional divisor $E$.  These are the only possible boundary components because the boundary components of $\partial\poly$ are its facet hyperplanes, $\pi$ is an isomorphism outside of $E$, and $\pi(\widetilde\poly) = \poly$.

\smallskip
In the case of (i), by the universal property of blow-ups, the strict transform $\widetilde{\PP H}$ is the blow-up $\operatorname{Bl}_{\PP W \cap \PP H} \PP H$, whose exceptional divisor $E_H$ is $E\cap \widetilde {\PP H}$, and whose blow-down map $\pi_H: \widetilde{\PP H} \to \PP H$ is the restriction $\pi|_{\widetilde{\PP H}}$.
Let $\widetilde{\poly_{H}}$ be the strict transform of $\poly_{H}$ under $\pi_H$, i.e.\ the Euclidean closure of $\pi_H^{-1}(\intr(\poly_{H}))$ in $\widetilde{\PP H}(\RR)$.
By the induction hypothesis, the pair $(\widetilde{\PP H} , \widetilde{\poly_H})$ is a positive geometry.

To see that $(\widetilde{\PP H} , \widetilde{\poly_H})$ is a boundary component of $(X, \widetilde\poly)$ as a positive geometry, we first note that $\widetilde{\poly_{H}}$ is the Euclidean closure of the interior of $\widetilde{\poly} \cap \widetilde{\PP H}$ in $\widetilde{\PP H}(\RR)$.  Indeed, outside of the codimension 1 closed locus $E_H$, the map $\pi$ identifies $(\widetilde\poly \cap \widetilde{\PP H}) \setminus E_H$ with $\poly_H \setminus \PP W$, so that $\widetilde{\poly} \cap \widetilde{\PP H}$ and $\widetilde{\poly_{H}}$ share the interior $\pi_H^{-1}(\intr(\poly_H))$.
Finally, to verify the appropriate residue property, note that the Zariski open subset $X\setminus E$ of $X$ intersects $\widetilde{\PP H}$ in a Zariski open subset, and $\pi$ is an isomorphism on $X\setminus E$.  Hence, when restricted to $X\setminus E$, we find
\[
\Res_{\widetilde{\PP H}} (\pi^* \Omega(\PP V, \poly)) = 
    \pi_{H}^* \big( \Res_{\PP H} \Omega(\PP V, \poly)\big) = 
    \pi_{H}^* \big( \Omega(\PP H, \poly_{H}) \big) = 
    \Omega(\widetilde{\PP H}, \widetilde{\poly_{H}}),
\]
where the first equality follows because taking residues commutes with pulling back by isomorphisms, and the second equality follows because $(\PP V, \poly)$ is a positive geometry with $(\PP H,\poly_H)$ as a boundary component.
We have an equality of rational forms $\Res_{\widetilde{\PP H}} (\pi^* \Omega(\PP V, \poly)) = \Omega(\widetilde{\PP H}, \widetilde{\poly_{H}})$ when restricted to the Zariski open subset $X\setminus E$, which extends to an equality on all of $X$ (see Lemma~\ref{lem:uniqueextend}).

\smallskip
In the case of (ii), if $\dim \poly_W < \dim \PP W$, then by Proposition~\ref{prop:Elessdim} the exceptional divisor $E$ is not a boundary component, and $\pi^*\Omega(\PP V, \poly)$ has no pole along $E$.
If $\dim \poly_W = \dim \PP W$, then by Proposition~\ref{prop:Eequaldim}, the exceptional divisor $E$ is a boundary component with $E_{\geq 0} = E \cap \widetilde\poly$, whose canonical form is the residue $\operatorname{Res}_E(\pi^*\Omega(\PP V, \poly))$.

\smallskip
We have shown that the form $\pi^*\Omega(\PP V, \poly)$ satisfies the recursive property stated in Definition~\ref{defn:positive-geometry}.
For its uniqueness, if $\Omega'$ is another form with the same recursive property, then the difference $\pi^* \Omega(\PP V, \poly) - \Omega'$, having no poles, is a holomorphic $n$-form on $X$.
This difference is necessarily zero: $X$ has no nonzero holomorphic $n$-forms because $X$ is rational and the geometric genus is a birational invariant of smooth projective varieties \cite[II.8.19]{hartshorne1977algebraic}.
\end{proof}

The rest of this section proves Propositions~\ref{prop:Elessdim} and \ref{prop:Eequaldim}.
For both propositions, we may assume that the codimension of $\PP W \subset \PP V$ is at least 2, since otherwise the statements hold trivially.

\begin{proof}[Proof of Proposition~\ref{prop:Elessdim}]
Recall that $\widetilde\poly \cap E \cong \poly_W \times \poly^W \subset \PP W \times \PP(V/W)$ by Proposition~\ref{prop:exceptionalproduct}. Since $\dim \poly^W = \dim \PP(V/W)$ as long as $\PP W \cap \poly \neq \emptyset$, we find that $\widetilde\poly \cap E$ is full-dimensional in $E$ if and only if $\dim \poly_W = \dim \PP W$.  This proves the first statement.

For the second statement, suppose $\dim \poly_W < \dim \PP W$, and let us fix a triangulation $\TT$ of $\poly$, so that Proposition~\ref{prop:canonical-form-polytope}\ref{prop:canonical-form-triangulation} gives us $\Omega(\PP V, \poly) = \sum_{T\in \TT} \Omega(\PP V, T)$.
We claim that $\pi^*\Omega(\PP V, T)$ has no pole along $E$ for every $T\in \TT$, which implies that $\pi^*(\PP V, \poly)$ has no pole along $E$.
For proof of the claim, let $\partial T$ be the collection of $\dim \PP V +1 $ many facet hyperplanes of a simplex $T\in \TT$.
Since $\dim(T \cap \PP W) < \dim \PP W$, we find that
\[
(\text{the number of hyperplanes in $\partial T$ containing $\PP W$}) < \operatorname{codim}_{\PP V} \PP W.
\]
The claim now follows from the following general Lemma.
\end{proof}

\begin{lem}\label{lem:nopoles}
Let $Y$ be a smooth irreducible closed subvariety of codimension at least 2 in a smooth complex variety $X$ of dimension $n$.  Let $\pi: X' = \operatorname{Bl}_Y X \to X$ be the blow-up with exceptional divisor $E$.
For a collection $D$ of (distinct) smooth irreducible divisors $D_1, \dotsc, D_m$ in $X$, let
\[
\operatorname{mult}_Y D = \text{the number of $D_i$ containing $Y$}.
\]
If $\operatorname{mult}_Y D < \operatorname{codim}_X Y$, then for any meromorphic $n$-form $\Omega$ on $X$ with at most simple poles along $D_1, \dotsc, D_m$, the pullback $\pi^*\Omega$ has no pole along $E$.
\end{lem}

\begin{proof}
Since $\pi$ restricts to an isomorphism on $\pi^{-1}(X\setminus Y)$, it suffices to verify that, for any point $y\in Y$ and an analytic neighborhood $U\ni y$ in $X$, the restriction $\pi^*\Omega|_{\pi^{-1}(U)}$ has no pole along $E$.
Since $Y$ and $D_1, \dotsc, D_m$ are smooth, shrinking $U$ if necessary, we may choose coordinates $x_1, \dotsc, x_n$ of $U$ to have $U\subseteq \AA^n$ and $y = (0, \dots, 0)$ such that $Y|_U$ is the subspace $L = \{x_{k+1} = \cdots = x_n = 0\}$ and each $D_i|_U$ is either empty or is the hyperplane $\{H_i = 0\}$ for some linear function $H_i$.
Let us say $D_1|_U, \dotsc, D_\ell|_U$ are nonempty.  By the assumption that $\operatorname{mult}_Y D < \operatorname{codim}_X Y = n - k$, at most $n-k-1$ many of the functions $H_1, \dotsc, H_\ell$ satisfy $H_i | _L \equiv 0$.

We now carry out the blow-up computation.  We have
$$\operatorname{Bl}_L U= \left\{ \big((x_1,\dots,x_{n}), [Y_{k+1}:\cdots:Y_{n}] \big)\mid x_iY_j = x_jY_i \; \text{for all } i,j \right\} \; \subset \; U \times \PP^{n-k-1} 
    $$
    and $\pi$ is just projection onto the first factor.
We work in the affine open chart 
    $$U' = \{Y_{n} \neq 0\} \subset \operatorname{Bl}_L U,$$
    where $U'$ is considered as an open subset of $\AA^n$ via local coordinates $x_1,\dots,x_k,y_{k+1},\dots, y_{n-1}, x_{n}$ given by $y_i = \frac{Y_i}{Y_n}$. A point 
    $$(x_1,\dots,x_k,y_{k+1},\dots, y_{n-1}, x_{n}) \in U' $$
    corresponds to 
    $$\big((x_1,\dots,x_k,y_{k+1}x_n,\dots,y_{n-1}x_n, x_{n}), [y_{k+1}:\cdots:y_{n-1}:1] \big) \in \operatorname{Bl}_L U.$$
    In this chart, the blow-down map
    $\pi :  U' \to U$  sends 
    \begin{equation}\label{eq:blow-down-coords}
        (x_1,\dots,x_k,y_{k+1},\dots, y_{n-1}, x_{n}) \mapsto (x_1,\dots,x_k,y_{k+1}x_n,\dots,y_{n-1}x_n, x_{n}).
    \end{equation}
Since $\Omega$ has at most simple poles along the $D_i$, we have $\Omega|_U = \frac{f(x)}{H_1\dotsm H_\ell} dx_1\dotsm dx_n$ for a holomorphic function $f$ on $U \subseteq \AA^n$.  Thus, writing $g(x,y)$ for the function $f(x_1, \dotsc, x_k, y_{k+1}x_n, \dotsc, y_{n-1}x_n, x_n)$ on $U'$, we have
\begin{align*}
        \big(\pi^*\, \Omega\big) \big|_{U'} 
        &= \frac{g(x,y)\ dx_1\wedge \cdots \wedge dx_k \wedge d(y_{k+1}x_n) \wedge \cdots \wedge d(y_{n-1}x_n) \wedge dx_{n}}{\prod_{i=1}^{\ell} H_i(x_1,\dots,x_k,y_{k+1}x_n,\dots,y_{n-1}x_n, x_{n})}\\
        &= \frac{g(x,y) \cdot x_n^{n-k-1} \ dx_1\wedge \cdots \wedge dx_k \wedge dy_{k+1} \wedge \cdots \wedge dy_{n-1} \wedge dx_{n}}{\prod_{i=1}^{\ell} H_i(x_1,\dots,x_k,y_{k+1}x_n,\dots,y_{n-1}x_n, x_{n})}.
\end{align*}
Since the exceptional divisor $E$ on our chart $U'$ is given by $\{x_n = 0\}$, we now verify that $\big(\pi^*\, \Omega\big) \big|_{U'}$ has no poles along $\{x_n = 0\}$.
As at most $n-k-1$ many of the $H_i$ satisfy $H_i|_L \equiv 0$, we find that at most $n-k-1$ many of the $H_i$ satisfy
    $$x_n \text{ divides } H_i(x_1,\dots,x_k,y_{k+1}x_n,\dots,y_{n-1}x_n, x_{n}).$$
Therefore, we have at most $n-k-1$ factors of $x_n$ in the denominator of $\big(\pi^*\, \Omega\big) \big|_{U'}$ above, which are canceled by $x_n^{n-k-1}$ in the numerator.
\end{proof}

For the proof of Proposition~\ref{prop:Eequaldim}, we prepare by first considering the case when $\poly$ is a simplex.
One may deduce Lemma~\ref{lem:blow-up-face-of-simplex} by appealing to the fact that projective normal toric varieties are positive geometries \cite[Section 3.6]{NimaLam17}, but we provide an elementary proof.

\begin{lem}\label{lem:blow-up-face-of-simplex}
    Proposition~\ref{prop:Eequaldim} holds when $\poly$ is a simplex.
\end{lem}

\begin{proof}
By a linear change of coordinates, we can assume that $\poly$ is the standard simplex $\Delta^n  = \{[X_0: X_1 : \dotsc : X_n] \mid X_i \geq 0 \text{ for all $i$}\} \subset \PP^n(\RR)$, and that the linear subspace $\PP W$ we are blowing-up is the subspace $\{X_{k+1} = \dots = X_{n} = 0\} \subset \PP^n$.
Note that both the face $\poly_W \subset \PP W$ and the normal polytope $\poly^W \subset \PP(V/W)$ are the standard simplices in $\PP^{k}$ and $\PP^{n-k-1}$, respectively.
Performing a linear change of coordinates if necessary, we may moreover assume that the orientation of the interior of $\poly = \Delta^n \subset \PP^n(\RR)$ is compatible with the standard orientation on the real points of the affine chart $\AA^n \cong U = \{X_0 \neq 0\} \subset \PP^n$ with coordinates $x_i = \frac{X_i}{X_0}$.

Let $\Omega(\Delta^n)$ be the canonical form of $(\PP^n, \Delta^n)$.  Under the blow-up $\pi: \operatorname{Bl}_{\PP W} \PP^n \to \PP^n$, and the identification $E \cong \PP^k \times \PP^{n-k-1}$, we show that $\pi^* \Omega(\Delta^n)$ has a simple pole along $E$, and that $\operatorname{Res}_E \pi^*\Omega(\Delta^n) = \Omega(\Delta^k) \wedge \Omega(\Delta^{n-k-1})$.
By Lemma~\ref{lem:uniqueextend}, it suffices to show these on a Zariski dense open subset $U'$ of $\operatorname{Bl}_{\PP W} \PP^n$ such that $U'\cap E$ is nonempty (and hence dense in $E$).
To this end, consider the affine chart $\AA^n \cong U = \{X_0 \neq 0\} \subset \PP^n$, over which we have $\operatorname{Bl}_L \AA^n \cong \pi^{-1}(U)$ where $L$ is the subspace $\{x_{k+1} = \dotsc = x_n = 0\}$.  
Recalling that
\[
\operatorname{Bl}_L \AA^n = \left\{ \big((x_1,\dots,x_{n}), [Y_{k+1}:\cdots:Y_{n}] \big)\mid x_iY_j = x_jY_i \; \text{for all } i,j \right\} \; \subset \; \AA^n \times \PP^{n-k-1},
\]
we consider the affine open chart $U' = \{Y_n \neq 0\}$ in $\operatorname{Bl}_{L} \AA^n$.
As in the proof of Lemma~\ref{lem:nopoles}, let $x_1,\dots,x_k,y_{k+1},\dots, y_{n-1}, x_{n}$ be local coordinates on $U'$ given by $y_i = \frac{Y_i}{Y_n}$.  The exceptional divisor $E$ is given by $\{x_n = 0\}$ in this chart, and the blow-down map
    $\pi :  \AA^n \cong U' \to U\cong \AA^n$ is as in~\eqref{eq:blow-down-coords}. 
    %$$(x_1,\dots,x_k,y_{k+1},\dots, y_{n-1}, x_{n}) \mapsto (x_1,\dots,x_k,y_{k+1}x_n,\dots,y_{n-1}x_n, x_{n}).$$
    
The determinant of the Jacobian matrix of this map is a power of $x_n$, and is thus positive when $x_n >0$.  Hence, because the points in $\intr(\poly)$ and $\pi^{-1}(\intr(\poly))$ satisfy $x_n > 0$, we consider the standard orientation on $U'(\RR)$, given by the ordered coordinates $(x_1, \dotsc, x_k, y_{k+1},\dots, y_{n-1}, x_{n})$, which is compatible via $\pi$ to the standard orientation on $U(\RR)$ that we picked, given by the ordered coordinates $(x_1, \dotsc, x_n)$.

\smallskip
By Proposition~\ref{prop:canonical-form-polytope}\ref{prop:canonical-form-general-simplex}, the canonical form $\Omega = \Omega(\PP^n, \Delta^n)$ is given by 
    $\Omega|_{U} = \frac{dx_1}{x_1} \wedge \cdots \wedge \frac{dx_n}{x_n}.$
Thus, the pullback $\pi^* \Omega$ in the local coordinates on $U'$ is
    \begin{align*}
        \big(\pi^* \, \Omega \big) \big|_{U'} &= \frac{dx_1}{x_1} \wedge \cdots \wedge \frac{dx_k}{x_k} \wedge \frac{d(y_{k+1}x_n)}{y_{k+1}x_n} \wedge \cdots \wedge \frac{d(y_{n-1}x_n)}{y_{n-1}x_n} \wedge \frac{dx_n}{x_n}\\
        &= \frac{dx_1}{x_1} \wedge \cdots \wedge \frac{dx_k}{x_k} \wedge \frac{dy_{k+1}}{y_{k+1}} \wedge \cdots \wedge \frac{dy_{n-1}}{y_{n-1}} \wedge \frac{dx_n}{x_n}.
    \end{align*}
In particular, the pullback has a simple pole along $\{x_n = 0\} = E \cap U'$, as desired.  We compute the residue along $E$:
    \begin{align*} \big(\Res_E (\pi^* \, \Omega ) \big)\big|_{E\cap U'} &= 
    \Res_{\{x_n = 0\}} \big( (\pi^* \, \Omega)|_{U'}\big) \\
    &= \left(\frac{dx_1}{x_1} \wedge \cdots \wedge \frac{dx_k}{x_k}\right) \wedge \left(\frac{dy_{k+1}}{y_{k+1}} \wedge \cdots \wedge \frac{dy_{n-1}}{y_{n-1}}\right)\\
    &= \Omega(\PP^k, \Delta^{k})\big|_{X_0 \neq 0} \wedge \Omega(\PP^{n-k-1},\Delta^{n-k-1}) \big|_{Y_n \neq 0}\\
    &= \big(\Omega(\poly_W) \wedge \Omega(\poly^W)\big) \big|_{E \cap U'},
    \end{align*}
where the third and fourth equalities are justified as follows.
For $\PP^k$ with homogeneous coordinates $X_0, \dots,X_k$, the form $\frac{dx_1}{x_1} \wedge \cdots \wedge \frac{dx_k}{x_k}$ is the canonical form $\Omega(\PP^k, \Delta^k)$ restricted to the chart $\AA^k \cong \{X_0 \neq 0\} \subset \PP^k$ with the standard orientation.
Similarly, for $\PP^{n-k-1}$ with homogeneous coordinates $Y_{k+1}, \dots,Y_n$, the form $\frac{dy_{k+1}}{y_{k+1}} \wedge \cdots \wedge \frac{dy_{n-1}}{y_{n-1}}$ is the canonical form $\Omega(\PP^{n-k-1}, \Delta^{n-k-1})$ restricted to the chart $\AA^{n-k-1} \cong \{Y_n \neq 0\} \subset \PP^k$ with the standard orientation. 
Under the natural identification $E \cong \PP^k \times \PP^{n-k-1}$, our open chart $E \cap U' \subset E$ corresponds precisely to the product of these charts $\AA^k \times \AA^{n-k-1} \subset \PP^k \times \PP^{n-k-1}$.
Hence, the standard orientations on the two affine charts $\AA^k$ and $\AA^{n-k-1}$ induce the standard orientation on $\AA^k \times \AA^{n-k-1} = \AA^{n-1} \cong E \cap U'$, which is the orientation induced by the restriction of our standard orientation on $U'$ to the boundary $E \cap U' = \{x_n = 0\}$ (see Remark~\ref{rem:orientation}).
\end{proof}

\begin{proof}[Proof of Proposition~\ref{prop:Eequaldim}]
Fix a triangulation $\TT$ of $\poly$.
By Proposition~\ref{prop:canonical-form-polytope}\ref{prop:canonical-form-triangulation}, we have
\[
\Omega(\PP V, \poly) = \sum_{T\in \TT} \Omega(\PP V, T).
\]
Let $\TT^{\text{full}} = \{ T\in \TT \mid \dim T_W = \dim \PP W\}$.
By Proposition~\ref{prop:Elessdim}, if $T\not\in \TT^{\text{full}}$ then $\pi^*\Omega(\PP V, T)$ has no pole along the exceptional divisor $E$.  Thus, applying Lemma~\ref{lem:blow-up-face-of-simplex} to simplices in $\TT^{\text{full}}$, we find
\[
\Res_E (\pi^*\Omega(\PP V, \poly)) = \sum_{T \in \TT^{\text{full}}} \Res_E \big(\pi^*\Omega(\PP V, T)\big) = \sum_{T \in \TT^{\text{full}}} \Omega(\PP W, T_W) \wedge \Omega(\PP(V/W), T^W).
\]
Since $\dim \PP W = \dim \poly_W$, Lemma~\ref{lem:induced-triangulations}\ref{lem:triangulation-of-face} implies that the collection $\TT_W = \{T_W \mid T\in \TT^{\text{full}}\}$ is a triangulation of $\poly_W$, and Lemma~\ref{lem:induced-triangulations}\ref{lem:triangulation-of-normal} implies that for any simplex $F$ in the triangulation $\TT_W$ of $\poly_W$, the collection $\TT^{W,F} = \{T^W \mid T\in \TT^{\text{full}} \text{ satisfying } T_W = F\}$ is a triangulation of $\poly^W$.
Thus, we compute that
    \begin{align*}
        \Res_E (\pi^*\Omega(\PP V, \poly)) &= \sum_{T \in \TT^{\text{full}}} \Omega(\PP W, T_W) \wedge \Omega(\PP(V/W), T^W)\\
        &=\sum_{F \in \TT_W} \sum_{\substack{T \in \TT^{\text{full}}\\ T_W = F}} \Omega\left(\PP W, T_W\right) \wedge \Omega\left(\PP(V/W), T^W\right) \\
        &= \sum_{F \in \TT_W} \Omega(\PP W, F) \wedge \left( \sum_{\substack{T \in \TT^{\text{full}}\\ T_W = F}} \Omega\left(\PP(V/W), T^W\right)\right)\\
        &= \sum_{F \in \TT_W} \Omega(\PP W, F) \wedge \Omega\left(\PP(V/W), \poly^W\right) \\
        &= \left(\sum_{F \in \TT_W} \Omega(\PP W, F)\right) \wedge \Omega\left(\PP(V/W), \poly^W\right) \\
        &= \Omega\left(\PP W, \poly_W\right) \; \wedge \; \Omega\left(\PP(V/W), \poly^W\right).
    \end{align*}
Here, for the fourth and sixth equalities, we used Proposition~\ref{prop:canonical-form-polytope}\ref{prop:canonical-form-triangulation} applied to the triangulations $\TT_W$ and $\TT^{W,F}$ of $(\PP W, \poly_W)$ and $(\PP(V/W), \poly^W)$, respectively, which are positive geometries because $\poly_W \subset \PP W$ and $\poly^W \subset \PP(V/W)$ are full-dimensional polytopes.

Lastly, that $(E, E_{\geq 0})$ is a positive geometry with the canonical form $\operatorname{Res}_E(\pi^*\Omega(\PP V, \poly))$
now follows from Lemma \ref{lem:productposgeom} below, which states that products of positive geometries are positive geometries.
\end{proof}

\begin{lem}\label{lem:productposgeom}\cite[Section 2.3]{NimaLam17}
If $(X,X_{\geq 0})$ and $(Y,Y_{\geq 0})$ are positive geometries with canonical forms $\Omega(X_{\geq 0})$ and $\Omega(Y_{\geq 0})$, respectively, then $(X \times Y, X_{\geq 0} \times Y_{\geq 0})$ is a positive geometry with canonical form $\Omega(X_{\geq 0}) \wedge \Omega(Y_{\geq 0})$.
\end{lem}
We will use Lemma \ref{lem:productposgeom} again in Theorem~\ref{thm:maingeneral}. 
\section{Blowing up a sequence of faces}\label{section:blow-up-sequence}

We prove Theorem~\ref{thm:intro_mainthm}, our  main result, in this section.  In Section~\ref{section:wonderful-boundaries}, we prepare by recording facts about the boundary structure of a wonderful compactification $X^\build$ and a wondertope $\widetilde\poly^\build$.
In Section~\ref{section:main-thm-proof}, we prove the main theorem using Theorem~\ref{thm:fundcomp} and induction.

\subsection{Wonderful compactifications and their boundaries} \label{section:wonderful-boundaries}

Let $\build$ be a set of proper linear subvarieties of $\PP V$.
We fix an ordering $\build = \{\PP W_1, \dotsc, \PP W_k\}$ such that $i\leq j$ if $W_i \subseteq W_j$, and consider the sequential blow-up
\[
\pi_\build : \PP V^\build = \operatorname{Bl}_{\PP W_k}(\dotsb (\operatorname{Bl}_{\PP W_2}(\operatorname{Bl}_{\PP W_1}\PP V))\dotsb ) \to \PP V.
\]
(As before, we repurpose the notation to write $\PP W_2$ also for its strict transform in $\operatorname{Bl}_{\PP W_1} \PP V$, and similarly for $\PP W_3, \dotsc, \PP W_k$).

\medskip
Assume now that $\build$ is a building set (as defined in the introduction).  In this case, De Concini and Procesi called $X^\build$ the \emph{wonderful compactification} of the arrangement complement $\PP V \setminus (\bigcup \build)$, and showed that its boundary has the following structure.
For $\PP W \in \build$, define
\begin{align*}
\build_W &= \{F \cap \PP W \mid F\in \build \text{ and } F\not\supseteq \PP W\} \quad\text{and}\\
\build^W &= \{\PP(W'/W) \subseteq \PP(V/W) \mid \PP W' \in \build \text{ and } W'\supseteq W\},
\end{align*}
and let $E_W \subset X^\build$ be (the strict transform of) the exceptional divisor of the blow-up at $\PP W$.

\begin{prop}\label{prop:wonderful}\cite[Theorem 4.3]{de1995wonderful}
For $\PP W\in \build$, both $\build_W$ and $\build^W$ are building sets in $\PP W$ and $\PP(V/W)$, respectively, and there is a natural isomorphism
$$ E_W \cong \PP W^{\build _ W} \times  \PP(V / W)^{\build ^{W}}$$
along with a commuting diagram of natural maps
    $$\begin{tikzcd}
    X^\build \arrow[d, "\pi_\build"]
        & &  E_W \cong \PP W^{\build_W} \times  \PP(V / W)^{\build ^{W}} \arrow[ll, "\supset"'] \arrow[d, "\pi_{B_W} \times \pi_{\build^{W}}"] \\
    \PP V  
        & \Bl_{\PP W} \PP V \arrow[l, "\pi_W"]
        & E \cong \PP W \times \PP(V/W) \arrow[l, "\supset"']
    \end{tikzcd}.$$
\end{prop}

In analogy with Proposition \ref{prop:wonderful}, we show that a wondertope has a similar boundary structure, as follows.
Let $\poly \subset \PP V$ be a full-dimensional polytope, and assume further now that $\build$ satisfies the property:
\begin{quote}
for every $F\in \build$, the intersection $F\cap \poly$ is a (possibly empty) face of $\poly$.
\end{quote}
Let
$
\widetilde\poly^\build$ be the Euclidean closure of $\pi_\build^{-1}(\intr(\poly))$ in $X^\build(\RR)$. This is technically not a wondertope since we have not yet assumed that all facet hyperplanes of $\poly$ are in $\build$.
In the following two lemmas, we will record some properties of the boundary structure of $\widetilde\poly^\build$.

\smallskip
For $\PP H \subset \PP V$ a hyperplane such that $\poly_H = \poly \cap \PP H$ is a facet of $\poly$, we define
\[
\build_H = \{F \cap \PP H \mid F \in \build \text{ and } F \not\supseteq \PP H\},
\]
whose ordering is inherited from the order on $\build$.
Let $\widetilde{\PP H}$ be the strict transform of $\PP H$ under $\pi_\build$, which by the universal property of blow-ups is isomorphic to the sequential blow-up \[ \pi_{\build_H}: \PP H^{\build_H} \to \PP H. \]

We define
\begin{align*}
\widetilde{\PP H}_{\geq 0} &= \text{the Euclidean closure of the interior of $\widetilde\poly^\build \cap \widetilde{\PP H}$ inside $\widetilde{\PP H}(\RR)$}, \quad\text{and}\\
\widetilde{\poly_H}^{\build_ H}&= \text{the Euclidean closure of $\pi_{\build_H}^{-1}(\intr(\poly_H))$ in $\PP H^{\build_H}(\RR)$}.
\end{align*}
We caution that $\build_H$ may not be a building set in $\PP H$, so we will soon impose an additional condition on $\build$ (see Definition~\ref{defn:well-adapted}) when we prove Theorem~\ref{thm:intro_mainthm}.
Moreover, we caution that taking Euclidean closure of the interior in the definition of $\widetilde{\PP H}_{\geq 0}$ is necessary for the following lemma: otherwise, the intersection $\widetilde\poly^\build \cap \widetilde{\PP H}$ (without taking the closure of the interior) may not be equal to $\widetilde{\poly_H}^{\build_H}$, so the following lemma would fail. This is illustrated in Example \ref{eg:pyramid1}. 

\begin{lem}\label{lemma:facet-transform}
Under the isomorphism $\widetilde{\PP H} \cong \PP H^{\build_H}$, we have $\widetilde{\PP H}_{\geq 0} = \widetilde{\poly_H}^{\build_ H}$.
\end{lem}

\begin{proof}
Since the exceptional locus intersects $\widetilde{\PP H}$ in a subset of codimension at least $1$, we can detect (a dense open subset of) the interior of $\widetilde\poly^\build \cap \widetilde{\PP H}$ inside $\widetilde{\PP H}(\RR)$ by looking away from the exceptional locus. We note that $\pi_\build$ is an isomorphism away from the exceptional locus. Therefore, outside of the exceptional locus, we know that $\widetilde\poly^\build \cap \widetilde{\PP H}$ is given by $\pi_\build^{-1}(\poly \cap H) = \pi_\build^{-1}(\poly_H)$.
    So a dense open subset of $\intr(\widetilde\poly^\build \cap \widetilde{\PP H})$
    is given by $\pi_\build^{-1}(\intr(\poly_H))$.
    We know that, restricted to $\widetilde{\PP H}$, our wonderful compactification $\pi_{\build}$ looks like $\pi_{\build_H} : \widetilde{\PP H} \cong \PP H^{\build_H} \to \PP H$, so $\pi_\build^{-1}(\intr(\poly_H)) = \pi_{\build_H}^{-1}(\intr(\poly_H)) = \intr(\widetilde{\poly_H}^{\build_H})$. Therefore, 
    \[\widetilde{\PP H}_{\geq 0} = \overline{\intr(\widetilde\poly^\build \cap \widetilde{\PP H})}^{\text{closure in } \widetilde{\PP H}} = 
    \overline{\intr(\widetilde{\poly_H}^{\build_H})}^{\text{closure in }\PP H^{\build_H}}
    =\widetilde{\poly_H}^{\build_H}.\qedhere
    \]
\end{proof}

    Let $\PP W \in \build$ be a linear subspace intersecting our polytope $\poly$ along a face $\poly_W = \poly \cap \PP W$, and let $E_W$ be (the strict transform of) the exceptional divisor in $X^\build$ corresponding to $\PP W$.
Let
\[
E_{W,\geq 0} = \text{the Euclidean closure of the interior of $\widetilde\poly^\build \cap E_W$ inside $E_W(\RR)$}.
\]
    
\begin{lem} \label{lem:exceptional-transform}
Under the isomorphism $ E_W \cong (\PP W)^{\build _ W} \times  \PP(V / W)^{\build ^{W}}$, we have
    $$E_{W,\geq 0} =
        \widetilde{\poly_W}^{\build_ W} \times \widetilde{\poly^W}^{\build^W}.$$
In particular, the subset $E_{W,\geq 0}$ is full-dimensional in $E_W$ if and only if $\dim \poly_W = \dim \PP W$.
\end{lem}
\begin{proof}
    As in Lemma~\ref{lemma:facet-transform}, because the rest of the exceptional locus intersects $E_W$ in a subset of codimension at least $1$, we can detect a dense open subset of $\intr(\widetilde\poly^\build \cap E_W)$ by looking away from this additional exceptional locus. 
    Away from the exceptional locus corresponding to $\build\setminus \{\PP W\}$, our blow-up $\pi_{\build}$ looks like (a dense open subset of) a single blow-up
    $$\pi_W: \Bl_{\PP W} \PP V \to \PP V,$$
    with exceptional divisor $E \cong \PP W \times \PP(V/W)$. Therefore, we know that, away from the exceptional locus of $\build\setminus \{\PP W\}$, $\widetilde\poly^\build \cap E_W$ is given by $\widetilde\poly^{\{\PP W\}} \cap E$. So (a dense open subset of) the interior $\intr(\widetilde\poly^\build \cap E_W)$ is given by (a dense open subset of) $\intr(\poly_W) \times \intr(\poly^W) \subset \PP W \times \PP(V/W) \cong E$, according to Proposition~\ref{prop:exceptionalproduct}. In particular, this is empty if $\poly_W$ is not full-dimensional in $\PP W$. So let's assume that $\dim \poly_W = \dim \PP W$.

    Under the canonical map 
    $$\pi_{B_W} \times \pi_{\build^{W}}: \left(E_W \cong \PP W^{\build_W} \times  \PP(V / W)^{\build^{W}}\right) \to \big(E \cong \PP W \times \PP(V/W)\big),$$
    the inverse image of $\intr(\poly_W) \times \intr(\poly^W)$ is precisely $\intr\left(\widetilde{\poly_W}^{\build_W} \right) \times \intr\left(\widetilde{\poly^W}^{\build^W}\right)$, and so
\begin{multline*}
E_{W,\geq 0} = \overline{\intr\left(\widetilde\poly^W \cap E_W\right)}^{\text{closure in } E_W} \\
= \overline{\intr\left(\widetilde{\poly_W}^{\build_W} \right) \times \intr\left(\widetilde{\poly^W}^{\build^W}\right)}^{\text{closure in } \PP W^{\build_W} \times \PP(V/W)^{\build^W}} \\
= \widetilde{\poly_W}^{\build_W} \times \widetilde{\poly^W}^{\build^W}.
\qedhere
\end{multline*}
\end{proof}

\subsection{Proof of the main theorem}\label{section:main-thm-proof}

We are now ready to prove Theorem~\ref{thm:intro_mainthm}.
We in fact prove a stronger version, stated as follows.
Let $\poly \subset \PP V$ be a full-dimensional polytope.
Let $\build$ be a set of proper linear subvarieties of $\PP V$ satisfying the following properties stated in the introduction:
\begin{itemize}
\item for every $F\in \build$, the intersection $F \cap \poly$ is a (possibly empty) face of $\poly$, and
\item $\build$ is a building set.
\end{itemize}
Furthermore, we impose the following relaxation of the condition \ref{item:hyperplane-condition} stated in the introduction.

\begin{defn}\label{defn:well-adapted}
We say that such a set $\build$ is \emph{well-adapted to $\poly$} if either $\dim \PP V = 1$ or, otherwise, $\build$ further satisfies the recursive property:
\begin{itemize}
\item for every facet hyperplane, i.e.\ a hyperplane $\PP H \subset \PP V$ such that $\poly_H = \poly \cap \PP H$ is a facet of $\poly$, the set $\build_H$ is a building set in $\PP H$ and is well-adapted to $\poly_H$.
\end{itemize}
\end{defn}

We prove the following theorem. 

\begin{thm}\label{thm:maingeneral}
Let $\build$ be well-adapted to $\poly$.  Then, the pair $(X^\build, \widetilde \poly^\build)$ is a positive geometry whose canonical form is the pullback $\pi_\build^*\Omega(\PP V, \poly)$ of the canonical form of $(\PP V, \poly)$.  The boundary components of $(X^\build, \widetilde\poly^\build)$ are the exceptional divisors $E_F$ for $F\in \build$ such that $\dim(\poly \cap F) = \dim F$, together with the strict transforms of the facet hyperplanes of $\poly$ that are not in $\build$.
\end{thm}

\begin{proof}[Proof of Theorem~\ref{thm:intro_mainthm} from Theorem~\ref{thm:maingeneral}]
If $\build$ contains all the facet hyperplanes of $\poly$, then $\build_H$ is a building set in $\PP H$ by Proposition~\ref{prop:wonderful}. Moreover, $\build_H$ contains all the facet hyperplanes of $\poly_H$, since a facet hyperplane of $\poly_H$ is an intersection of $\PP H$ with another facet hyperplane of $\poly$.
Hence, the set $\build$ is well-adapted to $\poly$ in this case.
In particular, Theorem~\ref{thm:maingeneral} contains Theorem~\ref{thm:intro_mainthm} as a special case.
\end{proof}

\begin{rem}
If $\build$ consists of linear subspaces that intersect $\poly$ along (possibly empty) faces and that are closed under intersection (i.e.\ if $F,F' \in \build$ then $F\cap F' \in \build$), then $\build$ is a building set, and $\build_H$ is intersection-closed for any hyperplane $\PP H$. Hence, the set $\build$ is well-adapted to $\poly$.
In this intersection-closed setting, Brown \cite[Section 5]{Brown23} calls the wondertope $\widetilde\poly^\build$ a ``blow-up of a polyhedral linear configuration'' and describes the structure of its boundary divisors.
\end{rem}

Theorem~\ref{thm:maingeneral} fails without the hypothesis that $\build$ is well-adapted to $\poly$, even if $\build$ is a building set; see Example~\ref{eg:pyramid3}.

\begin{proof}[Proof of Theorem~\ref{thm:maingeneral}]
We induct on $\dim \PP V$.  The base case $\dim \PP V = 1$ holds trivially since the blow-down map is an isomorphism onto all of $\PP V$ in that case.

\smallskip
For the general case,
we must verify that the residues of $\pi^*\Omega(\PP V, \poly)$ are the canonical forms of the boundary components of $\partial \widetilde\poly$.
The boundary components come in two cases:
(i) the strict transform $\widetilde {\PP H}$ of a hyperplane $\PP H\subset \PP V$ such that $\poly_H$ is a facet of $\poly$ and $\PP H \notin \build$, and
(ii) the exceptional divisor $E_W$ for $\PP W \in \build$.
That there are no more boundary components follows from the fact that the boundary components of $\partial\poly$ are its facet hyperplanes, $\pi_\build$ is an isomorphism outside of the exceptional locus $\bigcup_{\PP W \in \build} E_W$,  and $\pi_\build(\widetilde\poly^\build) = \poly$.

\smallskip
In the case of (i), by the universal property of blow-ups, we have $\widetilde{\PP H} \cong \PP H^{\build_H}$, with $\pi_{\build H}$ identified with the restriction $\pi_{\build}|_{\widetilde{\PP H}}$.
By the induction hypothesis, the pair $(\widetilde{\PP H}, \widetilde{\poly_H}^{\build_H})$ is a positive geometry with canonical form $\pi_{\build_H}^* \Omega(\PP H, \poly_H)$.
The pair is also a boundary component of $(\PP V^\build, \widetilde\poly^\build)$ by Lemma~\ref{lemma:facet-transform}.
To verify that $\operatorname{Res}_{\widetilde{\PP H}} \pi_\build^*\Omega(\PP V, \poly) = \pi_{\build_H}^* \Omega(\PP H, \poly_H)$, by Lemma~\ref{lem:uniqueextend} it suffices to do so on the Zariski dense open subset $U = \PP V^\build \setminus (\bigcup_{\substack{\PP W \in \build}} E_W)$, which restricts to a Zariski dense open subset in $\widetilde{\PP H}$.
Indeed, since $\pi_\build$ is an isomorphism on $U$, we find that over $U$, we have
\[
\Res_{\widetilde{\PP H}} (\pi_\build^* \Omega(\PP V, \poly)) = 
    \pi_{\build_H}^* \big( \Res_{\PP H} \Omega(\PP V, \poly)\big) = 
    \pi_{\build_H}^* \big( \Omega(\PP H, \poly_{H}) \big),
\]
as desired.

\smallskip
In the case of (ii), we first note that it suffices to consider only $\PP W\in \build$ such that $\dim \poly_W = \dim \PP W$.
Indeed, if $\dim \poly_W < \dim \PP W$ by Proposition~\ref{prop:Elessdim} and the commuting diagram in Proposition~\ref{prop:wonderful}, we find that $\pi^*\Omega(\PP V, \poly)$ has no pole along $E_W$.  Meanwhile, Lemma~\ref{lem:exceptional-transform} implies that $E_{W,\geq 0}$ is not full-dimensional if $\dim \poly_W < \dim \PP W$.

When $\dim \poly_W = \dim \PP W$, 
by the induction hypothesis and Lemma~\ref{lem:productposgeom},
the pair
\[
\left(\PP W^{\build _W} \times \PP(V/W)^{\build^W}, \widetilde{\poly_W}^{\build_W}\times \widetilde{\poly^W}^{\build^W} \right)
\]
is a positive geometry.

By Proposition~\ref{prop:wonderful} and Lemma~\ref{lem:exceptional-transform}, we have an isomorphism of pairs
\[
(E_W, E_{W, \geq 0}) \cong \left(\PP W^{\build _W} \times \PP(V/W)^{\build^W}, \widetilde{\poly_W}^{\build_W}\times \widetilde{\poly^W}^{\build^W}\right),
\]
and in particular $(E, E_{W,\geq 0})$ is a boundary component of $(\PP V^\build, \widetilde\poly^\build)$.

To verify that $\Omega(E_W, E_{W,\geq 0}) = \operatorname{Res}_{E_W} \pi_{\build}^* \Omega(\PP V, \poly)$, by Lemma~\ref{lem:uniqueextend} it suffices to do so on the Zariski dense open subset $U_W = \PP V^\build \setminus (\bigcup_{\substack{\PP W'\in \build \\ \PP W' \neq \PP W}}E_{W'})$, which restricts to a Zariski dense open subset in $E_W$.  
Indeed, on $U_W$, Proposition~\ref{prop:wonderful} identifies $E_W$ with the exceptional divisor $E$ of the blow-up $\pi_W: \operatorname{Bl}_{\PP W} \PP V \to \PP V$, and hence, over $U_W$, we have
\begin{align*}
\operatorname{Res}_{E_W} \pi_{\build}^* \Omega(\PP V, \poly) &= \operatorname{Res}_E \pi_W^* \Omega(\PP V, \poly)\\
&= \Omega\left(\PP W \times \PP(V/W), \poly_W \times \poly^W\right) \\
&= \Omega\left(\PP W^{\build _W} \times \PP(V/W)^{\build^W}, \widetilde{\poly_W}^{\build_W}\times \widetilde{\poly^W}^{\build^W}\right)
\end{align*}
where the second equality follows from Theorem~\ref{thm:fundcomp}, and the last equality follows from that the pair $(E, E_{\geq 0})$ in Theorem~\ref{thm:fundcomp} coincides with the pair $(E_W, E_{W, \geq 0})$ on $U_W$.

\smallskip
Lastly, the uniqueness of the canonical form follows from the rationality of $X^\build$, similarly to the argument at the end of the proof of Theorem~\ref{thm:fundcomp}.
\end{proof}

\section{Examples and questions}\label{section:examples}

In Section~\ref{section:braid}, we give a detailed account of how our Theorem~\ref{thm:intro_mainthm} recovers a result about the moduli space of pointed stable rational curves.
In Section~\ref{section:nonexamples}, we collect several examples where a blow-up of a positive geometry does not result in a positive geometry.
In Section~\ref{section:detailed-example}, we detail an example that both illustrates Theorem~\ref{thm:fundcomp} and provides ways to produce pathologies.
We finish with an example with a view towards a generalization of Theorem~\ref{thm:maingeneral} to a setting beyond polytopes.

\subsection{Motivating example: $\overline{M}_{0, n+1}$ and the braid arrangement}\label{section:braid}

Let $\overline{M}_{0,n+1}$ be the moduli space of stable $(n+1)$-pointed rational curves; we point to \cite[Chapter 1]{KV07} as a reference.
In \cite[Proposition 8.2]{AHLcluster}, the authors show that 
$(\overline{M}_{0,n+1}, (\overline{M}_{0,n+1})_{\geq 0})$ 
is a positive geometry, using theory developed in \cite[Proposition 2.7]{brown2010algebra}. 
Our work provides a self-contained proof of this fact by viewing $\overline{M}_{0,n+1}$ as a certain wonderful compactification of the braid arrangement complement.

\subsubsection{$M_{0,n+1}$ as the braid arrangement complement}
The space $M_{0,n+1}$ is the moduli space of $n+1$ distinct points $z_0,z_1,\dots,z_{n}$ in $\PP^1$, up to projective linear transformations of $\PP^1$. Fixing $(z_0, z_{n-1}, z_n) = (0,1,\infty)$ allows us to write
\begin{align*}
    M_{0,n+1} &= \big \{ (z_1,z_2, \dots, z_{n-2}) \in (\PP^1)^{n-2}: z_i \neq 0,1,\infty \text{ and } z_i \neq z_j \text{ for } i \neq j \big \}\\
    &\cong \big \{ (z_1,z_2, \dots, z_{n-2}) \in \AA^{n-2}: z_i \neq 0,1 \text{ and } z_i \neq z_j \text{ for } i \neq j \big \}.
\end{align*}

Let $V = \RR^n / \RR\cdot(1, 1, \dotsc, 1)$. Recall that the \emph{rank $(n-1)$ braid arrangement} $A_{n-1}$ consists of the image of the hyperplanes $\{ x_i - x_j = 0 \} \subset \RR^n$ in $\PP V$ for $1 \leq i < j \leq n$ (see for example~\cite[~Ch.~1]{stanley2007introduction}.)

It is well-known that $M_{0,n+1}$ is isomorphic to the complement $C(A_{n-1}) := \PP V \setminus A_{n-1}$ of the braid arrangement. We will go through this isomorphism in detail because it is useful in deriving the \emph{Parke-Taylor form} from Corollary \ref{cor:m0n-positivegeometry}.

Consider $V \cong \RR^{n-1}$ with basis given by the images $\overline{e_i} \in V$ of the standard basis vectors $e_i \in \RR^n$ for $1\leq i < n$. The projectivization $\PP V \cong \PP^{n-2}$ thus consists of non-zero points $y_1 \overline{e_1} + \cdots + y_{n-1}\overline{e_{n-1}}$, modulo scalar multiplication. The image of the hyperplanes $\{x_i-x_j = 0\} \subset \RR^n$ in $\PP V$ are the hyperplanes $\{y_i = 0\}$ for $1\leq i\leq n-1$ and $\{y_i - y_j = 0\}$ for $1\leq i < j\leq n-1$. Therefore, the complement $C(A_{n-1})$ can be written as
\begin{align*}
    C(A_{n-1}) 
    &= \big \{ [y_1 \overline{e_1} + \cdots + y_{n-1}\overline{e_{n-1}}] \in \PP V :  y_i \neq 0 \text{ and } y_i \neq y_j \textrm{ for } i \neq j \big \}\\
    &\cong \big\{ [ y_1: \cdots : y_{n-1} ] \in \PP^{n-2} :  y_i \neq 0 \text{ and } y_i \neq y_j \textrm{ for } i \neq j \big \}.
\end{align*}
The isomorphism $C(A_{n-1}) \cong M_{0,n+1}$ then follows via the change of coordinates $z_i = y_i/y_{n-1}$ for $1\leq i \leq n-2$. Indeed, we obtain:
\begin{align*}
    C(A_{n-1})
    & \cong \left\{ \left[z_1: \cdots: z_{n-2}: 1\right] \in \PP^{n-2} :  z_i \neq 0,1 \text{ and } z_i \neq z_j \textrm{ for } i \neq j \right \}\\
    &\cong \big \{ (z_1,z_2, \dots, z_{n-2}) \in \RR^{n-2}: z_i \neq 0,1 \text{ and } z_i \neq z_j \text{ for } i \neq j \big \}.
\end{align*}

\subsubsection{Wonderful Compactification} \label{section:minimal-building-set}
Under the correct choice of building set for $C(A_{n-1})$, the corresponding wonderful compactification is isomorphic to the Deligne--Knudsen--Mumford compactification $\overline{M}_{0,n+1}$ \cite[Remark 4.3.(3)]{de1995wonderful}. This building set is called the \emph{minimal building set} $\BB^{\textrm{min}}$, defined as 
\[ \BB^{\textrm{min}} := \bigg \{ \PP W : W = \{ x_{j_1} = x_{j_2} = \cdots = x_{j_k}  \} \text{ for } 0 \leq j_1 < \cdots < j_k \leq n \textrm{ and } k \geq 2 \bigg \}.  \]
There is a nice combinatorial restatement of this story. Let $\LL(A_{n-1})$ be the poset with elements given by intersections of hyperplanes in $A_{n-1}$, ordered by reverse inclusion. It is in fact a lattice, and is well known to be isomorphic to $\Pi_n$, the lattice of set partitions of $[n]:= \{ 1, 2, \cdots, n \}$, ordered by reverse refinement.
Under this isomorphism, the set partition $\pi$ corresponds to the subspace $W$ where $x_i = x_j$ in $W$ if and only if $i,j$ are in the same block in $\pi$. For example, the set partition $\{  1 3 \ | \ 2  4 6 \  | \ 5 \} $ corresponds to the subspace in $\LL(A_{5})$ where $x_{1} = x_{3}$ and $x_2 = x_4 = x_6$.

\begin{figure}[!h]
\begin{minipage}{0.6\textwidth}
	\centering
 \resizebox{9cm}{!}{
\begin{tikzpicture}[x=0.75pt,y=0.75pt,yscale=-1,xscale=1]

%Straight Lines [id:da2690659488355902] 
\draw    (120,180) -- (320,260) ;
\draw [shift={(120,180)}, rotate = 21.8] [color={rgb, 255:red, 0; green, 0; blue, 0 }  ][fill={rgb, 255:red, 0; green, 0; blue, 0 }  ][line width=0.75]      (0, 0) circle [x radius= 3.35, y radius= 3.35]   ;
%Straight Lines [id:da25674655173825844] 
\draw    (200,180) -- (320,260) ;
\draw [shift={(200,180)}, rotate = 33.69] [color={rgb, 255:red, 0; green, 0; blue, 0 }  ][fill={rgb, 255:red, 0; green, 0; blue, 0 }  ][line width=0.75]      (0, 0) circle [x radius= 3.35, y radius= 3.35]   ;
%Straight Lines [id:da6949950022004908] 
\draw    (280,180) -- (320,260) ;
\draw [shift={(280,180)}, rotate = 63.43] [color={rgb, 255:red, 0; green, 0; blue, 0 }  ][fill={rgb, 255:red, 0; green, 0; blue, 0 }  ][line width=0.75]      (0, 0) circle [x radius= 3.35, y radius= 3.35]   ;
%Straight Lines [id:da6826744329505839] 
\draw    (360,180) -- (320,260) ;
\draw [shift={(360,180)}, rotate = 116.57] [color={rgb, 255:red, 0; green, 0; blue, 0 }  ][fill={rgb, 255:red, 0; green, 0; blue, 0 }  ][line width=0.75]      (0, 0) circle [x radius= 3.35, y radius= 3.35]   ;
%Straight Lines [id:da03857813353289752] 
\draw    (440,180) -- (320,260) ;
\draw [shift={(440,180)}, rotate = 146.31] [color={rgb, 255:red, 0; green, 0; blue, 0 }  ][fill={rgb, 255:red, 0; green, 0; blue, 0 }  ][line width=0.75]      (0, 0) circle [x radius= 3.35, y radius= 3.35]   ;
%Straight Lines [id:da5545998322989484] 
\draw    (520,180) -- (320,260) ;
\draw [shift={(320,260)}, rotate = 158.2] [color={rgb, 255:red, 0; green, 0; blue, 0 }  ][fill={rgb, 255:red, 0; green, 0; blue, 0 }  ][line width=0.75]      (0, 0) circle [x radius= 3.35, y radius= 3.35]   ;
\draw [shift={(520,180)}, rotate = 158.2] [color={rgb, 255:red, 0; green, 0; blue, 0 }  ][fill={rgb, 255:red, 0; green, 0; blue, 0 }  ][line width=0.75]      (0, 0) circle [x radius= 3.35, y radius= 3.35]   ;
%Straight Lines [id:da6987784554627305] 
\draw    (80,100) -- (120,180) ;
%Straight Lines [id:da6309857316302835] 
\draw    (480,100) -- (520,180) ;
\draw [shift={(480,100)}, rotate = 63.43] [color={rgb, 255:red, 0; green, 0; blue, 0 }  ][fill={rgb, 255:red, 0; green, 0; blue, 0 }  ][line width=0.75]      (0, 0) circle [x radius= 3.35, y radius= 3.35]   ;
%Straight Lines [id:da36835555152991306] 
\draw    (240,100) -- (200,180) ;
%Straight Lines [id:da8059504914785415] 
\draw    (320,100) -- (440,180) ;
\draw [shift={(320,100)}, rotate = 33.69] [color={rgb, 255:red, 0; green, 0; blue, 0 }  ][fill={rgb, 255:red, 0; green, 0; blue, 0 }  ][line width=0.75]      (0, 0) circle [x radius= 3.35, y radius= 3.35]   ;
%Straight Lines [id:da4878920044751329] 
\draw    (400,100) -- (360,180) ;
%Straight Lines [id:da9148585897084128] 
\draw    (160,100) -- (120,180) ;
%Straight Lines [id:da7792378354765463] 
\draw    (560,100) -- (520,180) ;
\draw [shift={(560,100)}, rotate = 116.57] [color={rgb, 255:red, 0; green, 0; blue, 0 }  ][fill={rgb, 255:red, 0; green, 0; blue, 0 }  ][line width=0.75]      (0, 0) circle [x radius= 3.35, y radius= 3.35]   ;
%Straight Lines [id:da40926070722273] 
\draw    (80,100) -- (200,180) ;
\draw [shift={(80,100)}, rotate = 33.69] [color={rgb, 255:red, 0; green, 0; blue, 0 }  ][fill={rgb, 255:red, 0; green, 0; blue, 0 }  ][line width=0.75]      (0, 0) circle [x radius= 3.35, y radius= 3.35]   ;
%Straight Lines [id:da1857403557557481] 
\draw    (160,100) -- (360,180) ;
\draw [shift={(160,100)}, rotate = 21.8] [color={rgb, 255:red, 0; green, 0; blue, 0 }  ][fill={rgb, 255:red, 0; green, 0; blue, 0 }  ][line width=0.75]      (0, 0) circle [x radius= 3.35, y radius= 3.35]   ;
%Straight Lines [id:da3644861419541805] 
\draw    (240,100) -- (440,180) ;
\draw [shift={(240,100)}, rotate = 21.8] [color={rgb, 255:red, 0; green, 0; blue, 0 }  ][fill={rgb, 255:red, 0; green, 0; blue, 0 }  ][line width=0.75]      (0, 0) circle [x radius= 3.35, y radius= 3.35]   ;
%Straight Lines [id:da7690826764641078] 
\draw    (400,100) -- (280,180) ;
\draw [shift={(400,100)}, rotate = 146.31] [color={rgb, 255:red, 0; green, 0; blue, 0 }  ][fill={rgb, 255:red, 0; green, 0; blue, 0 }  ][line width=0.75]      (0, 0) circle [x radius= 3.35, y radius= 3.35]   ;
%Straight Lines [id:da45017787726508063] 
\draw    (320,100) -- (120,180) ;
%Straight Lines [id:da35089739488042637] 
\draw    (320,100) -- (520,180) ;
%Straight Lines [id:da98853239239975] 
\draw    (400,100) -- (440,180) ;
%Straight Lines [id:da8141899630569712] 
\draw    (480,100) -- (200,180) ;
%Straight Lines [id:da13157323138111232] 
\draw    (480,100) -- (360,180) ;
%Straight Lines [id:da5255109668440471] 
\draw    (560,100) -- (280,180) ;
%Straight Lines [id:da6453115926914129] 
\draw    (321,19.5) -- (80,100) ;
\draw [shift={(321,19.5)}, rotate = 161.53] [color={rgb, 255:red, 0; green, 0; blue, 0 }  ][fill={rgb, 255:red, 0; green, 0; blue, 0 }  ][line width=0.75]      (0, 0) circle [x radius= 3.35, y radius= 3.35]   ;
%Straight Lines [id:da5283905219019334] 
\draw    (321,19.5) -- (160,100) ;
%Straight Lines [id:da14671861558189936] 
\draw    (321,19.5) -- (240,100) ;
%Straight Lines [id:da31776017040037996] 
\draw    (321,19.5) -- (320,100) ;
%Straight Lines [id:da5684761235925027] 
\draw    (321,19.5) -- (400,100) ;
%Straight Lines [id:da7384832191617999] 
\draw    (321,19.5) -- (480,100) ;
%Straight Lines [id:da7766595023695629] 
\draw    (321,19.5) -- (560,100) ;
%Straight Lines [id:da4202518772497972] 
\draw    (80,100) -- (280,180) ;

% Text Node
\draw (75,180) node [anchor=north west][inner sep=0.75pt]   [align=left] {\textcolor[rgb]{0.82,0.01,0.11}{$\displaystyle 12|3|4$}};
% Text Node
\draw (235,180) node [anchor=north west][inner sep=0.75pt]   [align=left] {$\displaystyle \textcolor[rgb]{0.82,0.01,0.11}{13|2|4}$};
% Text Node
\draw (158,180) node [anchor=north west][inner sep=0.75pt]   [align=left] {\textcolor[rgb]{0.82,0.01,0.11}{$\displaystyle 23|1|4$}};
% Text Node
\draw (443,180) node [anchor=north west][inner sep=0.75pt]   [align=left] {$ $$\displaystyle \textcolor[rgb]{0.82,0.01,0.11}{14|2|3}$};
% Text Node
\draw (529,180) node [anchor=north west][inner sep=0.75pt]   [align=left] {$\displaystyle \textcolor[rgb]{0.82,0.01,0.11}{24|1|3}$};
% Text Node
\draw (369,180) node [anchor=north west][inner sep=0.75pt]   [align=left] {$\displaystyle \textcolor[rgb]{0.82,0.01,0.11}{34|1|2}$};
% Text Node
\draw (40,95) node [anchor=north west][inner sep=0.75pt]  [color={rgb, 255:red, 255; green, 255; blue, 255 }  ,opacity=1 ] [align=left] {$\displaystyle \textcolor[rgb]{0.82,0.01,0.11}{123|4}$};
% Text Node
\draw (120,95) node [anchor=north west][inner sep=0.75pt]   [align=left] {$\displaystyle 12|34$};
% Text Node
\draw (200,95) node [anchor=north west][inner sep=0.75pt]   [align=left] {$\displaystyle 13|24$};
% Text Node
\draw (275,90) node [anchor=north west][inner sep=0.75pt]   [align=left] {$\displaystyle \textcolor[rgb]{0.82,0.01,0.11}{124|3}$};
% Text Node
\draw (405,95) node [anchor=north west][inner sep=0.75pt]   [align=left] {$\displaystyle \textcolor[rgb]{0.82,0.01,0.11}{134|2}$};
% Text Node
\draw (490,95) node [anchor=north west][inner sep=0.75pt]   [align=left] {$\displaystyle \textcolor[rgb]{0.82,0.01,0.11}{234|1}$};
% Text Node
\draw (565,95) node [anchor=north west][inner sep=0.75pt]   [align=left] {$\displaystyle 13|24$};
% Text Node
\draw (340,10) node [anchor=north west][inner sep=0.75pt]   [align=left] {$\displaystyle \textcolor[rgb]{0.82,0.01,0.11}{1234}$};
% Text Node
\draw (336,259) node [anchor=north west][inner sep=0.75pt]   [align=left] {$\displaystyle {1|2|3|4}$};
\end{tikzpicture}}
    \end{minipage}  
    \caption{The lattice of flats $\LL(A_3) \cong \Pi_4$. The flats in the minimal building set $\BB^{\textrm{min}}$ are drawn in red.}\label{figure:braidlattice}
     \end{figure}
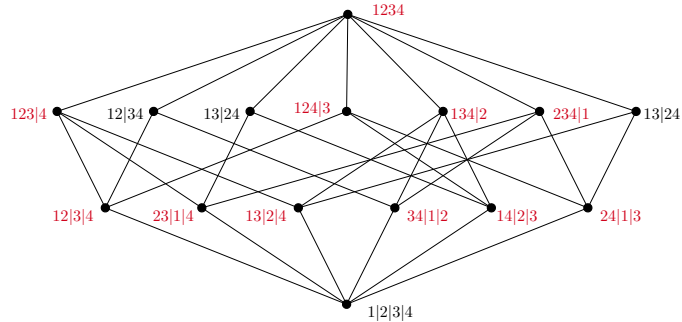

We will repurpose notation and refer interchangeably to $W \in \LL(A_{n-1})$ and the corresponding set partition in $\Pi_n$. In this language, the minimal building set is
\[ \BB^{\textrm{min}} = \bigg\{ \{ j_1, j_2, \cdots, j_k \ | \ i_1 \ | \cdots | \ i_\ell \}: k \geq 2\bigg\}. \] Figure \ref{figure:braidlattice} shows $\LL(A_3) \cong \Pi_4$, with the elements of $\BB^{\textrm{min}}$ colored in red.

The boundary components of $\overline{M}_{0,n+1}$ are well-studied. In particular, there are $2^n - n -2$ divisors, each of which factors as a product
$\overline{M}_{0,k} \times \overline{M}_{0,n+3-k}$ for some $3\leq k \leq n$. In the language of Section~\ref{section:wonderful-boundaries}, this says that for any $\PP W \in \BB^{\textrm{min}},$ both $\PP W^{(\BB^{\textrm{min}})_W}$ and $\PP(V/W)^{(\BB^{\textrm{min}})^W}$ are isomorphic to copies of $\overline{M}_{0,k}$ for some $3\leq k \leq n$. From the perspective of hyperplane arrangements, this is equivalent to the statement that the restriction and contraction arrangements of $A_{n-1}$ by $W \in \BB^{\textrm{min}}$ are always isomorphic to a smaller braid arrangement.

\subsubsection{Wondertopes and the Parke-Taylor form}
There are $n!/2$ regions in $C(A_{n-1}) \cong M_{0,n+1}$. 
The symmetric group acts transitively on these regions by permuting the standard basis vectors of $\RR^n$, which gives an automorphism of $\PP V$ preserving the hyperplanes in $A_{n-1}$. Thus, without loss of generality, we need only consider a single region in either space. 

Each region in $M_{0,n+1}$ corresponds to an ordering of the points $z_0,z_1,\dots,z_{n}$ on $\PP^1(\RR)$, up to orientation.
As before, fix $(z_0,z_{n-1},z_{n}) = (0,1,\infty)$ and consider the region %$(M_{0,n+1})_{>0}$ corresponding to the ordering
\begin{equation} \label{eq:posregion-M0n} 
(M_{0,n+1})_{>0} = \{ 0 < z_1 < z_2 < \cdots < z_{n-2} < 1\} . 
\end{equation}
To see this region as a simplex, define new coordinates $x_i$ for $1 \leq i \leq n-2$ by:
\[ x_i := \begin{cases}
    z_1 & \text{if } i = 1\\
    z_{i} - z_{i-1} & \text{if } 2 \leq i \leq n-2.
\end{cases}  \]
Under this linear change of coordinates, $M_{0,n+1}$ becomes
{\small \begin{align*}
    M_{0,n+1} &= \big \{ (z_1,z_2, \dots, z_{n-2}) \in \RR^{n-2}: z_i \neq 0,1 \text{ and } z_i \neq z_j \text{ for } i \neq j \big \}\\
    &= \big\{(x_1,x_2, \dots, x_{n-2}) \in \RR^{n-2} : x_1+x_2+\cdots+x_i \neq 1 \text{ and } x_{i} + x_{i+1}+\cdots+x_j \neq 0 \text{ for all } i\leq j\big\}.
\end{align*}}
Then the region in \eqref{eq:posregion-M0n} is the interior of a simplex in $\RR^{n-2}$, defined via the following inequalities:
$$(M_{0,n+1})_{>0} = \big\{ (x_1,x_2, \dots, x_{n-2}) \in \RR^{n-2} : x_i > 0 \text{ for all } i \text{ and } x_1+x_2+\cdots + x_{n-2} < 1 \big\}.$$

Considering $\RR^{n-2}$ as an affine chart of $\PP^{n-2}(\RR)$, the Euclidean closure $(M_{0,n+1})_{\geq 0}$ of $(M_{0,n+1})_{>0}$ becomes a simplex in $\PP^{n-2}(\RR)$. According to Proposition~\ref{prop:canonical-form-polytope} \ref{prop:canonical-form-general-simplex}, $\big(\PP^{n-2}, (M_{0,n+1})_{\geq 0}\big)$ is a positive geometry with canonical form given (up to sign) by
\begin{equation}\label{eq:canonicalform-simplex} \Omega\big(\PP^{n-2}, (M_{0,n+1})_{\geq 0}\big) = \frac{dx_1 dx_2 \cdots dx_{n-2}}{x_1\cdots x_{n-2} (1-x_1-\cdots -x_{n-2})}.\end{equation}

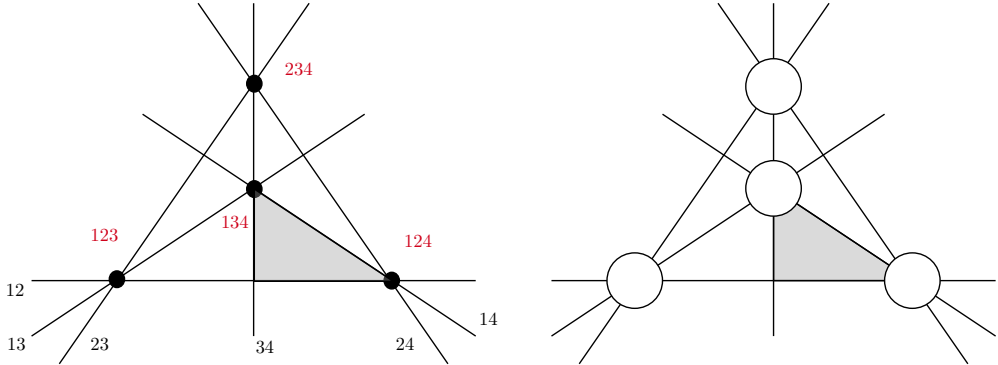
\begin{figure}[!h]
	\begin{minipage}{0.45\textwidth}
 \centering
    \scalebox{.7}{\input{./BEPV_Figures_arXiv/m05beforelabelled.tex}}
    \end{minipage}
	\begin{minipage}{0.45\textwidth}
	\centering
 \scalebox{.7}
    {\input{./BEPV_Figures_arXiv/m05after.tex}}
    \end{minipage}
    \caption{A region (a simplex) in $M_{0,5} \cong C(A_3)$ in $\PP^2$ (left) and its wondertope $(\overline{M}_{0,5})_{\geq 0}$ (an associahedron) in $\overline{M}_{0,5}$ (right)}
    \label{fig:n=5example}
\end{figure}

Recall that $\overline{M}_{0,n+1}$ is the wonderful compactification of $M_{0,n+1}$ with respect to the minimal building set $\build^{\textrm{min}}$ described in Section~\ref{section:minimal-building-set}, and let $(\overline{M}_{0,n+1})_{\geq 0}$ be the strict transform of $(M_{0,n+1})_{\geq 0}$ under this compactification.
Figure~\ref{fig:n=5example} shows an affine-local picture of $C(A_3) \cong M_{0,5}$ and $\overline{M}_{0,5}$, together with the regions $(M_{0,5})_{\geq 0}$ and $(\overline{M}_{0,5})_{\geq 0}$. (The singletons are dropped from the partition notation for simplicity.)

By Theorem \ref{thm:intro_mainthm}, 
$\big(\overline{M}_{0,n+1}, (\overline{M}_{0,n+1})_{\geq 0}\big)$ is a positive geometry with canonical form
$$\Omega\big(\overline{M}_{0,n+1}, (\overline{M}_{0,n+1})_{\geq 0} \big) = \pi^* \, \Omega\big(\PP^{n-2}, (M_{0,n+1})_{\geq 0}\big),$$
where $\pi$ denotes the blow-down map, as usual.
Therefore, on an appropriate local chart where $\pi$ is  the identity (i.e. away from the exceptional locus), the canonical form $\Omega\big(\overline{M}_{0,n+1}, (\overline{M}_{0,n+1})_{\geq 0}\big)$ is given by \eqref{eq:canonicalform-simplex}. Substituting back the original $z_i$ coordinates gives 
\begin{equation}\label{eq:parketaylorform}\Omega\big(\overline{M}_{0,n+1}, (\overline{M}_{0,n+1})_{\geq 0}\big) = \frac{dz_1 \ dz_2 \cdots dz_{n-2} }{(z_1-z_0)(z_2 - z_1)(z_3-z_2) \cdots (z_{n-1}-z_{n-2})},\end{equation}
where as before, $(z_0, z_{n-1}, z_{n}) = (0,1,\infty)$.
The form in \eqref{eq:parketaylorform} is known in the physics literature as the \emph{Parke-Taylor form} \cite[\S 8.2]{AHLcluster}.
The Parke-Taylor form can also be written in terms of \emph{dihedral} coordinates introduced by Brown in \cite{brown2009multiple}; see  \cite[\S 2]{delpezzo} for a detailed example of $\overline{M}_{0,5}$.

It turns out that $(\overline{M}_{0,n+1})_{\geq 0}$ is homeomorphic, as a stratified space, to the face stratification of the $n$-\emph{associahedron} \cite{AHLcluster}, a convex $(n-2)$-dimensional polytope whose vertices are in bijection with triangulations of an $(n+1)$-gon. The associahedron plays an important role in algebraic combinatorics, topology, and representation theory; an interested reader should look at \cite{postnikov2009permutohedra} and the references therein. For example, the wondertope $(\overline{M}_{0,5})_{\geq 0}$ in Figure~\ref{fig:n=5example} is combinatorially isomorphic to the $4$-associahedron (a pentagon), whose vertices are in bijection with the triangulations of a $5$-gon (also in this case a pentagon).

\subsection{Non-examples}\label{section:nonexamples}

We collect a series of examples where a blow-up of a positive geometry is not a positive geometry, illustrating the necessity of the conditions in Theorem~\ref{thm:maingeneral}.
To this end, we first note that $(\PP^1, \PP^1(\RR))$ is not a positive geometry.  Indeed, the nonnegative part $\PP^1(\RR)$ has no boundary, so that the canonical form has no poles, but $\PP^1$ has no nonzero $1$-form with no poles.
The pair $(\PP^1, \PP^1(\RR))$ is instead a \emph{pseudo-positive geometry} as in \cite[Section 2.2]{NimaLam17}.

\begin{eg}\label{eg:triangle}
We first give an example where we blow up a point in the boundary of a polytope which is not itself a face. Consider a triangle $T$ in $\PP^2(\RR)$ and a point $p$ in the interior of an edge.  Let $\pi: \operatorname{Bl}_p\ \PP^2 \to \PP^2$ be the blow-up.

Note that $T$ contains all normal directions (up to sign) at $p.$
Thus, the Euclidean closure of $\pi^{-1}(\intr(T))$, denoted $\widetilde T$, contains the entirety of the exceptional divisor $E(\RR)$.  See Figure \ref{fig:triangle} for an illustration.
In particular, we find that $E_{\geq 0}$, the Euclidean closure in $E(\RR)$ of the interior of $\widetilde T \cap E$, is all of $E(\RR) \cong \PP^1(\RR)$, so that the boundary component $(E,E_{\geq 0}) = (\PP^1, \PP^1(\RR))$ is not a positive geometry.  The pair $(\operatorname{Bl}_p \PP^2, \widetilde T)$ is thus not a positive geometry.

\begin{figure}[!h]
	\begin{minipage}{0.4\textwidth}
 \centering
    \scalebox{0.7}{\input{./BEPV_Figures_arXiv/trianglebefore.tex}}
    \end{minipage}
	\begin{minipage}{0.5\textwidth}
	\centering
    \scalebox{0.6}{\input{./BEPV_Figures_arXiv/triangleafter.tex}}
    \end{minipage}
    \caption{A triangle with point $p$ (left) and its blow-up (right)}
    \label{fig:triangle}
\end{figure}
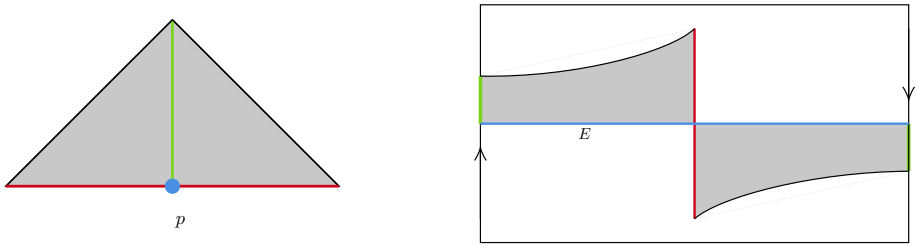
\end{eg}

The next two examples illustrate that Theorem~\ref{thm:maingeneral} does not generalize easily to the setting where the initial positive geometry is $(\PP^n, \poly)$ for $\poly$ a polyhedral complex instead of a polytope.

\begin{eg}\label{eg:complex}
Consider the positive geometry consisting of the three triangles in $\PP^2(\RR)$ illustrated in Figure~\ref{fig:trianglecomplex}. Let us consider blowing up the point $p$. 
As in Example \ref{eg:triangle}, every normal direction (up to sign) to $p$ is present in at least one polygon, so that once again we have $E_{\geq 0 } = E(\RR) \cong \PP^1(\RR)$ and so $(E,E_{\geq 0})$ is not a positive geometry.

In this example, every log resolution of our positive geometry is \emph{not} a positive geometry.
Indeed, let $D$ be the union of the three lines through $p$, which is not a simple normal crossing divisor.  By \cite[Proposition V.5.2]{hartshorne1977algebraic}, any log resolution of the pair $(\PP^2, D)$ must factor through the blow-up at $p$, which we have shown not to be a positive geometry.

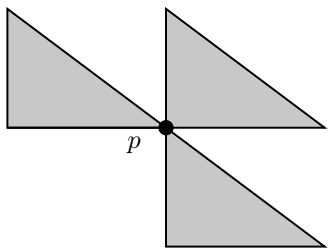
\begin{figure}[!h]
    \input{./BEPV_Figures_arXiv/trianglecomplex.tex}
    \caption{Three triangles meeting at a point}
    \label{fig:trianglecomplex}
\end{figure}

\end{eg}

\begin{eg}\label{eg:staircase}
Even if the polyhedral complex is connected in codimension $1,$ the resulting blow-up may not be a positive geometry. For example, consider four cubes arranged in a ``staircase'' and meeting at a point $p$, as in Figure \ref{fig:staircase}.
The boundary components of this polyhedral complex have a simple normal crossing at $p$, but the blow-up of $p$ is not a positive geometry.
Indeed, because the polyhedral complex contains all normal directions to $p$, we again find $E_{\geq 0} = E(\RR)$. Moreover, in this case, $E(\RR) \cong \PP^2(\RR)$ is not orientable, so we do not even get a pseudo-positive geometry.

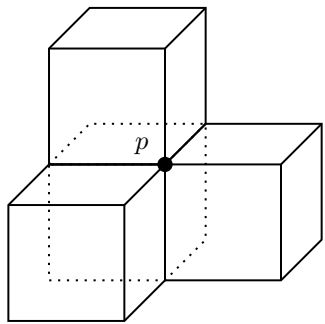
\begin{figure}[!h]
    \input{./BEPV_Figures_arXiv/staircase}
    \caption{Four cubes meeting at a point}
    \label{fig:staircase}
\end{figure}
\end{eg}

\begin{rem}
    Examples~\ref{eg:complex} and \ref{eg:staircase} hint at the fact that, if one hopes to prove a more general statement about a positive geometry $(X,X_{\geq 0})$ remaining a positive geometry after blowing-up a boundary stratum, some ``convexity condition'' on the semialgebraic set $X_{\geq 0}$ is required. We caution that the notion of ``corners'' in \cite[Definition~2.16]{Brown-Dupont} does not capture this necessary behavior in our setting as, for instance, the point $p$ in Example~\ref{eg:staircase} is a corner.
\end{rem}

\subsection{A detailed example} \label{section:detailed-example}

We feature one example in detail, which will illustrate behaviors such as:
\begin{enumerate}[label=(\roman*)]
\item \label{item:behavior-i} Without the condition that a collection $\build$ includes all facet hyperplanes of $\poly$, the boundary components of the pair $(\PP V^\build, \widetilde\poly^\build)$ may not be simple normal crossing, even though the pair is a positive geometry.
\item \label{item:behavior-ii} As cautioned above Lemma~\ref{lemma:facet-transform}, for a facet hyperplane $\PP H$ of $\poly$, the intersection $\widetilde\poly^\build \cap \widetilde{\PP H}$ may not be equal to $\widetilde{\poly_H}^{\build_H}$.
\item \label{item:behavior-iii} Without the condition that a collection $\build$ is well-adapted to a polytope $\poly$, even if $\build$ is a building set, the pair $(\PP V^\build, \widetilde\poly^\build)$ may not be a positive geometry.
\end{enumerate}
All of our examples will be variations on a square pyramid $\poly$ in $\PP^3(\RR).$

\begin{eg}\label{eg:pyramid1}
This example will illustrate the behaviors \ref{item:behavior-i} and \ref{item:behavior-ii} mentioned above.
Let $\poly$ be a square pyramid in $\PP^3(\RR)$.
Let $p$ be the cone point of $\poly$. Consider the blow-up $\text{Bl}_{Y} \PP^3,$ where $Y$ is the line through points $p$ and $q$ as labelled in Figure \ref{fig:pyramid}. Note that the collection $\{Y\}$ is trivially a building set which is well-adapted to $\poly$.  However, this collection does not contain the facet hyperplanes of $\poly$.
In this case, the pair $(\Bl_Y \PP^3, \widetilde\poly)$ is a positive geometry by Theorem~\ref{thm:fundcomp}, but we claim that its boundary components do not form a simple normal crossing divisor, and that the face poset of $\widetilde\poly$ is not isomorphic to the face poset of any polytope.

Let $F_1, \dotsc, F_4$ be the four triangles of the pyramid $\poly$, with $F_1, F_2$ in the back and $F_3, F_4$ in the front.
The exceptional divisor $E$ of the blow-up is isomorphic to $Y \times \PP^1$.
Using Proposition~\ref{prop:exceptionalproduct} to compute the strict transforms $\widetilde\poly, \widetilde{F_1}, \dotsc, \widetilde{F_4}$ of the polytopes, we find that they are as pictured in the right side of Figure \ref{fig:pyramid}.
The nonnegative region $E_{\geq 0}$ is the quadrilateral in the back with the colored boundary.

 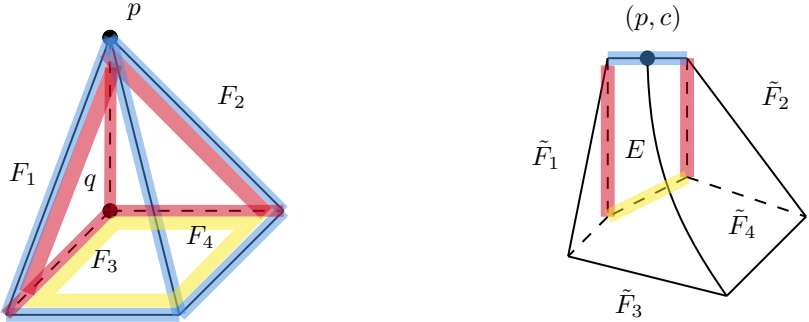
\begin{figure}[!h]
	\begin{minipage}{0.45\textwidth}
 \centering
    {\input{./BEPV_Figures_arXiv/pyramidbefore.tex}}
    \end{minipage}
	\begin{minipage}{0.45\textwidth}
	\centering
    {\input{./BEPV_Figures_arXiv/pyramidafter.tex}}
    \end{minipage}
    \caption{A square pyramid with cone point $p$ (left) and its blow-up (right)}
    \label{fig:pyramid}
	\end{figure}

Let $H_i$ be the hyperplane containing $F_i$, and $\widetilde{H_i}$ its strict transform in the blow-up.
Under the isomorphism $E \cong Y \times \PP^1$, we have $\widetilde{H_3} \cap E = \widetilde{H_4} \cap E = \{p\} \times \PP^1$, since $H_3$ and $H_4$ are hyperplanes both intersecting $Y$ only at $p$.
When restricted to $\widetilde\poly$, these intersections are pictured in blue in Figure \ref{fig:pyramid}, which is a union of boundaries of $\widetilde{F_3}$ and $\widetilde{F_4}$.
In particular, the positive geometry of $\widetilde\poly$ has three 2-dimensional boundary components $\widetilde{H_3}$, $\widetilde{H_4}$, and $E$ whose common intersection is the 1-dimensional $\{p\}\times \PP^1$.
Such an intersection is not simple normal crossing, and such behavior does not  happen in a 3-dimensional polytope: no three facets of any 3-dimensional polytope meet in a common line.

Lastly, we see that $\widetilde{\poly} \cap \widetilde{H_4}(\RR)$ is the region colored purple in Figure~\ref{fig:extrasegment}, which contains an additional segment in $E_{\geq 0}$ when compared to $\widetilde{\poly_{H_4}}$, colored green, illustrating the behavior~\ref{item:behavior-ii}.

\begin{figure}[!h]
	\begin{minipage}{0.45\textwidth}
 \centering
    {\input{./BEPV_Figures_arXiv/segmentno.tex}}
    \end{minipage}
	\begin{minipage}{0.45\textwidth}
	\centering
    {\input{./BEPV_Figures_arXiv/segmentyes.tex}}
    \end{minipage}
    \caption{The semialgebraic subsets  $\widetilde{\poly_{H_4}}$ (left) and $\widetilde{\poly} \cap \widetilde{H_4}(\RR)$ (right)}
    \label{fig:extrasegment}
	\end{figure}
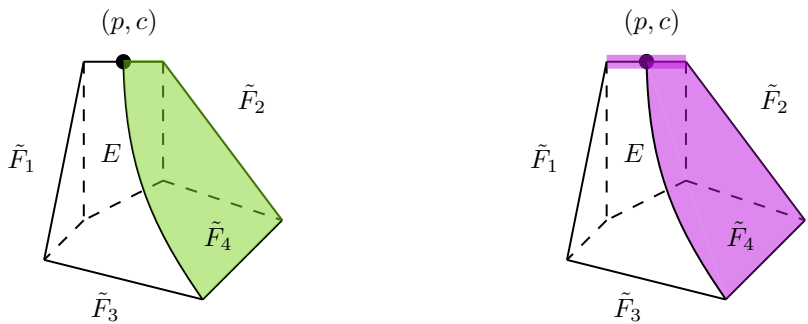
\end{eg}

 \begin{eg}\label{eg:pyramid2}
Using the notation from Example \ref{eg:pyramid1}, let us consider the collection $\{Y, Y'\}$ where $Y'$ is the line obtained by the intersection $H_3 \cap H_4$, i.e.\ the ``opposite line'' to the line $Y$ in the pyramid.
Note that this collection is not a building set, since the intersection $Y\cap Y' = \{p\}$ has codimension 
$3$ in $\PP^3$, whereas $Y$ and $Y'$ each have codimension $2$.
In the second blow-up of the sequential blow-up $\operatorname{Bl}_{Y'}\operatorname{Bl}_Y\PP^3$, the boundary component $(E, E_{\geq 0})$ of $(\operatorname{Bl}_Y\PP^3, \widetilde\poly)$ is blown-up at the point labelled $(p,c)$, which is not a positive geometry for the same reason that the blow-up in Example~\ref{eg:triangle} was not (i.e. because we are blowing-up a point in the interior of an edge).
In particular, the sequential blow-up is not a positive geometry.
 \end{eg}

 \begin{eg}\label{eg:pyramid3}
 We now give an example illustrating~\ref{item:behavior-iii}.  In other words, we show that $\build$ being a building set is insufficient, and the well-adapted condition on $\build$ is necessary for the blow-up to be a positive geometry. 

Consider a four-dimensional polytope $\poly$, which is a pyramid over a triangular prism.  Its set of vertices and a Schlegel diagram are given in Figure \ref{fig:schlegel}.

     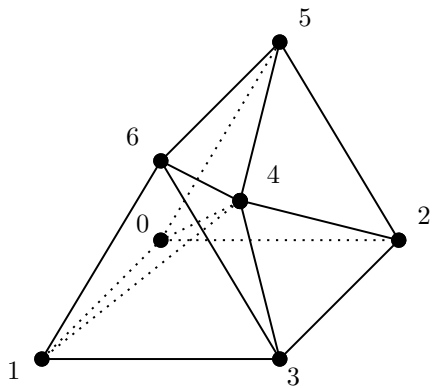
\begin{figure}[!h]
         \centering
         \input{./BEPV_Figures_arXiv/schlegel}
         \caption{A Schlegel diagram of $\poly$. Vertex $4$ is the top of the pyramid. 
         The ordered vertices can be given coordinates $\{0000,1000,0100,1100,0010,0001,1001\}$ in $\RR^4 \subset \PP^4(\RR).$}
         \label{fig:schlegel}
     \end{figure}
    Let $F_1, F_2,$ and $H$ be the linear spans of the vertices (considered as vectors in $\RR^4$ via $\RR^4 \cong \{X_0 \neq 0\} \subset \PP^4(\RR)$) labelled by $\{2,4,5\},\{1,4,6\},$ and $\{0,1,2,3,4\},$ respectively. Note that $\poly \cap \PP H$ recovers the square pyramid.  Consider the set $\build = \{F_1, F_2\},$ which is a building set; indeed, the codimensions of $F_1$ and $F_2$ add to the codimension of $F_1 \cap F_2.$ However, $\build_H = \{F_1 \cap \PP H, F_2 \cap \PP H\}$ recovers the two opposite lines of Example \ref{eg:pyramid2}, which we have seen is not a building set in $\PP H$, and does not result in a positive geometry.
 \end{eg}

\subsection{A further example} 

We provide one non-polytopal example where the blow-up of a boundary stratum results in a positive geometry.

\begin{eg}\label{eg:exc-not-product}
    Consider the semialgebraic set
    $$S := \left\{ [X:Y:Z:W] \in \PP^3(\RR) \mid W>0 \; \text{ and } \; Z, W-Z, W-Y, Y-X^2 \geq 0 \right\} \subset \PP^3(\RR).$$
    
    Consider the affine chart $U:= \{W \neq 0\} \subset \PP^3$ with dehomogenized coordinates $x = \frac{X}{W}, y = \frac{Y}{W},$ and $z = \frac{Z}{W}$. Our set $S$ is a full-dimensional semialgebraic subset of $U(\RR) \cong \AA^3(\RR)$, cut out by the following inequalities: 
    $$0\leq z \leq 1 \; \text{ and } \; x^2 \leq y \leq 1.$$

    The pair $(\PP^3, S)$ is a positive geometry with canonical form $\Omega = \Omega(\PP^3, S)$ given by
    $$\Omega\big|_{U} = \frac{2}{(y-x^2)(y-1)z(z-1)} dx \, dy\,dz.$$

    The subvariety $Y:=V(Z,YW - X^2) \subset \PP^3$ is the intersection of the boundary components $\{Z = 0\}$ and $\{YW = X^2\}$ of $(\PP^3, S)$. Figure~\ref{fig:exc-not-product} depicts the two boundary components and their intersection in our affine chart $U = \{W \neq 0\}$. The blow-up $\Bl_Y \PP^3$ has exceptional divisor
    $E = \PP\big(N_{Y\mid \PP^3}\big) \cong \Bl_{\text{pt}}\PP^2,$
    which is not isomorphic to a product of two projective spaces. Computations in affine-local coordinates show that $\left(\Bl_Y\PP^3, \widetilde{S}\right)$ is a positive geometry with canonical form $\pi^* \Omega$.

    \begin{figure}[!h]
	\begin{minipage}{0.45\textwidth}
 \centering
    \scalebox{.7}{\includegraphics[scale = .4]{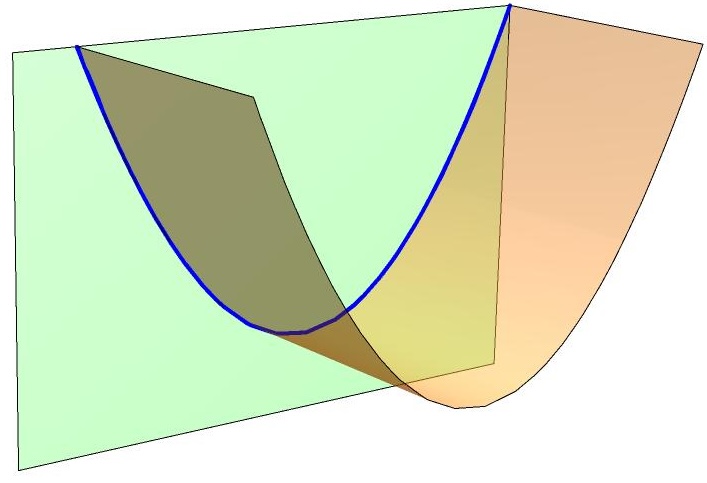}}
    \end{minipage}
	\begin{minipage}{0.45\textwidth}
	\centering
 \scalebox{.7}
    {\includegraphics[scale = .35]{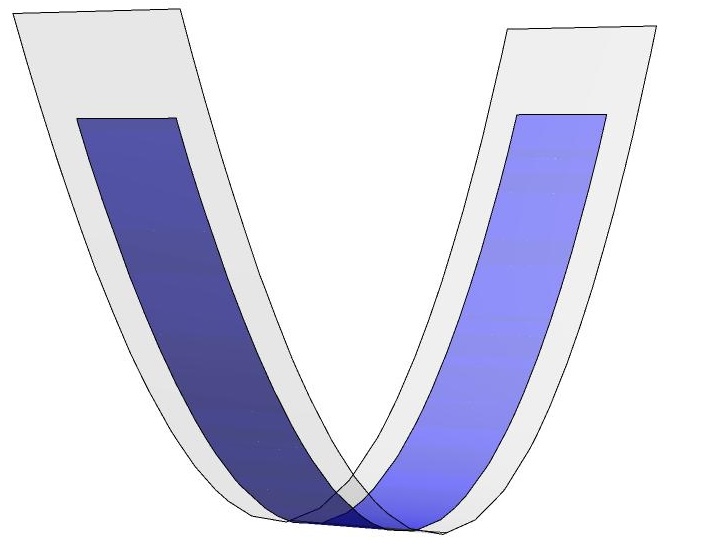}}
    \end{minipage}
    \caption{Left: The boundary components $\{z= 0\}$ in green and $\{y = x^2\}$ in orange, and their intersection $Y = V(z,y-x^2)$ in blue.
    Right: The exceptional divisor $E$ of $\Bl_{Y} \PP^3$ in grey and its non-negative part $E_{\geq 0}$ in blue.}
    \label{fig:exc-not-product}
\end{figure}
\end{eg}

\begin{rem}
    Example~\ref{eg:exc-not-product} indicates a possible future direction of study: one can investigate whether blow-ups of (higher dimensional) \emph{polypols} (see \cite{polypols}) are positive geometries. 
\end{rem}

\appendix
\section{Summary of combinatorial analogues}\label{section:combinatorialproduct}

We review combinatorial counterparts in matroid theory of the geometry of wonderful compactifications of hyperplane arrangements, and prove Theorem \ref{thm:nestedproduct}, the combinatorial counterpart of Proposition \ref{prop:wonderful}.
A more general version of Theorem \ref{thm:nestedproduct} can be found in \cite[Propositions 2.8.6-7]{bibby2022leray}. However, it takes some work to specialize their results to our context. We thus include our argument for completeness, with a view towards future directions (e.g. Question~\ref{question:matroidalwondertope}).

To motivate the combinatorial definitions to follow, we show in Table \ref{table:geometrytocombinatorics} the combinatorial objects of interest in relation to the geometry discussed in the body of the paper. The remainder of the section is then devoted to defining and discussing the right-hand column of  Table \ref{table:geometrytocombinatorics}.
For a reference on matroid theory, we point to \cite{Wel76, Oxl11}.

\begin{center}
\begin{table}[!h]
\begin{tabular}{|c|l|}\hline
{ \bf Geometry} & { \bf Combinatorics} \\ \hline\hline
Hyperplane arrangement $\arrangement$ & Matroid $M$ \\ \hline 
Lattice $\LL(\arrangement)$ of intersection subspaces of $\arrangement$ &  Lattice $\LL(M)$ of flats of $M$ \\ \hline
Building set $\BB$ as in the Introduction&  Combinatorial building set $\BB$  \\ 
& as in Definition \ref{def:buildingset_matroid}\\
\hline
Intersection lattice of exceptional divisors & Face lattice of $\NN(\BB)$  as in \\
in the wonderful compactification $X^\BB$& Definition \ref{def:nestedsetcomplex} \cite{feichtner2004incidence}\\ \hline
(Codimension $k$) intersections of $k$ exceptional divisors &  $(k-1)$-dimensional faces in $\NN(\BB)$ \\\hline
$E_W$, (the strict transform of) the exceptional divisor & Link of the nested set complex  \\
of the blow-up at $\PP W$ for $W \in \BB$ & $\NN(\BB)_W$ for $W \in \BB$ \\ 
\hline
Spaces $\PP (V/W)$ and $\PP W $& Restriction and contraction matroids \\
with building sets $\BB^W$ and $\BB_W$&$M^W$ and $M_W$ with building sets\\
& $\BB^W$ and $\BB_W$ as in Definition \ref{buildsetrestrictcontract}\\ \hline
Proposition \ref{prop:wonderful} \cite[Theorem 4.3]{de1995wonderful}:  & Theorem \ref{thm:nestedproduct}:\\
product structure of $E_W$& product structure of $\NN(\BB)_W$\\
\hline
Wondertope and its canonical form & {\bf Unknown}: Question \ref{question:matroidalwondertope} \\\hline
\end{tabular}
\vspace{0.3cm}
\caption{The geometric objects discussed throughout the paper, and their analogues in the context of matroid theory.}\label{table:geometrytocombinatorics}
\vspace{-1cm}
\end{table}
\end{center}

\begin{defn}
    A \emph{matroid} $M = (E, \mathcal{I})$ consists of a finite set $E$ and a nonempty collection $\mathcal{I}$ of subsets of $E$ such that 
    \begin{enumerate}
        \item if $I \in \mathcal{I}$ and $J \subset I$ then $J \in \mathcal{I}$ as well, and 
        \item if $I, J \in \mathcal{I}$ and $|I| < |J|$, there is a $j \in J \setminus I $ so that $I \cup \{ j \} \in \mathcal{I}$. 
    \end{enumerate}
\end{defn}

The set $E$ is called the \emph{ground set} of $M$ and the elements $\mathcal{I}$ are the \emph{independent sets} of $M$. For any subset $X \subseteq E$, the \emph{rank} of $X$ is 
    \[ \mathrm{rk}_M(X) = \max\{ |I|: I \in \mathcal{I} \text{ and } I \subseteq X \}.\]
A subset $F \subseteq E$ is a \emph{flat} of $M$ if $ \mathrm{rk}_M(F)< \mathrm{rk}_M(F \cup \{ i \})$ for any $i \in E \setminus F$.

The set of all flats under inclusion forms a lattice $\LL(M)$ called the \emph{lattice of flats} of $M$. For any $F, F'$ in $\LL(M)$, their \emph{join} $F \vee F'$ is the smallest flat containing $F \cup F'$ and their \emph{meet} is $F \wedge F' = F \cap F'$. We write the maximal element of $\LL(M)$ as $\hat{1}$, the minimal element as $\hat{0}$, and incomparable elements $X,Y \in \LL(M)$ as $X \  ||  \ Y$.

Note that when the matroid $M$ is realizable as a subspace $V \subset \CC^{E}$ (i.e. when $E$ can be identified with a finite set of vectors in $V^\vee$) then $M$
defines a hyperplane arrangement $\arrangement$, 
where $i \in E$ corresponds to $i$-th coordinate hyperplane of $\CC^{E}$, intersected with $V$. In this case, the lattices $\LL(M)$ and the lattice of subspaces $\LL(\arrangement)$ under intersection (ordered by reverse inclusion) are isomorphic.

 In addition to $\LL(M)$, we will be interested in the \emph{intervals} in $\LL(M)$ generated by $X, Y \in \LL(M)$:
\[ [X,Y]:= \{ F \in \LL(M): X \leq F \leq Y \} \subset \LL(M). \]
For any subset of elements $\BB  \subset \LL(M)$, write 
\[ \BB_{\leq X}:= \{ F \in \BB: F \leq X \}, \]
and let $\max(\BB_{\leq X})$ be the set of maximal elements of $\BB_{\leq X}$. Note that if $X \in \BB$, then $\max(\BB_{\leq X}) = X$.

Given a flat $F \in \LL(M)$, one can define the \emph{restriction} and \emph{contraction} matroids $M^ F$ and $M_F$ on the ground sets $F$ and $E \setminus F$, respectively. These matroids have flats
\begin{align*}
    \LL(M^ {F}):=& \{ X: X \text{ is a flat of } M \text{ contained in } F \}, \\
    \LL(M_ F) :=& \{ X \setminus F: X \text{ is a flat of }M \text{ which contains } F\}.  
\end{align*}
and in fact there are poset isomorphisms $\LL(M^ {F}) \cong [\hat{0}, F]$ and $\LL(M_ F) \cong [F, \hat{1}].$ We will implicitly go back and forth along these isomorphisms. 

A key insight of De Concini--Procesi in  \cite{de1995wonderful} was that the irreducible components of the normal crossing divisor in a wonderful compactification of a hyperplane complement $C(\arrangement) = \PP V \setminus \PP \arrangement$ could be indexed by a subset of the lattice of intersections $\LL(\arrangement)$ called a \emph{building set}. This idea was then generalized by Feichtner--Kozlov \cite{feichtner2004incidence} to the context of an arbitrary semilattice. In what follows, we restrict our attention to the case of geometric lattices, i.e. the lattice of flats for some matroid $M$.
\begin{defn}[Building set] \label{def:buildingset_matroid}
    A set $\BB  \subset \LL(M) \setminus \{\hat{0}\}$ is a \emph{building set} of $\LL(M)$ if for any $X \in \LL(M) \setminus \{ \hat{0} \}$ with $\max(\BB_{\leq X}) = \{ X_1, \cdots X_\ell \}$, there is an isomorphism of posets 
    \begin{align*}
        \varphi_X:& \prod_{j=1}^{\ell} [ \hat{0}, X_j  ] \longrightarrow [ \hat{0}, X  ]
    \end{align*}
     where $\varphi_X$ sends $( \hat{0}, \cdots, X_j, \cdots \hat{0}  )$ to $X_j$ for each $j$. The set $\max(\BB_{\leq X})$ is called the \emph{set of factors} of $X$.
    \end{defn}
Observe that if $X \in \BB$, then $\varphi_X$ is the identity. The \emph{minimal building set} $\BB^{\textrm{min}}$ is the building set of $\LL(M)$ which includes $F$ if and only if $[\hat{0},F]$ cannot be written as a product of lower intervals in $\LL(M)$. The \emph{maximal building set} $\BB^{\max}$ is the building set which includes every element of $\LL(M) \setminus~\{ \hat{0} \}$. 
\begin{eg}[Boolean matroid]
   The \emph{Boolean matroid} on a set $S$ has flats consisting of all subsets of $S$. In this case, Definition \ref{def:buildingset_matroid} coincides with the definition of building set in \cite[Definition 7.1]{postnikov2009permutohedra}, which states that $\BB$ is a building set if (1) for $I, J \in \BB$, if $I \cap J \neq \emptyset$, then $I \cup J \in \BB$, and (2)  $\{ i \} \in \BB$ for every $i \in S$. The minimal building set $\BB^{\textrm{min}}$ consists of the atoms of $\LL(M)$, i.e. the singleton subsets of $S$. 
\end{eg}

From the data of a building set $\BB$, one can define a simplicial complex whose elements are the so-called \emph{nested sets}, defined below. Let $\BB^{\mathrm{top}}$ be the set of maximal elements of $\BB$; following Zelevinsky \cite{zelevinsky2006nested}, we will not include $\BB^{\mathrm{top}}$ in this simplicial complex. 

\begin{defn}[Nested set complex]\label{def:nestedsetcomplex}
For a building set $\BB$, a subset $S$ of  $\BB \setminus \BB^{\mathrm{top}}$ is \emph{nested} if for any set of incomparable elements $X_1, \dots, X_t$ in $S$ with $t \geq 2$, their join $X_1 \vee X_2 \vee \cdots \vee X_t$ is \emph{not} in $\BB$. We define the \emph{nested set complex} $\NN(\BB)$ of $\BB$ to be the simplicial complex on $\BB \setminus \BB^{\mathrm{top}}$ whose elements are the nested sets of $\BB$. 
\end{defn}

For instance, the vertex set (i.e. the set of 0-simplices) of $\NN(\BB)$ is $\BB \setminus \BB^{\textrm{top}}$.
Note that, in particular, $\{ X_1, X_2 \}$ is nested if $X_1 \leq X_2$; indeed, we emphasize that the join condition above is only applicable in the case of incomparable elements. Thus the 1-simplices of $\NN(\BB)$ are pairs $\{ X_1, X_2 \}$ such that one of the following holds: 
\begin{enumerate}
    \item $X_1 < X_2$,
    \item $X_2 < X_1$, or
    \item  $X_1 \vee X_2 \not \in \BB$.
\end{enumerate}

Given a flat $F \in \BB$, define the \emph{link} of $F$ in $\NN(\BB)$ as
\[ \NN(\BB)_F:= \{ I \in \NN(\BB): F \not \in I \text{ and } F \cup I \in \NN(\BB) \}. \]

Our goal is to relate $\NN(\BB)_F$ to nested sets corresponding to certain building sets in the lattices of flats of $M^ {F}$ and $M_ F$. In order to do so, we must first define the appropriate building sets.
\begin{defn}\label{buildsetrestrictcontract}
    Define the \emph{restriction of $\BB$ along $F$} as the set 
\[ \BB^ F:= \{ X \in \BB : \ X \leq F \}, \]
and the \emph{contraction of $\BB$ to $F$} as
\[ \BB_ F:= \{ X \vee F:  X \in \BB \text{ and }  X \not \leq F\}. \]
\end{defn}
We claim that $\BB^ F$ and $\BB_ F$ actually satisfy the conditions in Definition \ref{def:buildingset_matroid} for $\LL(M^ F)$ and $\LL(M_F)$, respectively. To verify this, we will first need Lemmas \ref{lemma:partitioncontraction} and \ref{lemma:uniqueness_join}.

\begin{lem}\label{lemma:partitioncontraction}
    The set $\BB_ F$ can be partitioned into the sets 
    \begin{align*}
        (\BB_ F)_1 &:= \{ X \in \BB: F \lneq X \}\\
        (\BB_ F)_2 & := \{ X \vee F: X \in \BB, \ X \  || \ F, \text{ and }  X \vee F \not \in \BB\} 
    \end{align*}
\end{lem}
\begin{proof}
For any $X \in \BB$ such that $X \not \leq F$, let $Y = X \vee F$. If $F \neq X$, then $Y = X \in (\BB_ F)_1$. If $F || X$ and $Y \in \BB$, then $F \lneq Y$ and so $Y  \in (\BB_ F)_1$. Otherwise, $Y \not \in \BB$ and so $Y \in (\BB_ F)_2$. Thus, every element of the form $X \vee F$ for $X \not \leq F$ and $X \in \BB$ is in exactly one of our two sets $(\BB_ F)_1$ and $(\BB_ F)_2$.
\end{proof}

\begin{lem}\label{lemma:uniqueness_join}
For a fixed $F \in \BB$, suppose there exist $X, X' \in \BB$ which are both incomparable to $F$, and for which $X \vee F = X' \vee F \not \in \BB$. Then $X = X'$. 
\end{lem}
\begin{proof}
Note that $\{ X, F \}$ is a nested set. By \cite[Proposition 2.8 (1),(2)]{feichtner2004incidence}, the set of factors of $X \vee F$ is precisely $\{ X, F \}$. Since $X \vee F = X' \vee F$, we have $X' < X \vee F$; Proposition 2.5 (1) of \cite{feichtner2004incidence} then says that there is a unique element $Y$ of the set of factors of $X \vee F$ such that $X' \leq Y$. By the incomparability assumption, $Y \neq F$ so we must have $Y = X$, and thus $X' \leq X$. Reversing the argument then shows that $X \leq X'$, and so they must be equal.
\end{proof}
Lemma \ref{lemma:uniqueness_join} thus shows that for any $Y = X \vee F \in (\BB_ F)_2$, the element $X$ is unique. We are now ready to show that $\BB^ F$ and $\BB_ F$ are indeed building sets.
Fix a building set $\BB$ of $\LL(M)$ and $F \in \BB$.

\begin{prop}\label{prop:building-sets}
    The sets $\BB^ F$ and $\BB_ F$ are building sets of $\LL(M^ F)$ and $\LL(M_ F)$, respectively. 
\end{prop}
\begin{proof}
Suppose $X \in \LL(M)$ and $X \leq F$, so $X\in [\hat{0}, F]$. Then 
\[ [\hat{0}, X] \cong \prod_{i=1}^k \ [\hat{0}, Y_i] \]
for $Y_i \in \max(\BB_{\leq X})$. Hence $Y_i \leq F$ as well, and so $Y_i \in \BB^ F$. Thus $\BB^ F$ is a building set for $\LL(M^ F)$.

Now consider the case that $X > F$, so that $X \in [F,\hat{1}]$. Again, assume that $\max(\BB_{\leq X}) = \{Y_1,\dots, Y_k\}$, so that there exists a poset isomorphism
	$\varphi: [\hat{0}, X] \xrightarrow{\cong} \prod_{i=1}^k \; [\hat{0}, Y_i]$ with $\varphi(Y_i) = (\hat{0},\dots,\hat{0},Y_i, \hat{0}, \dots,\hat{0})$. We will show that $[F, X] \cong \prod_{i=1}^k \; [F, Y_i \vee F]$.

Since $F < X$ and $F \in B$, there must exist some $i$ such that $F \leq Y_i$. This $i$ is unique.
Indeed, if $F \leq Y_i$ and $F \leq Y_j$ for $i \neq j$, then $\varphi(F) \leq \varphi(Y_i)$ and $\varphi(F) \leq \varphi(Y_j )$. This is
impossible because the only element smaller than both $\varphi(Y_i)$ and $\varphi(Y_j )$ is $(\hat{0}, . . . , \hat{0})$. Without loss of generality, assume $F \leq Y_1$ and $\varphi(F) = (F, \hat{0}, \dots, \hat{0}).$ In particular, this means that $\varphi$ restricts to an isomorphism
	$$[F,X] \cong [F, Y_1] \times \prod_{i > 1} \; [\hat{0}, Y_i].$$
	
	We now show that $[\hat{0}, Y_i] \cong [F, F \vee Y_i]$ for all $i > 1$, which proves our claim. Since $\varphi(F) = (F, \hat{0}, \dots, \hat{0})$ and $\varphi(Y_i) = (\hat{0},\dots,Y_i, \dots,\hat{0})$, then $\varphi(F\vee Y_i) = (F, \hat{0}, \dots, Y_i, \dots, \hat{0})$ because this is the lowest common upper bound of $\varphi(F)$ and $\varphi(Y_i)$ in $\prod_i \; [\hat{0}, Y_i]$. And so $\varphi$ restricts to an isomorphism 
	$$[F, F\vee Y_i] \; \cong \; \{F\} \times \{\hat{0}\} \times \dots \times [\hat{0}, Y_i] \times \dots \times \{\hat{0}\} \; \cong \; [\hat{0}, Y_i].$$
	This isomorphism sends $F \mapsto \hat{0}$ and $F \vee Y_i \mapsto Y_i$, so we get an isomorphism 
	$$[F, X] \cong \prod_i \; [F, Y_i \vee F]$$
	which sends each $Y_i \vee F$ to $(F, \dots, F, Y_i \vee F, F, \dots, F)$, as desired.

	Finally, the fact above that $\varphi(F\vee Y_i) = (F, \hat{0}, \dots, Y_i, \dots, \hat{0})$ proves that the $F \vee Y_i$'s are incomparable in $\LL(M)$ because so are their images under $\varphi$.
\end{proof}

\begin{eg}
Recall from Section \ref{section:braid} the description of $\LL(A_{n-1}) \cong \Pi_n$ and the minimal building set 
\[ \BB^{\textrm{min}} = \bigg \{\{ j_{1}, j_2 \cdots j_k |i_{1}| \cdots |i_{\ell} \}: k \geq 2 \bigg \}.\]
Consider the case of $\LL(A_5) \cong \Pi_6$ with building set $\BB = \BB^{\textrm{min}}$ and flat $F = \{ 123|4|5|6\}$. Then 
\[ \BB^{\{ 123|4|5|6\}} = \bigg \{ \{ 12|3|4|5|6 \} ,  \{ 13|2|4|5|6 \} ,  \{ 23|1|4|5|6 \},  \{ 123|4|5|6\} \bigg \}.   \]
Using the notation of Lemma \ref{lemma:partitioncontraction}, 
\begin{align*}
    \big( \BB_{\{ 123|4|5|6\}}\big)_1 =& \bigg \{ \{ 1234|5|6 \} ,  \{ 1235|4|6 \} ,  \{ 1236|4|5 \},  \{ 12345|6 \}, \{ 12346|5 \}, \{ 12356|4 \} , \{123456\} \bigg \},\\
    \big( \BB_{\{ 123|4|5|6\}} \big)_{2} =& \bigg \{ \{ 123|45|6 \} ,  \{ 123|4|56 \} ,  \{ 123|46|5 \},  \{ 123|456 \}  \bigg \},
\end{align*}  
and $\BB_{\{ 123|4|5|6\}} = \big( \BB_{\{ 123|4|5|6\}}\big)_1 \sqcup \big( \BB_{\{ 123|4|5|6\}}\big)_2$.
\end{eg}
	
We can now describe the relationship between the link of $\NN(\BB)$ and the sets $\BB^ F$ and $\BB_ F$.
\begin{thm}[Theorem \ref{thm:nestedproduct}] \label{thm:nestedsetiso}
    For any $F \in \BB$, there are isomorphisms of simplicial complexes
    \[ \NN(\BB)_F \cong \NN(\BB^ F) \times \NN(\BB_ F) \cong \NN(\BB^ F \times \BB_ F).\]
\end{thm}
Theorem \ref{thm:nestedsetiso} was proved by Zelevinksy in \cite{zelevinsky2006nested} when $M$ is the Boolean matroid; our proof is modeled on his argument.
\begin{proof}
For a fixed $F \in \BB$, write $\BB' = \BB^ F \times \BB_ F$ and $\BB_0$ as the vertex set of $\NN(\BB)_F$. It follows that $\BB_0$ consists of the elements 
\[ \BB_0 = \{ X \in \BB \setminus \{ F \}: \{ X, F \} \in \NN(\BB) \}, \]
which can be described as the disjoint union of 
\begin{align*}
    \BB_1:= & \{ X \in \BB \setminus \{ F, \hat{1} \}: F < X \}, \\
    \BB_2:= & \{ X \in \BB \setminus \{\hat{1}\} :  X \ || \  F, \ \  X \vee F \not \in \BB \}\\
    \BB_3:= & \{ X \in \BB \setminus \{ F \}: X < F \}.
\end{align*}
Given $X \in \BB_0$, define the map 
\[ \tau: X \mapsto \begin{cases} X & X \in \BB_1, \BB_3\\
X \vee F & X \in \BB_2.
\end{cases} \]
We claim that $\tau$ defines a bijection between the vertex sets of $\NN(\BB_0)$ and $\NN(\BB')$. First, 
 $\BB_3 = \BB^ F \setminus \{ F \}$ is the vertex set of $\NN(\BB^ F)$ since $F = (\BB^ F)^{\textrm{top}}$ and is therefore not included in the vertex set of $\NN(\BB^ F)$. 
Analogously, since $(\BB_ F)^{\textrm{top}} = \hat{1}$, we have that $\BB_1$ corresponds to the portion of the vertex set of $\NN(\BB_F)$ coming from $\BB_1$. Finally,
by Lemma \ref{lemma:uniqueness_join}, $\tau$ is a bijection between elements of $\BB_2$ and $(\BB_ F)_2 \setminus \{ \hat{1} \}$. Thus, $\tau$ defines a bijection between the vertex sets of $\NN(\BB_0)$ and $\NN(\BB')$. Write the inverse of this map as $\tau^{-1}$.

It remains to show that in fact, $\tau$ defines an isomorphism of simplicial complexes; in other words $I = \{ X_1, \cdots, X_k \} \in \NN(\BB)_F$ (equivalently $\{ X_1, \cdots X_k, F \} \in \NN(\BB)$) if and only if $\tau(I):= \{ \tau(X_1), \cdots, \tau(X_k) \} \in \NN(\BB')$. We will show the converse, namely that $I \not \in \NN(\BB)_F$ if and only if $\tau(I) \not \in \NN(\BB')$.

Consider $J = \{ Y_1 \cdots Y_k \} \not \in \NN(\BB')$, with $Y_i \in \BB'$.  Without loss of generality, we may assume that $Y_i \  || \  Y_j$ for all $i,j$, and thus $Y:= Y_1 \vee \cdots \vee Y_k \in \BB'$. If $Y \in \BB^ F$, then $\tau^{-1}(Y_i) = Y_i$ for all $i$, and so $\{ J, F \} \not \in \NN(\BB)$, since $Y \in \BB^ {F} \subset \BB$.

Suppose $Y \in \BB_ F$; by Lemma \ref{lemma:partitioncontraction}, one has $Y$ in either $(\BB_ F)_1$ or $(\BB_ F)_2$. We have the following cases:
\begin{itemize}
    \item If every $Y_i \in (\BB_ F)_1$, then it follows that $Y \in (\BB_ F)_1$,    so
    $\tau^{-1}(Y_i) = Y_i$, and again we have that $\{ J, F \} \not \in \NN(\BB)$ since $Y$ must be in $\BB_3 \subset \BB$. 
    \item If every $Y_i \in (\BB_ F)_2$, then $Y_i = X_i \vee F$ for $X_i \in \BB$ incomparable to $F$, and so $Y$ can be written as 
    \[ Y = (X_1 \vee F) \vee \cdots (X_k \vee F) = (X_1 \vee \cdots \vee X_k) \vee F. \]
    Since $(X_1 \vee \cdots \vee X_k) \  ||  \ F$, we have that $Y \in (\BB_ F)_2$. By Lemma \ref{lemma:partitioncontraction}, $(X_1 \vee \cdots \vee X_k) \in \BB$, and so $\{ X_1 \cdots X_k, F \}$ cannot be in $\NN(\BB)$.
    \item Suppose we have $Y_i$ in both $(\BB_ F)_1$ and $(\BB_ F)_2$, so that $Y$ can be written (up to relabeling) as 
    \begin{align*}
        Y = (X_1 \vee F) \vee \cdots \vee (X_\ell \vee F) \vee Y_{\ell + 1} \vee \cdots \vee Y_{k} &=
        (X_1 \vee \cdots \vee X_\ell) \vee F \vee  (Y_{\ell +1} \vee \cdots \vee Y_k)\\
        &=  (X_1 \vee \cdots \vee X_\ell) \vee  (Y_{\ell +1} \vee \cdots \vee Y_k)
    \end{align*} 
    where the $F$ disappears because $Y_j \vee F = Y_j$ for any $\ell + 1 \leq j \leq k$. Then because each $Y_j$ properly contains $F$, it must be that $Y \in (\BB_ F)_1$ since $Y$ cannot be written as $X \vee F$ as in Lemma \ref{lemma:partitioncontraction}. Thus for $\tau^{-1}(J) = \{ X_1 \cdots X_{\ell}, Y_{\ell +1}, \cdots, X_k \}$, we have 
    \[ \tau^{-1}(Y) =  (X_1 \vee \cdots \vee X_\ell) \vee (Y_{\ell +1} \vee \cdots \vee Y_k) = Y \in \BB_3 \subset \BB, \]
    and so again $\{ \tau^{-1}(J) \cup F \}$ cannot be nested. 
\end{itemize}
An analogous argument in the opposite direction then gives the claim.
\end{proof}

%bibliography
\small
\bibliography{literature}
\bibliographystyle{alpha}

\end{document}

%% file: BEPV_Figures_arXiv/prodbefore.tex
\tikzset{every picture/.style={line width=0.75pt}} %set default line width to 0.75pt        

\begin{tikzpicture}[x=0.75pt,y=0.75pt,yscale=-1,xscale=1]
%uncomment if require: \path (0,300); %set diagram left start at 0, and has height of 300

%Shape: Cube [id:dp09691646959291422] 
\draw   (260,99.5) -- (298.5,61) -- (400,61) -- (400,162.5) -- (361.5,201) -- (260,201) -- cycle ; \draw   (400,61) -- (361.5,99.5) -- (260,99.5) ; \draw   (361.5,99.5) -- (361.5,201) ;
%Straight Lines [id:da08480837794183949] 
\draw [color={rgb, 255:red, 126; green, 211; blue, 33 }  ,draw opacity=0.6 ][line width=4.5]    (260,81) -- (260,221) ;
%Shape: Rectangle [id:dp21285636428163657] 
\draw  [dash pattern={on 0.84pt off 2.51pt}] (300,61) -- (400,61) -- (400,161) -- (300,161) -- cycle ;
%Straight Lines [id:da05273092415937364] 
\draw  [dash pattern={on 0.84pt off 2.51pt}]  (300,161) -- (260,201) ;

\end{tikzpicture}

%% file: BEPV_Figures_arXiv/prodafter.tex
\tikzset{every picture/.style={line width=0.75pt}} %set default line width to 0.75pt        

\begin{tikzpicture}[x=0.75pt,y=0.75pt,yscale=-1,xscale=1]
%uncomment if require: \path (0,300); %set diagram left start at 0, and has height of 300

%Flowchart: Stored Data [id:dp7269199625853422] 
\draw  [fill={rgb, 255:red, 126; green, 211; blue, 33 }  ,fill opacity=0.5 ] (292,101.5) -- (292,200.73) .. controls (292,190.29) and (282.02,179.56) .. (269.71,176.75) .. controls (257.41,173.95) and (247.43,180.14) .. (247.43,190.57) -- (247.43,91.34) .. controls (247.43,80.9) and (257.41,74.71) .. (269.71,77.52) .. controls (282.02,80.32) and (292,91.06) .. (292,101.5) -- cycle ;
%Straight Lines [id:da3537784055600085] 
\draw    (300,61.72) -- (400,61.72) ;
%Straight Lines [id:da3264029746162196] 
\draw    (400,61.72) -- (400,161.72) ;
%Straight Lines [id:da330226688603406] 
\draw    (300,61.72) -- (290.81,70.81) -- (279.01,81.57) ;
%Straight Lines [id:da6459318668283519] 
\draw    (400,161.72) -- (360,201.72) ;
%Straight Lines [id:da05758485547477599] 
\draw    (400,61.72) -- (360,101.72) ;
%Straight Lines [id:da17838474716326524] 
\draw    (292,101.5) -- (360,101.72) ;
%Straight Lines [id:da9548544297578031] 
\draw    (360,101.72) -- (360,201.72) ;
%Straight Lines [id:da9885990516772846] 
\draw    (292,200.73) -- (360,201.72) ;
%Shape: Rectangle [id:dp5209181261038505] 
\draw  [dash pattern={on 0.84pt off 2.51pt}] (300,61.72) -- (400,61.72) -- (400,161.72) -- (300,161.72) -- cycle ;
%Straight Lines [id:da8025611599456548] 
\draw  [dash pattern={on 0.84pt off 2.51pt}]  (300,161.72) -- (279,180.5) ;
%Shape: Rectangle [id:dp6598751821105362] 
\draw  [draw opacity=0][fill={rgb, 255:red, 255; green, 255; blue, 255 }  ,fill opacity=1 ] (278.49,71.84) -- (278.49,210.33) -- (209.42,210.33) -- (209.42,71.84) -- cycle ;
%Curve Lines [id:da16285805687357346] 
\draw [color={rgb, 255:red, 74; green, 74; blue, 74 }  ,draw opacity=1 ] [dash pattern={on 0.84pt off 2.51pt}]  (254,101.5) .. controls (265,94.5) and (282,101.5) .. (292,119.5) ;
%Shape: Arc [id:dp419254555926731] 
\draw  [draw opacity=0] (253.62,179.09) .. controls (257.27,176.83) and (262.25,175.81) .. (267.79,176.55) .. controls (271.86,177.09) and (275.71,178.5) .. (279,180.5) -- (268.51,193.22) -- cycle ; \draw   (253.62,179.09) .. controls (257.27,176.83) and (262.25,175.81) .. (267.79,176.55) .. controls (271.86,177.09) and (275.71,178.5) .. (279,180.5) ;  
%Straight Lines [id:da2669731725707075] 
\draw    (252.39,79.36) -- (253.62,179.09) ;
%Shape: Arc [id:dp014329260958881762] 
\draw  [draw opacity=0] (252.39,79.36) .. controls (256.05,77.11) and (261.02,76.09) .. (266.57,76.82) .. controls (271.16,77.43) and (275.46,79.15) .. (279.01,81.57) -- (267.28,93.5) -- cycle ; \draw   (252.39,79.36) .. controls (256.05,77.11) and (261.02,76.09) .. (266.57,76.82) .. controls (271.16,77.43) and (275.46,79.15) .. (279.01,81.57) ;  
%Straight Lines [id:da022450610909868107] 
\draw [color={rgb, 255:red, 155; green, 155; blue, 155 }  ,draw opacity=1 ]   (279.01,81.57) -- (279,180.5) ;
%Curve Lines [id:da14333170380570803] 
\draw [color={rgb, 255:red, 74; green, 74; blue, 74 }  ,draw opacity=1 ] [dash pattern={on 0.84pt off 2.51pt}]  (253.62,130.09) .. controls (264.62,123.09) and (281.62,130.09) .. (291.62,148.09) ;

\end{tikzpicture}

%% file: BEPV_Figures_arXiv/m05beforelabelled.tex
\tikzset{every picture/.style={line width=0.75pt}} %set default line width to 0.75pt        

\begin{tikzpicture}[x=0.75pt,y=0.75pt,yscale=-1,xscale=1]
%uncomment if require: \path (0,300); %set diagram left start at 0, and has height of 300

%Straight Lines [id:da2595689200505573] 
\draw    (340,0) -- (160,260) ;
%Straight Lines [id:da20209357591905164] 
\draw    (260,0) -- (440,260) ;
%Straight Lines [id:da12227075149732447] 
\draw    (140,200) -- (460,200) ;
%Straight Lines [id:da12667411046467913] 
\draw    (300,0) -- (300,240) ;
%Straight Lines [id:da8386518145076799] 
\draw    (460,240) -- (220,80) ;
%Straight Lines [id:da2006804366522824] 
\draw    (380,80) -- (140,240) ;
%Flowchart: Connector [id:dp9187294539922067] 
\draw  [fill={rgb, 255:red, 0; green, 0; blue, 0 }  ,fill opacity=1 ] (196,199) .. controls (196,195.69) and (198.46,193) .. (201.5,193) .. controls (204.54,193) and (207,195.69) .. (207,199) .. controls (207,202.31) and (204.54,205) .. (201.5,205) .. controls (198.46,205) and (196,202.31) .. (196,199) -- cycle ;
%Flowchart: Connector [id:dp2294884600842544] 
\draw  [fill={rgb, 255:red, 0; green, 0; blue, 0 }  ,fill opacity=1 ] (394,200) .. controls (394,196.69) and (396.46,194) .. (399.5,194) .. controls (402.54,194) and (405,196.69) .. (405,200) .. controls (405,203.31) and (402.54,206) .. (399.5,206) .. controls (396.46,206) and (394,203.31) .. (394,200) -- cycle ;
%Flowchart: Connector [id:dp5621013357713643] 
\draw  [fill={rgb, 255:red, 0; green, 0; blue, 0 }  ,fill opacity=1 ] (295,58) .. controls (295,54.69) and (297.46,52) .. (300.5,52) .. controls (303.54,52) and (306,54.69) .. (306,58) .. controls (306,61.31) and (303.54,64) .. (300.5,64) .. controls (297.46,64) and (295,61.31) .. (295,58) -- cycle ;
%Flowchart: Connector [id:dp5965027221576544] 
\draw  [fill={rgb, 255:red, 0; green, 0; blue, 0 }  ,fill opacity=1 ] (295,134) .. controls (295,130.69) and (297.46,128) .. (300.5,128) .. controls (303.54,128) and (306,130.69) .. (306,134) .. controls (306,137.31) and (303.54,140) .. (300.5,140) .. controls (297.46,140) and (295,137.31) .. (295,134) -- cycle ;
%Shape: Right Triangle [id:dp8773757465425586] 
\draw  [fill={rgb, 255:red, 155; green, 155; blue, 155 }  ,fill opacity=0.37 ] (300.5,134) -- (400.5,200.5) -- (300.5,200.5) -- cycle ;

% Text Node
\draw (120,201) node [anchor=north west][inner sep=0.75pt]  [color={rgb, 255:red, 0; green, 0; blue, 0 }  ,opacity=1 ] [align=left] {$\displaystyle 12$};
% Text Node
\draw (121,240) node [anchor=north west][inner sep=0.75pt]  [color={rgb, 255:red, 0; green, 0; blue, 0 }  ,opacity=1 ] [align=left] {$\displaystyle 13$};
% Text Node
\draw (181,240) node [anchor=north west][inner sep=0.75pt]   [align=left] {$\displaystyle 23$};
% Text Node
\draw (461,222) node [anchor=north west][inner sep=0.75pt]   [align=left] {$\displaystyle 14$};
% Text Node
\draw (401,240) node [anchor=north west][inner sep=0.75pt]   [align=left] {$\displaystyle 24$};
% Text Node
\draw (181,161) node [anchor=north west][inner sep=0.75pt]  [color={rgb, 255:red, 208; green, 2; blue, 27 }  ,opacity=1 ] [align=left] {$\displaystyle 123$};
% Text Node
\draw (300,242) node [anchor=north west][inner sep=0.75pt]   [align=left] {$\displaystyle 34$};
% Text Node
\draw (407,166) node [anchor=north west][inner sep=0.75pt]  [color={rgb, 255:red, 208; green, 2; blue, 27 }  ,opacity=1 ] [align=left] {$\displaystyle 124$};
% Text Node
\draw (321,42) node [anchor=north west][inner sep=0.75pt]  [color={rgb, 255:red, 208; green, 2; blue, 27 }  ,opacity=1 ] [align=left] {$\displaystyle 234$};
% Text Node
\draw (275,152) node [anchor=north west][inner sep=0.75pt]  [color={rgb, 255:red, 208; green, 2; blue, 27 }  ,opacity=1 ] [align=left] {$\displaystyle 134$};

\end{tikzpicture}

%% file: BEPV_Figures_arXiv/m05after.tex
\tikzset{every picture/.style={line width=0.75pt}} 
\begin{tikzpicture}[x=0.75pt,y=0.75pt,yscale=-1,xscale=1]
%uncomment if require: \path (0,300); %set diagram left start at 0, and has height of 300

%Shape: Right Triangle [id:dp812420414095946] 
\draw  [fill={rgb, 255:red, 155; green, 155; blue, 155 }  ,fill opacity=0.37 ] (340,133.5) -- (440,200) -- (340,200) -- cycle ;
%Straight Lines [id:da26538011226539804] 
\draw    (380,0) -- (200,260) ;
%Straight Lines [id:da027923421114589897] 
\draw    (300,0) -- (391.03,131.48) -- (480,260) ;
%Straight Lines [id:da9655702540598985] 
\draw    (180,200) -- (500,200) ;
%Straight Lines [id:da29645066446593815] 
\draw    (340,0) -- (340,240) ;
%Straight Lines [id:da8758246422989321] 
\draw    (500,240) -- (260,80) ;
%Straight Lines [id:da919786210124061] 
\draw    (420,80) -- (180,240) ;
%Flowchart: Connector [id:dp7894051329213467] 
\draw  [fill={rgb, 255:red, 255; green, 255; blue, 255 }  ,fill opacity=1 ] (320,60) .. controls (320,48.95) and (328.95,40) .. (340,40) .. controls (351.05,40) and (360,48.95) .. (360,60) .. controls (360,71.05) and (351.05,80) .. (340,80) .. controls (328.95,80) and (320,71.05) .. (320,60) -- cycle ;
%Flowchart: Connector [id:dp5329741848287216] 
\draw  [fill={rgb, 255:red, 255; green, 255; blue, 255 }  ,fill opacity=1 ] (420,200) .. controls (420,188.95) and (428.95,180) .. (440,180) .. controls (451.05,180) and (460,188.95) .. (460,200) .. controls (460,211.05) and (451.05,220) .. (440,220) .. controls (428.95,220) and (420,211.05) .. (420,200) -- cycle ;
%Flowchart: Connector [id:dp21857012746323357] 
\draw  [fill={rgb, 255:red, 255; green, 255; blue, 255 }  ,fill opacity=1 ] (220,200) .. controls (220,188.95) and (228.95,180) .. (240,180) .. controls (251.05,180) and (260,188.95) .. (260,200) .. controls (260,211.05) and (251.05,220) .. (240,220) .. controls (228.95,220) and (220,211.05) .. (220,200) -- cycle ;
%Flowchart: Connector [id:dp7060659875846711] 
\draw  [fill={rgb, 255:red, 255; green, 255; blue, 255 }  ,fill opacity=1 ] (320,133.5) .. controls (320,122.45) and (328.95,113.5) .. (340,113.5) .. controls (351.05,113.5) and (360,122.45) .. (360,133.5) .. controls (360,144.55) and (351.05,153.5) .. (340,153.5) .. controls (328.95,153.5) and (320,144.55) .. (320,133.5) -- cycle ;

\end{tikzpicture}

%% file: BEPV_Figures_arXiv/trianglebefore.tex
\tikzset{every picture/.style={line width=0.75pt}} %set default line width to 0.75pt        

\begin{tikzpicture}[x=0.75pt,y=0.75pt,yscale=-1,xscale=1]
%uncomment if require: \path (0,300); %set diagram left start at 0, and has height of 300

%Shape: Triangle [id:dp9468642048723857] 
\draw  [fill={rgb, 255:red, 200; green, 200; blue, 200 }  ,fill opacity=1 ] (340,100) -- (460,220) -- (220,220) -- cycle ;
%Straight Lines [id:da29155449647885423] 
\draw [color={rgb, 255:red, 126; green, 211; blue, 33 }  ,draw opacity=1 ][line width=1.5]    (340,100) -- (340,220) ;
%Straight Lines [id:da8775149672294112] 
\draw    (340,100) -- (220,220) ;
%Straight Lines [id:da14870906316549615] 
\draw    (460,220) -- (340,100) ;
%Straight Lines [id:da9284322900280968] 
\draw [color={rgb, 255:red, 208; green, 2; blue, 27 }  ,draw opacity=1 ][line width=1.5]    (220,220) -- (460,220) ;
%Shape: Circle [id:dp7607255835798434] 
\draw  [draw opacity=0][fill={rgb, 255:red, 74; green, 144; blue, 226 }  ,fill opacity=1 ] (334.5,220) .. controls (334.5,216.96) and (336.96,214.5) .. (340,214.5) .. controls (343.04,214.5) and (345.5,216.96) .. (345.5,220) .. controls (345.5,223.04) and (343.04,225.5) .. (340,225.5) .. controls (336.96,225.5) and (334.5,223.04) .. (334.5,220) -- cycle ;

% Text Node
\draw (341,240) node [anchor=north west][inner sep=0.75pt]   [align=left] {$\displaystyle p$};

\end{tikzpicture}

%% file: BEPV_Figures_arXiv/triangleafter.tex
\tikzset{every picture/.style={line width=0.75pt}} %set default line width to 0.75pt        

\begin{tikzpicture}[x=0.75pt,y=0.75pt,yscale=-1,xscale=1]
%uncomment if require: \path (0,300); %set diagram left start at 0, and has height of 300

%Shape: Polygon [id:ds46295092390327464] 
\draw  [draw opacity=0][fill={rgb, 255:red, 200; green, 200; blue, 200 }  ,fill opacity=1 ] (520,140) -- (520,180) -- (340,220) -- (340,140) -- cycle ;
%Shape: Polygon [id:ds019659870862811957] 
\draw  [draw opacity=0][fill={rgb, 255:red, 200; green, 200; blue, 200 }  ,fill opacity=1 ] (340,60) -- (340,140) -- (160,140) -- (160,100) -- cycle ;
%Straight Lines [id:da08048060396169798] 
\draw [color={rgb, 255:red, 126; green, 211; blue, 33 }  ,draw opacity=1 ][line width=2.25]    (520,140) -- (520,180) ;
%Straight Lines [id:da8872005895998365] 
\draw [color={rgb, 255:red, 208; green, 2; blue, 27 }  ,draw opacity=1 ][line width=1.5]    (340,220) -- (340,60) ;
%Shape: Rectangle [id:dp13206390720625272] 
\draw   (160,40) -- (520,40) -- (520,240) -- (160,240) -- cycle ;
%Straight Lines [id:da8139155778103764] 
\draw [color={rgb, 255:red, 126; green, 211; blue, 33 }  ,draw opacity=1 ][line width=2.25]    (160,100) -- (160,113.5) -- (160,140) ;
%Straight Lines [id:da7123060149463483] 
\draw [color={rgb, 255:red, 74; green, 144; blue, 226 }  ,draw opacity=1 ][line width=1.5]    (160,140) -- (520,140) ;
%Straight Lines [id:da7719664170543487] 
\draw    (160,220) -- (160,162) ;
\draw [shift={(160,160)}, rotate = 90] [color={rgb, 255:red, 0; green, 0; blue, 0 }  ][line width=0.75]    (10.93,-4.9) .. controls (6.95,-2.3) and (3.31,-0.67) .. (0,0) .. controls (3.31,0.67) and (6.95,2.3) .. (10.93,4.9)   ;
%Straight Lines [id:da13780653253873598] 
\draw    (520,60) -- (520,118) ;
\draw [shift={(520,120)}, rotate = 270] [color={rgb, 255:red, 0; green, 0; blue, 0 }  ][line width=0.75]    (10.93,-4.9) .. controls (6.95,-2.3) and (3.31,-0.67) .. (0,0) .. controls (3.31,0.67) and (6.95,2.3) .. (10.93,4.9)   ;
%Curve Lines [id:da5990797764802824] 
\draw [fill={rgb, 255:red, 255; green, 255; blue, 255 }  ,fill opacity=1 ]   (160,100) .. controls (240,101.5) and (320,79.5) .. (340,60) ;
%Curve Lines [id:da6674771783418083] 
\draw [fill={rgb, 255:red, 255; green, 255; blue, 255 }  ,fill opacity=1 ]   (340,220) .. controls (361,201.5) and (440,180.5) .. (520,180) ;

% Text Node
\draw (241,142) node [anchor=north west][inner sep=0.75pt]   [align=left] {$\displaystyle E$};

\end{tikzpicture}

%% file: BEPV_Figures_arXiv/trianglecomplex.tex
\tikzset{every picture/.style={line width=0.75pt}} %set default line width to 0.75pt        

\begin{tikzpicture}[x=0.75pt,y=0.75pt,yscale=-1,xscale=1]
%uncomment if require: \path (0,300); %set diagram left start at 0, and has height of 300

%Shape: Right Triangle [id:dp19767903126884712] 
\draw  [fill={rgb, 255:red, 200; green, 200; blue, 200 }  ,fill opacity=1 ] (340,80) -- (420,140) -- (340,140) -- cycle ;
%Shape: Right Triangle [id:dp9376691451708083] 
\draw  [fill={rgb, 255:red, 200; green, 200; blue, 200 }  ,fill opacity=1 ] (340,140) -- (420,200) -- (340,200) -- cycle ;
%Shape: Right Triangle [id:dp7014153666605438] 
\draw  [fill={rgb, 255:red, 200; green, 200; blue, 200 }  ,fill opacity=1 ] (260,80) -- (340,140) -- (260,140) -- cycle ;
%Straight Lines [id:da47442411197749135] 
\draw    (260,140) -- (340,140) ;
\draw [shift={(340,140)}, rotate = 0] [color={rgb, 255:red, 0; green, 0; blue, 0 }  ][fill={rgb, 255:red, 0; green, 0; blue, 0 }  ][line width=0.75]      (0, 0) circle [x radius= 3.35, y radius= 3.35]   ;

% Text Node
\draw (319,143) node [anchor=north west][inner sep=0.75pt]   [align=left] {$\displaystyle p$};

\end{tikzpicture}

%% file: BEPV_Figures_arXiv/staircase.tex
\tikzset{every picture/.style={line width=0.75pt}} %set default line width to 0.75pt        

\begin{tikzpicture}[x=0.75pt,y=0.75pt,yscale=-1,xscale=1]
%uncomment if require: \path (0,300); %set diagram left start at 0, and has height of 300

%Shape: Cube [id:dp5709458123224793] 
\draw   (200.5,180.5) -- (221,160) -- (279.5,160) -- (279.5,218.5) -- (259,239) -- (200.5,239) -- cycle ; \draw   (279.5,160) -- (259,180.5) -- (200.5,180.5) ; \draw   (259,180.5) -- (259,239) ;
%Shape: Cube [id:dp48099429692556184] 
\draw   (221,101.5) -- (241.5,81) -- (300,81) -- (300,139.5) -- (279.5,160) -- (221,160) -- cycle ; \draw   (300,81) -- (279.5,101.5) -- (221,101.5) ; \draw   (279.5,101.5) -- (279.5,160) ;
%Shape: Cube [id:dp9037651998320962] 
\draw   (279.5,160) -- (300,139.5) -- (358.5,139.5) -- (358.5,198) -- (338,218.5) -- (279.5,218.5) -- cycle ; \draw   (358.5,139.5) -- (338,160) -- (279.5,160) ; \draw   (338,160) -- (338,218.5) ;
%Shape: Cube [id:dp3010758993629642] 
\draw  [dash pattern={on 0.84pt off 2.51pt}] (221,160) -- (241.5,139.5) -- (300,139.5) -- (300,198) -- (279.5,218.5) -- (221,218.5) -- cycle ; \draw  [dash pattern={on 0.84pt off 2.51pt}] (300,139.5) -- (279.5,160) -- (221,160) ; \draw  [dash pattern={on 0.84pt off 2.51pt}] (279.5,160) -- (279.5,218.5) ;
%Straight Lines [id:da48319775288098554] 
\draw    (279.5,160) -- (221,160) ;
\draw [shift={(279.5,160)}, rotate = 180] [color={rgb, 255:red, 0; green, 0; blue, 0 }  ][fill={rgb, 255:red, 0; green, 0; blue, 0 }  ][line width=0.75]      (0, 0) circle [x radius= 3.35, y radius= 3.35]   ;

% Text Node
\draw (263,145) node [anchor=north west][inner sep=0.75pt]   [align=left] {$\displaystyle p$};

\end{tikzpicture}

%% file: BEPV_Figures_arXiv/pyramidbefore.tex
\tikzset{every picture/.style={line width=0.75pt}} %set default line width to 0.75pt        

\begin{tikzpicture}[x=0.75pt,y=0.75pt,yscale=-1,xscale=1]
%uncomment if require: \path (0,300); %set diagram left start at 0, and has height of 300

%Straight Lines [id:da2178887142234629] 
\draw    (291.5,67) -- (239,207) ;
\draw [shift={(291.5,67)}, rotate = 110.56] [color={rgb, 255:red, 0; green, 0; blue, 0 }  ][fill={rgb, 255:red, 0; green, 0; blue, 0 }  ][line width=0.75]      (0, 0) circle [x radius= 3.35, y radius= 3.35]   ;
%Straight Lines [id:da5739169124698561] 
\draw    (291.5,67) -- (326.5,207) ;
%Straight Lines [id:da3704798278969309] 
\draw    (291.5,67) -- (379,154.5) ;
%Straight Lines [id:da7724454748306906] 
\draw  [dash pattern={on 4.5pt off 4.5pt}]  (291.5,67) -- (291.5,154.5) ;
\draw [shift={(291.5,154.5)}, rotate = 90] [color={rgb, 255:red, 0; green, 0; blue, 0 }  ][fill={rgb, 255:red, 0; green, 0; blue, 0 }  ][line width=0.75]      (0, 0) circle [x radius= 3.35, y radius= 3.35]   ;
\draw [shift={(291.5,67)}, rotate = 90] [color={rgb, 255:red, 0; green, 0; blue, 0 }  ][fill={rgb, 255:red, 0; green, 0; blue, 0 }  ][line width=0.75]      (0, 0) circle [x radius= 3.35, y radius= 3.35]   ;
%Straight Lines [id:da7984164748378139] 
\draw  [dash pattern={on 4.5pt off 4.5pt}]  (291.5,154.5) -- (239,207) ;
%Straight Lines [id:da6146311748954177] 
\draw  [dash pattern={on 4.5pt off 4.5pt}]  (291.5,154.5) -- (379,154.5) ;
%Straight Lines [id:da6984154011553895] 
\draw    (326.5,207) -- (239,207) ;
%Straight Lines [id:da08475251298432196] 
\draw    (379,154.5) -- (326.5,207) ;
%Straight Lines [id:da3599593668450074] 
\draw [color={rgb, 255:red, 208; green, 2; blue, 27 }  ,draw opacity=0.5 ][line width=6]    (289.1,74.2) -- (370.2,155.8) ;
%Straight Lines [id:da996876519524639] 
\draw [color={rgb, 255:red, 208; green, 2; blue, 27 }  ,draw opacity=0.5 ][line width=6]    (294.2,80.2) -- (249.4,195.8) ;
%Straight Lines [id:da39728036171125003] 
\draw [color={rgb, 255:red, 74; green, 144; blue, 226 }  ,draw opacity=0.5 ][line width=4.5]    (291.5,67) -- (239,207) ;
%Straight Lines [id:da46695522935547573] 
\draw [color={rgb, 255:red, 74; green, 144; blue, 226 }  ,draw opacity=0.5 ][line width=4.5]    (291.5,67) -- (379,154.5) ;
%Straight Lines [id:da370193731657235] 
\draw [color={rgb, 255:red, 74; green, 144; blue, 226 }  ,draw opacity=0.5 ][line width=5.25]    (239,207) -- (326.5,207) ;
%Straight Lines [id:da4887607353148631] 
\draw [color={rgb, 255:red, 74; green, 144; blue, 226 }  ,draw opacity=0.5 ][line width=5.25]    (326.5,207) -- (379,154.5) ;
%Straight Lines [id:da19439257197716042] 
\draw [color={rgb, 255:red, 74; green, 144; blue, 226 }  ,draw opacity=0.5 ][line width=5.25]    (291.5,67) -- (326.5,207) ;
%Straight Lines [id:da838494916952799] 
\draw [color={rgb, 255:red, 208; green, 2; blue, 27 }  ,draw opacity=0.5 ][line width=4.5]    (291.5,154.5) -- (379,154.5) ;
%Straight Lines [id:da741225505892283] 
\draw [color={rgb, 255:red, 208; green, 2; blue, 27 }  ,draw opacity=0.5 ][line width=4.5]    (239,207) -- (291.5,154.5) ;
%Straight Lines [id:da8974311620752046] 
\draw [color={rgb, 255:red, 208; green, 2; blue, 27 }  ,draw opacity=0.5 ][line width=4.5]    (291.8,82.2) -- (291.5,154.5) ;
%Shape: Parallelogram [id:dp6426060397236342] 
\draw  [color={rgb, 255:red, 248; green, 231; blue, 28 }  ,draw opacity=0.5 ][line width=5.25]  (294.2,160.6) -- (362.2,160.6) -- (323.8,199.8) -- (255.8,199.8) -- cycle ;

% Text Node
\draw (239,129) node [anchor=north west][inner sep=0.75pt]   [align=left] {$\displaystyle F_{1}$};
% Text Node
\draw (344,90) node [anchor=north west][inner sep=0.75pt]   [align=left] {$\displaystyle F_{2}$};
% Text Node
\draw (280,173) node [anchor=north west][inner sep=0.75pt]   [align=left] {$\displaystyle F_{3}$};
% Text Node
\draw (328.25,160.5) node [anchor=north west][inner sep=0.75pt]   [align=left] {$\displaystyle F_{4}$};
% Text Node
\draw (298.8,48.2) node [anchor=north west][inner sep=0.75pt]   [align=left] {$\displaystyle p$};
% Text Node
\draw (276.8,134.8) node [anchor=north west][inner sep=0.75pt]   [align=left] {$\displaystyle q$};

\end{tikzpicture}

 

%% file: BEPV_Figures_arXiv/pyramidafter.tex
\tikzset{every picture/.style={line width=0.75pt}} %set default line width to 0.75pt        

\begin{tikzpicture}[x=0.75pt,y=0.75pt,yscale=-1,xscale=1]
%uncomment if require: \path (0,300); %set diagram left start at 0, and has height of 300

%Straight Lines [id:da15410407277173932] 
\draw    (281,90) -- (261,190) ;
%Straight Lines [id:da9361121555665884] 
\draw    (321,90) -- (381,170) ;
%Straight Lines [id:da9975109665640945] 
\draw  [dash pattern={on 4.5pt off 4.5pt}]  (321,150) -- (381,170) ;
%Straight Lines [id:da15332774967071638] 
\draw    (341,210) -- (261,190) ;
%Straight Lines [id:da6208063155987579] 
\draw    (381,170) -- (341,210) ;
%Straight Lines [id:da330299078792216] 
\draw  [dash pattern={on 4.5pt off 4.5pt}]  (281,170) -- (261,190) ;
%Straight Lines [id:da6668820256988214] 
\draw  [dash pattern={on 4.5pt off 4.5pt}]  (321,150) -- (281,170) ;
%Straight Lines [id:da8392964013058353] 
\draw  [dash pattern={on 4.5pt off 4.5pt}]  (281,90) -- (281,170) ;
%Straight Lines [id:da2666645473358342] 
\draw  [dash pattern={on 4.5pt off 4.5pt}]  (321,90) -- (321,150) ;
%Straight Lines [id:da9556971667559309] 
\draw    (321,90) -- (281,90) ;
%Curve Lines [id:da36183889747842224] 
\draw    (301,90) .. controls (302,140) and (313,170) .. (341,210) ;
\draw [shift={(301,90)}, rotate = 88.85] [color={rgb, 255:red, 0; green, 0; blue, 0 }  ][fill={rgb, 255:red, 0; green, 0; blue, 0 }  ][line width=0.75]      (0, 0) circle [x radius= 3.35, y radius= 3.35]   ;
%Straight Lines [id:da357309798059953] 
\draw [color={rgb, 255:red, 208; green, 2; blue, 27 }  ,draw opacity=0.5 ][line width=5.25]    (281,93.67) -- (281,170) ;
%Straight Lines [id:da9424375349697551] 
\draw [color={rgb, 255:red, 208; green, 2; blue, 27 }  ,draw opacity=0.5 ][line width=5.25]    (321,90) -- (321,150) ;
%Straight Lines [id:da4909252626353522] 
\draw [color={rgb, 255:red, 74; green, 144; blue, 226 }  ,draw opacity=0.5 ][line width=5.25]    (321,90) -- (281,90) ;
%Straight Lines [id:da5811295010369295] 
\draw [color={rgb, 255:red, 248; green, 231; blue, 28 }  ,draw opacity=0.5 ][line width=5.25]    (281,170) -- (321,150) ;

% Text Node
\draw (241,130) node [anchor=north west][inner sep=0.75pt]   [align=left] {$\displaystyle \tilde{F}_{1}$};
% Text Node
\draw (288,130) node [anchor=north west][inner sep=0.75pt]   [align=left] {$\displaystyle E$};
% Text Node
\draw (357,100) node [anchor=north west][inner sep=0.75pt]   [align=left] {$\displaystyle \tilde{F}_{2}$};
% Text Node
\draw (340,165) node [anchor=north west][inner sep=0.75pt]   [align=left] {$\displaystyle \tilde{F}_{4}$};
% Text Node
\draw (283,205) node [anchor=north west][inner sep=0.75pt]   [align=left] {$\displaystyle \tilde{F}_{3}$};
% Text Node
\draw (288,62) node [anchor=north west][inner sep=0.75pt]   [align=left] {$\displaystyle ( p,c)$};

\end{tikzpicture}

%% file: BEPV_Figures_arXiv/segmentno.tex
\tikzset{every picture/.style={line width=0.75pt}} %set default line width to 0.75pt        

\begin{tikzpicture}[x=0.75pt,y=0.75pt,yscale=-1,xscale=1]
%uncomment if require: \path (0,300); %set diagram left start at 0, and has height of 300

%Straight Lines [id:da4114231298523543] 
\draw    (298,91) -- (278,191) ;
%Straight Lines [id:da6564986215997618] 
\draw    (338,91) -- (398,171) ;
%Straight Lines [id:da5645329729778803] 
\draw  [dash pattern={on 4.5pt off 4.5pt}]  (338,151) -- (398,171) ;
%Straight Lines [id:da4360756972508866] 
\draw    (358,211) -- (278,191) ;
%Straight Lines [id:da693897444657799] 
\draw    (398,171) -- (358,211) ;
%Straight Lines [id:da3067313995661176] 
\draw  [dash pattern={on 4.5pt off 4.5pt}]  (298,171) -- (278,191) ;
%Straight Lines [id:da4528551743088576] 
\draw  [dash pattern={on 4.5pt off 4.5pt}]  (338,151) -- (298,171) ;
%Straight Lines [id:da7415531217076259] 
\draw  [dash pattern={on 4.5pt off 4.5pt}]  (298,91) -- (298,171) ;
%Straight Lines [id:da7656194656443613] 
\draw  [dash pattern={on 4.5pt off 4.5pt}]  (338,91) -- (338,151) ;
%Straight Lines [id:da5325961638875116] 
\draw    (338,91) -- (298,91) ;
%Curve Lines [id:da37664662519867276] 
\draw    (318,91) .. controls (319,141) and (330,171) .. (358,211) ;
\draw [shift={(318,91)}, rotate = 88.85] [color={rgb, 255:red, 0; green, 0; blue, 0 }  ][fill={rgb, 255:red, 0; green, 0; blue, 0 }  ][line width=0.75]      (0, 0) circle [x radius= 3.35, y radius= 3.35]   ;
%Straight Lines [id:da3810666879513519] 
\draw [draw opacity=0][fill={rgb, 255:red, 126; green, 211; blue, 33 }  ,fill opacity=0.5 ][line width=1.5]    (358,210.4) -- (398,170.4) ;
%Straight Lines [id:da16835848894022354] 
\draw [draw opacity=0][fill={rgb, 255:red, 126; green, 211; blue, 33 }  ,fill opacity=0.5 ][line width=1.5]    (338,90.4) -- (398,170.4) ;
%Straight Lines [id:da8511331499305845] 
\draw [draw opacity=0][fill={rgb, 255:red, 126; green, 211; blue, 33 }  ,fill opacity=0.5 ][line width=1.5]    (338,90.4) -- (318,90.4) ;
%Curve Lines [id:da052534401715136014] 
\draw [draw opacity=0][fill={rgb, 255:red, 126; green, 211; blue, 33 }  ,fill opacity=0.5 ][line width=1.5]    (318,90.4) .. controls (319,140.4) and (330,170.4) .. (358,210.4) ;

%Shape: Polygon [id:ds8268037276090691] 
\draw  [draw opacity=0][fill={rgb, 255:red, 126; green, 211; blue, 33 }  ,fill opacity=0.5 ] (338,90.4) -- (398,170.4) -- (358,210.4) -- (318,90.4) -- cycle ;

% Text Node
\draw (258,131) node [anchor=north west][inner sep=0.75pt]   [align=left] {$\displaystyle \tilde{F}_{1}$};
% Text Node
\draw (305,131) node [anchor=north west][inner sep=0.75pt]   [align=left] {$\displaystyle E$};
% Text Node
\draw (374,101) node [anchor=north west][inner sep=0.75pt]   [align=left] {$\displaystyle \tilde{F}_{2}$};
% Text Node
\draw (357,170) node [anchor=north west][inner sep=0.75pt]   [align=left] {$\displaystyle \tilde{F}_{4}$};
% Text Node
\draw (300,206) node [anchor=north west][inner sep=0.75pt]   [align=left] {$\displaystyle \tilde{F}_{3}$};
% Text Node
\draw (305,63) node [anchor=north west][inner sep=0.75pt]   [align=left] {$\displaystyle ( p,c)$};

\end{tikzpicture}

%% file: BEPV_Figures_arXiv/segmentyes.tex
\tikzset{every picture/.style={line width=0.75pt}} %set default line width to 0.75pt        

\begin{tikzpicture}[x=0.75pt,y=0.75pt,yscale=-1,xscale=1]
%uncomment if require: \path (0,300); %set diagram left start at 0, and has height of 300

%Straight Lines [id:da3391760191869235] 
\draw    (318,111) -- (298,211) ;
%Straight Lines [id:da8873689448106391] 
\draw    (358,111) -- (418,191) ;
%Straight Lines [id:da21171960907717913] 
\draw  [dash pattern={on 4.5pt off 4.5pt}]  (358,171) -- (418,191) ;
%Straight Lines [id:da7782626834227745] 
\draw    (378,231) -- (298,211) ;
%Straight Lines [id:da2510811295977903] 
\draw    (418,191) -- (378,231) ;
%Straight Lines [id:da8869561392198967] 
\draw  [dash pattern={on 4.5pt off 4.5pt}]  (318,191) -- (298,211) ;
%Straight Lines [id:da43758722498994296] 
\draw  [dash pattern={on 4.5pt off 4.5pt}]  (358,171) -- (318,191) ;
%Straight Lines [id:da9459800322390518] 
\draw  [dash pattern={on 4.5pt off 4.5pt}]  (318,111) -- (318,191) ;
%Straight Lines [id:da7975384267012031] 
\draw  [dash pattern={on 4.5pt off 4.5pt}]  (358,111) -- (358,171) ;
%Straight Lines [id:da6460910211792301] 
\draw    (358,111) -- (318,111) ;
%Curve Lines [id:da705536206121621] 
\draw [fill={rgb, 255:red, 189; green, 16; blue, 224 }  ,fill opacity=0.5 ]   (338.11,110.89) .. controls (339.11,160.89) and (350.11,190.89) .. (378.11,230.89) ;
\draw [shift={(338.11,110.89)}, rotate = 88.85] [color={rgb, 255:red, 0; green, 0; blue, 0 }  ][fill={rgb, 255:red, 0; green, 0; blue, 0 }  ][line width=0.75]      (0, 0) circle [x radius= 3.35, y radius= 3.35]   ;
%Straight Lines [id:da5720042204141664] 
\draw [color={rgb, 255:red, 189; green, 16; blue, 224 }  ,draw opacity=0.5 ][line width=5.25]    (358,111) -- (318,111) ;
%Shape: Polygon [id:ds2769649433268119] 
\draw  [draw opacity=0][fill={rgb, 255:red, 189; green, 16; blue, 224 }  ,fill opacity=0.5 ] (358.11,110.89) -- (418.11,190.89) -- (378.11,230.89) -- (338.11,110.89) -- cycle ;

% Text Node
\draw (278,151) node [anchor=north west][inner sep=0.75pt]   [align=left] {$\displaystyle \tilde{F}_{1}$};
% Text Node
\draw (325,151) node [anchor=north west][inner sep=0.75pt]   [align=left] {$\displaystyle E$};
% Text Node
\draw (394,121) node [anchor=north west][inner sep=0.75pt]   [align=left] {$\displaystyle \tilde{F}_{2}$};
% Text Node
\draw (377,190) node [anchor=north west][inner sep=0.75pt]   [align=left] {$\displaystyle \tilde{F}_{4}$};
% Text Node
\draw (320,226) node [anchor=north west][inner sep=0.75pt]   [align=left] {$\displaystyle \tilde{F}_{3}$};
% Text Node
\draw (325,83) node [anchor=north west][inner sep=0.75pt]   [align=left] {$\displaystyle ( p,c)$};

\end{tikzpicture}

%% file: BEPV_Figures_arXiv/schlegel.tex
\tikzset{every picture/.style={line width=0.75pt}} %set default line width to 0.75pt        

\begin{tikzpicture}[x=0.75pt,y=0.75pt,yscale=-1,xscale=1]
%uncomment if require: \path (0,300); %set diagram left start at 0, and has height of 300

%Straight Lines [id:da5684582856664965] 
\draw    (240,220) -- (360,220) ;
\draw [shift={(360,220)}, rotate = 0] [color={rgb, 255:red, 0; green, 0; blue, 0 }  ][fill={rgb, 255:red, 0; green, 0; blue, 0 }  ][line width=0.75]      (0, 0) circle [x radius= 3.35, y radius= 3.35]   ;
%Straight Lines [id:da8788742906501231] 
\draw    (360,220) -- (420,160) ;
%Straight Lines [id:da3336505328012628] 
\draw  [dash pattern={on 0.84pt off 2.51pt}]  (240,220) -- (300,160) ;
%Straight Lines [id:da3138308975118722] 
\draw  [dash pattern={on 0.84pt off 2.51pt}]  (300,160) -- (420,160) ;
%Straight Lines [id:da1764294389070561] 
\draw  [dash pattern={on 0.84pt off 2.51pt}]  (300,160) -- (360,60) ;
\draw [shift={(300,160)}, rotate = 300.96] [color={rgb, 255:red, 0; green, 0; blue, 0 }  ][fill={rgb, 255:red, 0; green, 0; blue, 0 }  ][line width=0.75]      (0, 0) circle [x radius= 3.35, y radius= 3.35]   ;
%Straight Lines [id:da26451583455901306] 
\draw    (240,220) -- (300,120) ;
\draw [shift={(300,120)}, rotate = 300.96] [color={rgb, 255:red, 0; green, 0; blue, 0 }  ][fill={rgb, 255:red, 0; green, 0; blue, 0 }  ][line width=0.75]      (0, 0) circle [x radius= 3.35, y radius= 3.35]   ;
\draw [shift={(240,220)}, rotate = 300.96] [color={rgb, 255:red, 0; green, 0; blue, 0 }  ][fill={rgb, 255:red, 0; green, 0; blue, 0 }  ][line width=0.75]      (0, 0) circle [x radius= 3.35, y radius= 3.35]   ;
%Straight Lines [id:da0965278818552433] 
\draw    (360,60) -- (420,160) ;
\draw [shift={(420,160)}, rotate = 59.04] [color={rgb, 255:red, 0; green, 0; blue, 0 }  ][fill={rgb, 255:red, 0; green, 0; blue, 0 }  ][line width=0.75]      (0, 0) circle [x radius= 3.35, y radius= 3.35]   ;
%Straight Lines [id:da3023046298656691] 
\draw    (300,120) -- (360,60) ;
\draw [shift={(360,60)}, rotate = 315] [color={rgb, 255:red, 0; green, 0; blue, 0 }  ][fill={rgb, 255:red, 0; green, 0; blue, 0 }  ][line width=0.75]      (0, 0) circle [x radius= 3.35, y radius= 3.35]   ;
%Straight Lines [id:da580044035676378] 
\draw  [dash pattern={on 0.84pt off 2.51pt}]  (340,140) -- (240,220) ;
\draw [shift={(240,220)}, rotate = 141.34] [color={rgb, 255:red, 0; green, 0; blue, 0 }  ][fill={rgb, 255:red, 0; green, 0; blue, 0 }  ][line width=0.75]      (0, 0) circle [x radius= 3.35, y radius= 3.35]   ;
\draw [shift={(340,140)}, rotate = 141.34] [color={rgb, 255:red, 0; green, 0; blue, 0 }  ][fill={rgb, 255:red, 0; green, 0; blue, 0 }  ][line width=0.75]      (0, 0) circle [x radius= 3.35, y radius= 3.35]   ;
%Straight Lines [id:da20085460619217532] 
\draw    (420,160) -- (340,140) ;
\draw [shift={(340,140)}, rotate = 194.04] [color={rgb, 255:red, 0; green, 0; blue, 0 }  ][fill={rgb, 255:red, 0; green, 0; blue, 0 }  ][line width=0.75]      (0, 0) circle [x radius= 3.35, y radius= 3.35]   ;
\draw [shift={(420,160)}, rotate = 194.04] [color={rgb, 255:red, 0; green, 0; blue, 0 }  ][fill={rgb, 255:red, 0; green, 0; blue, 0 }  ][line width=0.75]      (0, 0) circle [x radius= 3.35, y radius= 3.35]   ;
%Straight Lines [id:da7720217014455926] 
\draw    (360,60) -- (340,140) ;
\draw [shift={(340,140)}, rotate = 104.04] [color={rgb, 255:red, 0; green, 0; blue, 0 }  ][fill={rgb, 255:red, 0; green, 0; blue, 0 }  ][line width=0.75]      (0, 0) circle [x radius= 3.35, y radius= 3.35]   ;
\draw [shift={(360,60)}, rotate = 104.04] [color={rgb, 255:red, 0; green, 0; blue, 0 }  ][fill={rgb, 255:red, 0; green, 0; blue, 0 }  ][line width=0.75]      (0, 0) circle [x radius= 3.35, y radius= 3.35]   ;
%Straight Lines [id:da6979985723538971] 
\draw    (340,140.5) -- (360,220) ;
\draw [shift={(360,220)}, rotate = 75.88] [color={rgb, 255:red, 0; green, 0; blue, 0 }  ][fill={rgb, 255:red, 0; green, 0; blue, 0 }  ][line width=0.75]      (0, 0) circle [x radius= 3.35, y radius= 3.35]   ;
\draw [shift={(340,140.5)}, rotate = 75.88] [color={rgb, 255:red, 0; green, 0; blue, 0 }  ][fill={rgb, 255:red, 0; green, 0; blue, 0 }  ][line width=0.75]      (0, 0) circle [x radius= 3.35, y radius= 3.35]   ;
%Straight Lines [id:da31845149506684634] 
\draw  [dash pattern={on 0.84pt off 2.51pt}]  (340,140.5) -- (300,160) ;
\draw [shift={(300,160)}, rotate = 154.01] [color={rgb, 255:red, 0; green, 0; blue, 0 }  ][fill={rgb, 255:red, 0; green, 0; blue, 0 }  ][line width=0.75]      (0, 0) circle [x radius= 3.35, y radius= 3.35]   ;
\draw [shift={(340,140.5)}, rotate = 154.01] [color={rgb, 255:red, 0; green, 0; blue, 0 }  ][fill={rgb, 255:red, 0; green, 0; blue, 0 }  ][line width=0.75]      (0, 0) circle [x radius= 3.35, y radius= 3.35]   ;
%Straight Lines [id:da5789665861721902] 
\draw    (300,120) -- (360,220) ;
%Straight Lines [id:da8064666050068923] 
\draw    (340,140.5) -- (300,120) ;
\draw [shift={(300,120)}, rotate = 207.14] [color={rgb, 255:red, 0; green, 0; blue, 0 }  ][fill={rgb, 255:red, 0; green, 0; blue, 0 }  ][line width=0.75]      (0, 0) circle [x radius= 3.35, y radius= 3.35]   ;
\draw [shift={(340,140.5)}, rotate = 207.14] [color={rgb, 255:red, 0; green, 0; blue, 0 }  ][fill={rgb, 255:red, 0; green, 0; blue, 0 }  ][line width=0.75]      (0, 0) circle [x radius= 3.35, y radius= 3.35]   ;

% Text Node
\draw (221,220) node [anchor=north west][inner sep=0.75pt]   [align=left] {$\displaystyle 1$};
% Text Node
\draw (286,146) node [anchor=north west][inner sep=0.75pt]   [align=left] {$\displaystyle 0$};
% Text Node
\draw (362,223) node [anchor=north west][inner sep=0.75pt]   [align=left] {$\displaystyle 3$};
% Text Node
\draw (428,142) node [anchor=north west][inner sep=0.75pt]   [align=left] {$\displaystyle 2$};
% Text Node
\draw (368,42) node [anchor=north west][inner sep=0.75pt]   [align=left] {$\displaystyle 5$};
% Text Node
\draw (352,121) node [anchor=north west][inner sep=0.75pt]   [align=left] {$\displaystyle 4$};
% Text Node
\draw (281,102) node [anchor=north west][inner sep=0.75pt]   [align=left] {$\displaystyle 6$};

\end{tikzpicture}